\documentclass[10pt,reqno]{amsart}
\usepackage{lipsum}
\usepackage{amsmath}
\usepackage{mathtools}
\usepackage{enumerate}
\usepackage{amssymb}
\usepackage{hyperref}
\usepackage[all]{xy}
\usepackage{tensor}
\usepackage{comment}
\usepackage{etoolbox}
\usepackage{tikz}
\usetikzlibrary{shapes.geometric}
\usetikzlibrary{cd}
\usetikzlibrary{calc}
\tikzcdset{every label/.append style = {font = \small}}

\newcommand{\C}{\mathcal{C}}

\newcommand{\I}{\mathcal{I}}

\newcommand{\K}{\mathcal{K}}

\newcommand{\Tc}{\mathcal{T}}
\newcommand{\Sc}{\mathcal{S}}

\newcommand{\e}{\epsilon}

\newcommand{\RC}{\mathcal{RC}}

\newcommand{\RTC}{\mathcal{RTC}}
\newcommand{\TA}{\mathcal{TA}}
\newcommand{\TC}{\mathcal{TC}}

\DeclareMathOperator{\CC}{\C \boxtimes \overline{\C}}
\DeclareMathOperator{\ann}{an}
\DeclareMathOperator{\cre}{cr}

\DeclareMathOperator{\Obj}{Obj}
\DeclareMathOperator{\Irr}{Irr}
\DeclareMathOperator{\Blank}{--}
\DeclareMathOperator{\End}{End}
\DeclareMathOperator{\Mod}{Mod}
\DeclareMathOperator{\Hom}{Hom}

\DeclareMathOperator{\Vect}{\normalfont\underline{Vect}}
\DeclareMathOperator{\fld}{\mathbb{K}}
\DeclareMathOperator{\Yon}{\yen}

\DeclareMathOperator{\id}{id}

\DeclareMathOperator{\tid}{\normalfont \textbf{1}}

\DeclareMathOperator{\vprod}{\cdot}

\newcommand{\tl}[1][\beta]{\text{\normalfont TL}}
\newcommand{\Rtl}[1][\beta]{\mathcal{R}\text{\normalfont TL}}
\newcommand{\tlred}[1][\beta]{\text{\normalfont TL}^\text{red}}
\newcommand{\Rtlred}[1][\beta]{\mathcal{R}\text{\normalfont TL}^\text{red}}

\newcommand{\dual}{^*}

\newcommand{\Sharp}{^\sharp}

\newcommand{\idem}{\varepsilon}
\newcommand{\defeq}{\vcentcolon=}
\newcommand{\LMod}[1]{\Mod \! \mbox{-}#1}

\theoremstyle{plain}
\newtheorem{THM}{Theorem}[section]
\newtheorem{PROP}[THM]{Proposition}
\newtheorem{LEMMA}[THM]{Lemma}

\newtheorem{COR}[THM]{Corollary}
\theoremstyle{definition}
\newtheorem{REM}[THM]{Remark}

\newtheorem{DEF}[THM]{Definition}
\newtheorem*{DEF*}{Definition}

\title{A Graphical Approach to the Drinfeld Centre}

\usepackage{amsaddr}
\author{Leonard Hardiman}
\address{Institut Camille Jordan\\
Universit\'e Claude Bernard Lyon 1\\
43 Boulevard du 11 novembre 1918\\
69622 Villeurbanne Cedex France\\
ORCID: 0000-0003-1986-6704}
\email{hardiman@math.univ-lyon1.fr}

\begin{document}

    \begin{abstract}

    Let $ \C $ be a spherical fusion category. The goal of this  article is to present the tube category of $ \C $, denoted $ \TC $, as giving a graphical perspective on the Drinfeld centre of $ \C $, denoted $ Z(\C) $. We then exploit this perspective to obtain an alternative proof of the equivalence between $ Z(\C) $ and $ \CC $ when $ \C $ is modular. Sections \ref{sec:prelim_lin_cat}-\ref{sec:prelim_fusion_cat} provide a survey of the required prerequisites; readers already familiar with graphical calculus in fusion categories should start at Section \ref{sec:tube_category}.

    \end{abstract}

\maketitle

\section{Introduction}

The graphical calculus of monoidal categories has played an increasingly important role in mathematical physics since its original conception by Pensrose~\cite{MR0281657} and subsequent formalisation due to Joyal and Street~\cite{MR1113284}. The Drinfeld centre of a monoidal category also occupies a prominent place within current research. In this article we consider the question of how to apply the tools offered by graphical calculus to computations within the Drinfeld centre in the case when the underlying category is a spherical fusion category.

Let $ \C $ be a spherical fusion category. The idea that graphical calculus in the Drinfeld centre of $ \C $, denoted $ Z(\C) $, should correspond to graphical calculus in $ \C $ drawn on a cylinder (or tube) goes back to Ocneanu's definition of the tube algebra~\cite{Ocneanu1993}. Although originally defined in an operator algebra context, the concept generalises to an arbitrary spherical fusion category and a result of Popa, Shlyakhtenko and Vaes \cite[Proposition 3.14]{MR3801484} proves that the category of representations of the tube algebra is equivalent to $ Z(\C) $. Independently, Kirillov Jr.\ worked on the string-net spaces introduced in the papers of Levin and Wen~\cite{MR2717368}. These spaces may be thought of as the $ \Hom $-spaces in categories whose graphical calculus is drawn on an arbitrary surface. In particular,~\cite[Theorem 6.4]{alex2011stringnet} proves that, when the surface is a cylinder, one recovers a category whose idempotent completion is $ Z(\C) $.

The purpose of this article is to provide a `from the ground up' exposition of this perspective by analysing the so called \emph{tube category}, first introduced in~\cite{hardiman_king}. This construction may be viewed as a middle ground between the tube algebra and string-net approaches. Although the resulting object is a category rather than an algebra, the definition more closely resembles that of the tube algebra. In particular, considering the endomorphism algebra of a generating object in the tube category directly recovers Ocneanu's tube algebra. Working with a category, however, allows us to describe the idempotents which appear in Kirillov's work more easily.

The first four sections provide a gentle introduction to all the required pre-requisites. The tube category is then defined and it is shown that its category of representations (or equivalently, its idempotent completion) is equivalent to $ Z(\C) $. This equivalence stems from the fact that, when $ \C $ is a spherical fusion category, the Yoneda embedding is an equivalence and the data required to extend the image of $ X $ under the Yoneda embedding to $ \TC $ corresponds to a half braiding on $ X $. To present an alternative perspective on this equivalence, for an object $ (X,\tau) $ in $ Z(\C) $ we describe an idempotent
	\begin{align*}
	\e_{\tau} = \frac{1}{d(\C)}
        \bigoplus\limits_S d(S)
		\begin{array}{c}
			\begin{tikzpicture}[scale=0.25,every node/.style={inner sep=0,outer sep=-1}]
			\node (v1) at (0,4) {};
			\node (v4) at (0,-4) {};
			\node (v2) at (4,0) {};
			\node (v3) at (-4,0) {};
			\node (v5) at (-2,2) {};
			\node (v6) at (2,-2) {};
			\node (v7) at (2,2) {};
			\node (v11) at (-2,-2) {};
			\node [draw,diamond,outer sep=0,inner sep=.5,minimum size=28,fill=white] (v9) at (0,0) {$ \tau_{S} $};
			\node at (3,3) {$ X $};
			\node at (-3,-3) {$ X $};
			\node at (-3,3) {$ S $};
			\node at (3,-3) {$ S $};
			\draw [thick] (v9) edge (v6);
			\draw [thick] (v9) edge (v7);
			\draw [thick] (v9) edge (v11);
			\draw [thick] (v5) edge (v9);
			\draw[very thick, red]  (v1) edge (v3);
			\draw[very thick, red]  (v2) edge (v4);
			\draw[very thick]  (v1) edge (v2);
			\draw[very thick]  (v3) edge (v4);
			\end{tikzpicture}
		\end{array}
    \in \End_{\TC}(X)
	\end{align*}
which yields the corresponding representation under the Yoneda embedding. Finally, the $ \TC $ framework is exploited to provide an alternative proof of the equivalence between $ Z(\C) $ and $ \CC $ when $ \C $ is a modular tensor category (this result is originally due to M\"uger~\cite{MR1990929}). This is achieved by applying graphical calculus to the idempotent
	\begin{align*}
	\e_X^Y = \frac{1}{d(\C)} \bigoplus\limits_S d(S)
		\begin{array}{c}
			\begin{tikzpicture}[scale=0.2,every node/.style={inner sep=0,outer sep=-1}]
			\node (v1) at (0,5) {};
			\node (v4) at (0,-5) {};
			\node (v2) at (5,0) {};
			\node (v3) at (-5,0) {};
			\node (v5) at (-2.5,2.5) {};
			\node (v6) at (2.5,-2.5) {};
			\node (v7) at (1.5,3.5) {};
			\node (v11) at (3.5,1.5) {};
			\node (v12) at (-1.5,-3.5) {};
			\node (v10) at (-3.5,-1.5) {};
			\node (v13) at (2.5,-2.5) {};
			\draw [thick] (v7) edge (v10);
			\draw [line width =0.5em,white] (v5) edge (v13);
			\draw [thick] (v5) edge (v13);
			\draw [line width =0.5em,white] (v11) edge (v12);
			\draw [thick] (v11) edge (v12);
			\node at (2.75,4.75) {$ X $};
			\node at (4.75,2.75) {$ Y $};
			\node at (-3.75,3.75) {$ S $};
			\node at (3.75,-3.5) {$ S $};
			\draw[very thick, red]  (v1) edge (v3);
			\draw[very thick, red]  (v2) edge (v4);
			\draw[very thick]  (v1) edge (v2);
			\draw[very thick]  (v3) edge (v4);
			\end{tikzpicture}
		\end{array} \in \End_{TC}(XY)
	\end{align*}
which, under the Yoneda embedding, corresponds to the image of $ X \boxtimes Y $ under the canonical functor $ \Phi\colon \CC \to Z(\C) = \RTC $. In particular it is shown that the complete set of simples $ \{I \boxtimes J \}_{I,J \in \Irr(\C)} $ in $ \CC $ maps to a complete set of simples in $ \RTC = Z(\C) $.     

\vspace{0.8em}

\noindent \textbf{Acknowledgements.} The author thanks Alastair King for his guidance during the period this work was carried out. He is also grateful to Ingo Runkel for multiple helpful conversations.

\section{Preliminaries on $ \fld $-linear Categories} \label{sec:prelim_lin_cat}

We start by introducing the basics on linear categories. The material within this section is completely standard; our treatment is designed to render this article as self-contained as possible.

Let $ \fld $ be a field. From now on all categories are assumed to be \emph{linear categories over $ \fld $}, i.e.\ the $ \Hom $-spaces are finite dimensional vector spaces over $ \fld $ and composition is bilinear. For example, the category of finite dimensional vector spaces $ \Vect $ is a linear category as $ \Hom(V,W) $ is a $ \dim(V)\times\dim(W) $-dimensional vector space. Furthermore all functors are \emph{linear functors}, i.e.\ functors between two linear categories such that the corresponding maps between $ \Hom $-spaces are linear.

We say that an object $ Z $ is a \emph{direct sum} of two objects $ X $ and $ Y $ if $ \Hom_\C(Z,A) $ is naturally identified with $ \Hom_\C(X,A) \oplus \Hom_\C(Y,A) $ and $ \Hom_\C(A,Z) $ is naturally identified with $ \Hom_\C(A,X) \oplus \Hom_\C(A,Y) $. Similarly, we say that an object $ Z $ is the \emph{product} of a vector space $ V $ and and object $ X $ if $ \Hom_{\C}(Z, Y) $ is naturally identified with $ V\dual \otimes \Hom_{\C}(X, Y) $ and $ \Hom_{\C}(Y, Z) $ is naturally identified with $ V \otimes \Hom_{\C}(Y, X) $. Note that both the direct sum of $ X $ and $ Y $ and the product of $ V $ with $ X $ are unique; this can also be seen as a consequence of the Yoneda Lemma. They are denoted $ X \oplus Y $ and $ V \vprod X $ respectively.

In practice we never consider the question of whether or not direct sums or products exist. This is due to the fact that, if they do not exist, they may be formally added unambiguously. The following lemma captures the relationship between products and direct sums.

	\begin{LEMMA} \label{lem:basis_decomposition}

	Let be $ V $ in $ \Vect $, let $ {b} $ be a basis of $ V $ and let $ X $ be in $ \C $. Then $$ V \vprod X = \bigoplus\limits_b X. $$

	\proof

	Let $ {b^*} $ be the dual basis to $ {b} $. The maps $ \bigoplus_b b \otimes \id_X \in \Hom_{\C}\big(\bigoplus_b X, V \vprod X \big) $ and $ \bigoplus_b b^* \otimes \id_X \in \Hom_{\C}\big( V \vprod X, \bigoplus_b X \big) $ are inverse to one another.	\endproof

	\end{LEMMA}

We recall that an object $ X $ in $ \C $ is called \emph{simple} if $ X $ has no proper subobjects. Schur's Lemma implies that, for a simple object $ S $ in $ \C $, $ \End_{\C}(S) $ is a division algebra over $ \fld $. A \emph{complete set of simples} is a set $ \Irr(\C) $ such that for all $ I $ in $ \Irr(\C) $, $ I $ is a simple object in $ \C $ and for all simple object $ S $ in $ \C $ there exists a unique $ I \in \Irr(\C) $ such that $ \Hom_{\C}(S,I) \neq 0 $. We call an object $ X $ \emph{Schurian} if $ \End_{\C}(X) = \fld $. We call a category Schurian if all of its simple objects are Schurian. In particular if $ \fld $ is algebraically closed then $ \C $ is Schurian. Finally, we recall that a category $ \C $ is called \emph{semisimple} if every object in $ \C $ is a direct sum of finitely many simple objects. The following canonical decomposition of an object in a semisimple category will be of great importance for the remainder of this article.

    	\begin{PROP} \label{prop:simple_decomposition}

    	Let $ \C $ be a semisimple, Schurian category, let $ \Irr(\C) $ be a complete set of simples and let $ X $ be an object in $ \C $. Then
    		\begin{align*}
    		X = \bigoplus\limits_{S \in \Irr(\C)} \Hom_{\C}(S,X) \vprod S.
    		\end{align*}
    	\proof

    	As $ X $ is semisimple we have
    		\begin{align*}
    		X = \bigoplus\limits_{i \in \I} X_i
    		\end{align*}
    	where the $ X_i $ are simple objects and $ \I $ is an indexing set. We therefore have
    		\begin{align*}
    		\bigoplus\limits_{S \in \Irr(\C)} \Hom_{\C}(S,X) \vprod S = \bigoplus\limits_{\substack{i \in \I \\ S \in \Irr(C)}} \Hom_{\C}(S,X_i) \vprod S = \bigoplus\limits_{\substack{i \in \I_{S} \\ S \in \Irr(C)}} \Hom_{\C}(S,X_i) \vprod S
    		\end{align*}
    	where $ \I_S = \{i \in \I \mid X_i \cong S \} $. As $ \C $ is Schurian $ \Hom_{\C}(S,X_i) $ is one-dimensional and the canonical morphism
    		\begin{align*}
    		\id \in \End(\Hom_{\C}(S,X_i))
    	    &= \Hom_{\C}(S,X_i)\dual \otimes \Hom_{\C}(S,X_i) \\
    	    &= \Hom_{\C}(\Hom_{\C}(S,X_i) \vprod S, X_i)
    		\end{align*}
    	is an isomorphism by Schur's Lemma. Therefore
    		\begin{align*}
    		\bigoplus\limits_{\substack{i \in \I_{S} \\ S \in \Irr(C)}} \Hom_{\C}(S,X_i) \vprod S = \bigoplus\limits_{\substack{i \in \I_{S} \\ S \in \Irr(C)}} X_i = X.
    		\end{align*}
    	\endproof

    	\end{PROP}

Although our focus will mainly be on semisimple categories throughout this article, we shall come across certain cagoteries which are only semisimple `in spirit'. To illustrate this point we consider again $ \Vect $, the category of finite dimensional vector spaces. Let $ V $ be in $ \Vect $ and suppose $ \idem \in \End(V) $ is an idempotent, i.e. an endomorphism satisfying $ \idem^2 = \idem $. In this case $ V_\idem = \{ v \in V \mid \idem(v) = v \} $ forms a subspace of $ V $ and $ \idem $ may be thought of as a projection from $ V $ onto $ V_\idem $ composed with an inclusion of $ V_\idem $ back into $ V $. Furthermore this correspondence between idempotents on $ V $ and split subspaces of $ V $ defines a bijection. To say this more generally we made the following definition,

	\begin{DEF}

	Let $ \idem \in \End_{\C}(X) $ be an idempotent in a category $ \C $. An \emph{image object} for $ \idem $ is a object $ X_\idem $ in $ \C $ together with morphisms $ \pi \colon X \to X_\idem $ and  $ i \colon X_\idem \to X $ such that  $ i \circ \pi = \idem $ and $ \pi \circ i = \id_{X_\idem} $.

	\end{DEF}

    \begin{REM}

	Suppose $ (X_1, \pi_1, i_1) $ and $ (X_2, \pi_2, i_2) $ are two image objects for an idempotent $ \idem \in \End_{\C}(X) $. Then $ \pi_2 \circ i_1 $ and $ \pi_1 \circ i_2 $ give inverse isomorphisms between $ X_1 $ and $ X_2 $. Furthermore, as the diagram
	   \begin{center}
		    \begin{tikzcd}[row sep=huge]
			& X \arrow[shift left=0.5ex]{dl}{\pi_1} \arrow[shift right=0.5ex,swap]{dr}{\pi_2} & \\
			|[alias=X_1]| X_1\arrow[to=X_2,shift left=0.5ex]{}{\pi_2 \circ i_1} \arrow[shift left=0.5ex]{ur}{i_1} && |[alias=X_2]| X_2 \arrow[to=X_1, shift left=0.5ex]{}{\pi_1 \circ i_2} \arrow[shift right=0.5ex,swap]{ul}{i_2}
			\end{tikzcd}
		\end{center}
	commutes, image objects of $ \idem $ are unique as summands of $ X $.

	\end{REM}

Now let $ \C $ be any (linear) category. If there exists an image object for every idempotent in $ \C $ we say that $ \C $ is \emph{idempotent complete}. As this property fails for many categories it can be desirable to fully embed $ \C $ into another category $ \overline{\C} $ that is idempotent complete.

	\begin{DEF}

	An \emph{idempotent completion} of a category $ \C $ is a category $ \overline{\C} $ together with a covariant functor $ \Psi \colon \C \to \overline{\C} $ such that

		\begin{itemize}

		\item $ \Psi $ is fully faithful.

		\item $ \overline{\C} $ is idempotent complete.

		\item For every object $ X $ in $ \overline{\C} $ there exists an idempotent $ \idem $ in $ \C $ such that $ X $ is an image object for $ \Psi(\idem) $.

		\end{itemize}

	\end{DEF}

	\begin{REM}

	Idempotent completions are unique up to equivalence of categories \cite[Section 5.1.4]{Lurie09}.

	\end{REM}

Pleasingly, the Yoneda Lemma provides a universal way for realising the idempotent completion of any category $ \C $. Let $ \RC $ denote the category of \emph{contravariant} functors from $ \C $ into $ \Vect $, we call this the category of representations of $ \C $. We consider the functor
	\begin{align*}
	\Yon \colon \C &\to \RC \\
	X &\mapsto X\Sharp
	\end{align*}
where $ X\Sharp = \Hom_{\C}(\Blank,X) $.

It is a well known corollary of the Yoneda Lemma that $ \Yon $ is fully faithful and is therefore referred to as the \emph{Yoneda embedding}. Furthermore for any idempotent $ \idem\in\End_\C(X) $, there is, in $ \RC $, a subfunctor $ (X,\idem)\Sharp\leq X\Sharp $ given by
	\begin{align*}
	(X,\idem)\Sharp \colon \C &\to \Vect \\
	Y &\mapsto \Hom_{\C}(Y, \idem) \defeq \{ f \in \Hom_{\C}(Y, X) \mid \idem \circ f = f \} \\
	(\alpha \colon Y \to Z) &\mapsto (f \mapsto \alpha \circ f)
	\end{align*}
which is an image object for $ \Yon(\idem) $ in $ \RC $. Indeed, $ (X,\idem)\Sharp $ is a summand of $ X\Sharp $. This image object exists because $ \RC $ is idempotent complete, even if $ \C $ may not be. Concretely, $ (X,\idem)\Sharp(Y) $ is the image of $ \idem\Sharp_Y=\idem_* $, which is an idempotent endomorphism of $ X\Sharp(Y)=\Hom_\C(Y,X) $. The naturality of $ \idem\Sharp $, i.e.\ the fact that $ \idem_* $ commutes with $ \phi^* $ for any $ \phi\colon Z\to Y $, makes $ (X,\idem)\Sharp $ a functor.

 Let $ \overline{\C}_{\Yon} $ be the full subcategory of $ \RC $ spanned by object of the form $ (X,\idem)\Sharp $. We note that $ \Yon $ factors through $ \overline{\C}_{\Yon} $ as $ (X,\id_X)\Sharp = \Yon(X) $ for all $ X $ in $ \C $. We therefore obtain the following result.

 	\begin{PROP} \label{prop:yoneda_gives_idem_completion}

	The category $ \overline{\C}_{\Yon} $ together with the Yoneda embedding is an idempotent completion of $ \C $.

	\end{PROP}

As one should expect, if $ \C $ is semisimple we automatically have $ \C = \overline{\C}_{\Yon} $. However, if $ \C $ is also Schurian then it also coincides with its entire category of representations.

    \begin{PROP} \label{prop:yon_equivalence}

    Let $ \C $ be a semisimple Schurian category with finitely many isomorphism classes of simple objects. Then the Yoneda embedding
        \begin{align*}
        \C \to \RC\\
        X \mapsto X\Sharp
        \end{align*}
    is an equivalence.

    \proof

    As the Yoneda embedding is fully faithful we only have to show that it is essentially surjective. For $ F $ in $ \RC $ and $ X $ in $ \C $, we have
        \begin{align*}
        F(X) &= \bigoplus\limits_{S \in \Irr(\C)} F(S) \otimes \Hom_C(S,X)^* \\
        &= \bigoplus\limits_{S \in \Irr(\C)} F(S) \otimes \Hom_C(X,S) \\
        &= \bigoplus\limits_{S \in \Irr(\C)} F(S) \otimes S\Sharp(X)
        \end{align*}
    where the first equality uses the semisimplicity of $ \C $ and the contravariance of $ F $ and the second equality uses the fact $ S $ is Schurian.    \endproof

    \end{PROP}

Let $ \C $ be a semisimple category together with a complete set of simples $ \Irr(\C) $. Suppose we choose an element $ I \in \Irr(\C) $ and consider the full subcategory of $ \C $ whose objects are non-isomorphic to $ I $. Clearly this new category fails to be semisimple, however the missing simple objects may still be detected by considering the idempotent endomorphisms of any object that has a proper summand isomorphic to $ I $ (this is what was mean before by `semisimple in spirit'). There is, therefore, a notion analogous to a complete set of simples for idempotents.

	\begin{DEF} \label{def:primitive_idempotents}

	A \emph{set of primitive orthogonal idempotents} in a linear category $ \C $ is a set of idempotents $ \I $ in $ \C $ such that
		\begin{align*}
		\Hom_{\C}(\idem,\idem') =
			\begin{cases}
			\fld \quad &\text{if $ \idem = \idem' $} \\
			0 \quad &\text{else.}
			\end{cases}
		\end{align*}
	A set of primitive orthogonal idempotents is called \emph{complete} if we have
		\begin{align*}
		\bigoplus\limits_{\idem \in \I} \Hom_{\C}(X,\idem) \otimes \Hom_{\C}(\idem,Y)  = \Hom_{C}(X,Y)
		\end{align*}
	for all $ X,Y $ in $ \C $.

	\end{DEF}

The proof Proposition~\ref{prop:yon_equivalence} shows that the Yoneda embedding maps a complete set of Schurian simples in $ \C $ to a complete set of Schurian simples
in $ \RC $. The corresponding claim for a complete set of primitive orthogonal idempotents also holds.

	\begin{PROP}

	\label{prop:complete_idem_implies_complete_simples}

	Let $ \C $ be a linear category with a complete set of primitive orthogonal idempotents $ \I $. Then $ \RC $ is a semisimple Schurian category and $ \{(X_\idem,\idem)\Sharp \}_{\idem \in \I} $ forms a complete set of simples in $ \RC $ (where $ \idem \in \End_{\C}(X_\idem) $).

	\proof

	It is straightforward to check that the set $ \{(X_\idem,\idem)\Sharp \}_{\idem \in \I} $ contains distinct simple Schurian objects in $ \RC $. The condition that $ \I $ is a complete set of orthogonal idempotents may be rephrased as
		\begin{align*}
		Y\Sharp = \bigoplus\limits_{\idem \in \I} \Hom_{\C}(\idem,Y) \vprod (X_\idem,\idem)\Sharp
		\end{align*}
	for all $ Y $ in $ \C $. Then, using a similar argument to the proof of Proposition~\ref{prop:yon_equivalence}, we take $ F $ in $ \RC $, $ Y $ in $ \C $ and compute,
	 	\begin{align*}
		F(Y) &= \Hom_{\RC}(Y\Sharp,F) \\
		&= \bigoplus\limits_{\idem \in \I} \Hom_{\C}(\idem,Y)\dual \otimes \Hom_{\RC}((X_\idem,\idem)\Sharp,F)\\
		&= \bigoplus\limits_{\idem \in \I} \Hom_{\C}(Y,\idem) \otimes \Hom_{\RC}((X_\idem,\idem)\Sharp,F)\\
		&= \bigoplus\limits_{\idem \in \I} \Hom_{\RC}((X_\idem,\idem)\Sharp,F) \otimes (X_\idem,\idem)\Sharp(Y)
	 	\end{align*}
	as desired.	\endproof

	\end{PROP}

\section{Preliminaries on Monoidal Categories}  \label{sec:prelim_monoidal_cat}

We now turn our attention to the required preliminaries on monoidal (or tensor) categories. As before the material covered is completely standard and is included for the sake of being self-contained.

For any of the fundamental definitions relating to monoidal categories see~\cite[Chapter 4]{Etingof15}. Throughout this article we shall make extensive use of a graphical notation, due to Roger Penrose~\cite{MR0281657}, which is extremely well suited to rendering computations within monoidal categories more accessible. In particular, its 2-dimentionality is very much a consequence of the fact that a monoidal category is a special kind of 2-category (i.e.\ one with a single object).

The conventions for this graphical notation are as follows. To represent a morphism $ \alpha \in \Hom_{\C}(X,Y) $ we draw a strand labelled $ X $, a strand labelled $ Y $ and a connection between them labelled $ \alpha $ as follows
	\begin{align*}
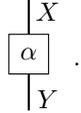
		\begin{array}{c}
			\begin{tikzpicture}[scale=0.15,every node/.style={inner sep=0,outer sep=-1}]
			\node [draw,outer sep=0,inner sep=2,minimum size=15] (v9) at (0,-3) {$ \alpha $};
			\node at (1.7,0.8) {$X$};
			\node at (1.7,-7) {$Y$};
			\draw[thick]  (0,2) edge (v9);
			\draw[thick]  (v9) edge (0,-8);
			\end{tikzpicture}
		\end{array}.
	\end{align*}
Composition is then depicted by vertical juxtaposition and the monoidal product by horizontal juxtaposition:
	\begin{align*}
	\beta \circ \alpha =
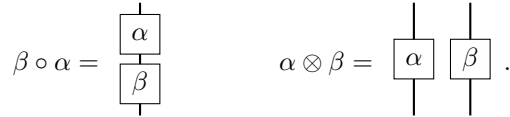
		\begin{array}{c}
			\begin{tikzpicture}[scale=0.15,every node/.style={inner sep=0,outer sep=-1}]
			\node [draw,outer sep=0,inner sep=2,minimum size=15] (v9) at (0,-0.8) {$ \alpha $};
			\node [draw,outer sep=0,inner sep=2,minimum size=15] (v10) at (0,-5.2) {$ \beta $};
			\draw[thick]  (0,2) edge (v9);
			\draw[thick]  (v9) edge (v10);
			\draw[thick]  (v10) edge (0,-8);
			\end{tikzpicture}
		\end{array}
	\qquad \qquad
	\alpha \otimes \beta =
		\begin{array}{c}
			\begin{tikzpicture}[scale=0.15,every node/.style={inner sep=0,outer sep=-1}]
			\node [draw,outer sep=0,inner sep=2,minimum size=15] (v9) at (0,-3) {$ \alpha $};
			\node [draw,outer sep=0,inner sep=2,minimum size=15] (v10) at (5,-3) {$ \beta $};
			\draw[thick]  (0,2) edge (v9);
			\draw[thick]  (v9) edge (0,-8);
			\draw[thick]  (5,2) edge (v10);
			\draw[thick]  (v10) edge (5,-8);
			\end{tikzpicture}
		\end{array}.
	\end{align*}
As the monoidal product is depicted by horizontal juxtaposition the associativity maps are implicitly used but not depicted. Similarly any strand labelled by the tensor identity is not drawn, therefore the unit isomorphisms are also implicitly used but not depicted.

Our goal with the graphical approach is to manipulate diagrams in a manner which leaves the corresponding morphisms unchanged. The most basic manipulation is given by bending the strands around. The resulting strand will no longer be indexed by the original object, however, but by its (left or right) \emph{dual}. To say this more precisely, for an object $ X $ in $ \C  $ a right dual to $ X $ is an object $ X^\vee $ together with morphisms
	\begin{align*}
	\cre_X =
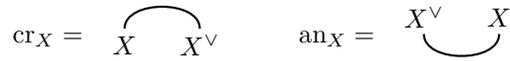
		\begin{array}{c}
			\begin{tikzpicture}
			\node at (-0.5,-0.1) {$X$};
			\node at (0.5,-0.1) {$X^\vee$};
			\node (v1) at (-0.5,0) {};
			\node (v2) at (0.5,0) {};
			\draw[thick] (v1) to[out=90,in=90] (v2);
			\end{tikzpicture}
		\end{array}
	\qquad \ann_X =
		\begin{array}{c}
			\begin{tikzpicture}
			\node at (-0.5,0.1) {$X^\vee$};
			\node at (0.5,0.1) {$X$};
			\node (v1) at (-0.5,0) {};
			\node (v2) at (0.5,0) {};
			\draw[thick] (v1) to[out=-90,in=-90] (v2);
			\end{tikzpicture}
		\end{array}
	\end{align*}
which satisfy the intuitive graphical equation
    \begin{align} \label{eq:both_S}
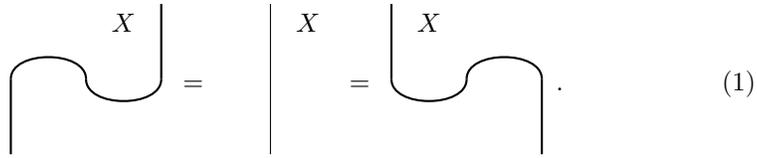
        \begin{array}{c}
            \begin{tikzpicture}[scale = 0.5, xscale = -1]
            \draw[thick] (-1,2) -- (-1,0);
            \draw[thick] (-1, 0) to[out=-90,in=-90] (1, 0);
            \draw[thick] (1, 0) to[out=90,in=90] (3, 0);
            \draw[thick] (3,0) -- (3,-2);
            \node at (0,1.5) {$ X $};
            \end{tikzpicture}
        \end{array}
    =
        \begin{array}{c}
            \begin{tikzpicture}[scale = 0.5]
            \draw (0,2) -- (0,-2);
            \node at (1,1.5) {$ X $};
            \node at (-1,1.5) {};
            \end{tikzpicture}
        \end{array}
    =
        \begin{array}{c}
            \begin{tikzpicture}[scale = 0.5]
            \draw[thick] (-1,2) -- (-1,0);
            \draw[thick] (-1, 0) to[out=-90,in=-90] (1, 0);
            \draw[thick] (1, 0) to[out=90,in=90] (3, 0);
            \draw[thick] (3,0) -- (3,-2);
            \node at (0,1.5) {$ X $};
            \end{tikzpicture}
        \end{array}.
    \end{align}
A left dual to $ X $ is an object $ ^\vee X $ together with morphism $ \cre_{^\vee X} $ and $ \ann_{^\vee X} $ such that $ X $ (together with the same morphisms) is a right dual to $ ^\vee X $. If $ X $ in $ \C $ admits a right (resp. left) dual, then the dual is unique \cite[Proposition 2.10.5]{Etingof15}. Note that, for $ X $ and $ Y $ in $ \C $, the canonical isomorphism
    \begin{align*}
    &(X \otimes Y)^\vee \xrightarrow{\id \otimes \cre_X} (X \otimes Y)^\vee \otimes X \otimes X^\vee \\
    \xrightarrow{\id \otimes \cre_Y \otimes \id} &(X \otimes Y)^\vee \otimes X \otimes Y \otimes Y^\vee \otimes X^\vee \xrightarrow{\ann_{X \otimes Y} \otimes \id} Y^\vee \otimes X^\vee
    \end{align*}
identifies $ (X \otimes Y)^\vee $ and $ Y^\vee \otimes X^\vee $; $ ^\vee(X\otimes Y) $ and $ ^\vee Y \otimes ^\vee \! X $ may be identified similarly. Furthermore, the canonical isomorphism
    \begin{align*}
    \Hom_{\C}(X \otimes Y, Z) &\to \Hom_{\C}(X, Z \otimes Y^\vee)\\
    g &\mapsto (g \otimes \id_{Y\vee}) \circ (\id_X \otimes \ \cre_Y)
    \end{align*}
identifies $ \Hom_{\C}(X \otimes Y, Z) $ and $ \Hom_{\C}(X, Z \otimes Y^\vee) $; the equalities
\begin{align*}
    &\Hom_{\C}(X, Y \otimes Z) = \Hom_{\C}(Y^\vee \otimes X, Z), \\
    &\Hom_{\C}(X \otimes \prescript{\vee}{}{Y}, Z) = \Hom_{\C}(X, Z \otimes Y),\\
    &\Hom_{\C}(Y \otimes X, Z) = \Hom_{\C}(X, \prescript{\vee}{}{Y} \otimes Z)
    \end{align*}
may be derived similarly.  The following rephrasing of these equalities in terms of the Yoneda embedding will prove useful.

    \begin{LEMMA} \label{lem:duals_and_yoneda}

    Let $ \C $ be a rigid monoidal category. Then we have
        \begin{align*}
        Y\Sharp \circ (X^\vee \otimes \Blank) = (X \otimes Y)\Sharp = X\Sharp \circ (\Blank \otimes \prescript{\vee}{}{Y}).
        \end{align*}
    \proof

    The natural isomorphisms described above give the desired result. \endproof

    \end{LEMMA}

A category is called \emph{rigid} is every object admits a left and right dual. There are several simplifications (or complexifications depending on ones perspective!) to this very general setup. For instance a rigid category may be self dual. In this case every object comes with the a pair of maps realising it as a left (and therefore also right) dual to itself. For an example of such a category we can simply consider the category of vector spaces with a preferred basis. However, this is not the simplification we are interested in as there is no issue, from a graphical perspective, in the label of a strand changing when we bend it. However, it is less desirable that the new label should depend on the direction of our bend. In other words we need to identify the left and right duals. Such an identification is called a \emph{pivotal structure}.

	\begin{DEF} \label{def:pivotal_structure}

	Let $ \C $ be a rigid category. A pivotal structure on $ \C $ is a choice of natural isomorphism
		\begin{align*}
		\delta_X \colon ^\vee X \to X^\vee
		\end{align*}
	such that $ \delta_{X \otimes Y} = \delta_Y \otimes \delta_X. $ By convention $ \delta_X $ is suppressed from the graphical notation. A rigid category equipped with a pivotal structure is called a \emph{pivotal category}.

	\end{DEF}

    \begin{REM}

    Although the identity does give a pivotal structure on any self dual category, other choices can be interesting. As an example, we may consider the semisimple quotient of the Temperley-Lieb category with loop parameter $ -[2]_q $ where $ q $ is a primitive $ 2h $ root of unity for some $ h \in \mathbb{Z}_{\geq 2} $. This is a self dual category and is known to be equivalent to the category of integrable highest weight modules at level $ h-2 $ of the affine Lie algebra $ A_1^{(1)} $ as a monoidal category. However, for this equivalence to be pivotal one must equip the Temperley-Lieb category with a non trivial pivotal structure. For further details on this subtlety of pivotal structures see \cite{Snyder09}.

    \end{REM}

The presence of a pivotal structure allows us to `close up' objects in the following sense. For any $ X $ in $ \C $ we may consider the endomorphism of the identity given by
    \begin{align} \label{diag:right_trace}
        \begin{array}{c}
            \begin{tikzpicture}[scale=0.15,every node/.style={inner sep=0,outer sep=-1}]
            \node at (4.5,0) {$ X^\vee $};
            \node at (-3,0) {$ X $};
            \node (v8) at (-1,3) {};
            \node (v10) at (2,3) {};
            \node (v9) at (-1,-3) {};
            \node (v11) at (2,-3) {};
            \draw[thick]  (v8) edge (v9);
            \draw[thick] (v8) to[out=90,in=90] (v10);
            \draw[thick]  (v10) edge (v11);
            \draw[thick] (v9) to[out=-90,in=-90] (v11);
            \end{tikzpicture}
        \end{array}.
    \end{align}
However, under no additional assumptions, \eqref{diag:right_trace} is not necessarily equal to
	\begin{align} \label{diag:left_trace}
		\begin{array}{c}
			\begin{tikzpicture}[scale=0.15,every node/.style={inner sep=0,outer sep=-1}, xscale=-1]
			\node at (4.5,0) {$ X^\vee $};
			\node at (-3,0) {$ X $};
			\node (v8) at (-1,3) {};
			\node (v10) at (2,3) {};
			\node (v9) at (-1,-3) {};
			\node (v11) at (2,-3) {};
			\draw[thick]  (v8) edge (v9);
			\draw[thick] (v8) to[out=90,in=90] (v10);
			\draw[thick]  (v10) edge (v11);
			\draw[thick] (v9) to[out=-90,in=-90] (v11);
			\end{tikzpicture}
		\end{array}.
	\end{align}
Within the scope of this article, however, we shall often consider the case when these two diagrams do coincide. In this case, i.e. when \eqref{diag:right_trace} and \eqref{diag:left_trace} define the same element of $ \End_{\C}(\tid) $ for all $ X $ in $ \C $, the category is called \emph{spherical}. Furthermore, said element is then called the \emph{dimension} of $ X $ and is denoted $ d(X) $.

There is a one final embellishment of monoidal cagoteries that we shall consider in this article. Our end goal is for the diagrammatic calculus to be well defined up to \emph{ribbon tangle isotopy}. As ribbon tangles are fundamentally 3-dimensional objects this will require that we consider a special class of 3-categories. However, just as a 2-category with 1 object may be though of as monoidal 1-category, the relevant class of 3-categories may be though of as the class of monoidal 1-categories with even more additional structure. That additional structure is a \emph{braiding} and is given by a collection of natural isomorphisms $ \sigma_{X,Y} \colon X \otimes Y \to Y \otimes X $, such that the diagrams
	\begin{align*}
		\begin{array}{ccc}
		X \otimes Y \otimes Z & \xrightarrow{\sigma_{X,Y\otimes Z}} & Y \otimes Z \otimes X \\
		{\scriptstyle \sigma_{X,Y} \otimes \tid} \searrow & & \nearrow { \scriptstyle \tid \otimes  \sigma_{X,Z}} \\
		& Y \otimes X \otimes Z
		\end{array}
	\end{align*}
and
	\begin{align*}
		\begin{array}{ccc}
		X \otimes Y \otimes Z & \xrightarrow{\sigma_{X \otimes Y, Z}} & Y \otimes Z \otimes X \\
		{\scriptstyle \tid \otimes \sigma_{Y,Z}} \searrow & & \nearrow { \scriptstyle \sigma_{X,Z} \otimes \tid} \\
		& Y \otimes X \otimes Z
		\end{array}
	\end{align*}
commute. These conditions are often referred to as the \emph{hexagon identities} (our diagrams are triangular as we have suppressed the associativity isomorphisms). In graphical notation the braiding is depicted by the over-crossing,
	\begin{align*}
		\begin{array}{c}
			\begin{tikzpicture}[scale = 0.5]
			\draw[thick] (1,1) -- (-1,-1);
			\draw[line width = 0.3cm, white] (-1, 1) -- (1,-1);
			\draw[thick] (-1, 1) -- (1,-1);
			\node at (-1, 1.6) {$ X $};
			\node at (1, 1.6) {$ Y $};
			\end{tikzpicture}
		\end{array}.
	\end{align*}
The hexagon identities guarantee that this notation is consistent with horizontal juxtaposition depicting the monoidal product. Naturality of the braiding allows morphisms to pass over and under strands, i.e.\
	\begin{align*}
		\begin{array}{c}
			\begin{tikzpicture}[scale = 0.5]
			\draw[thick] (-1,3) -- (-1,1.5);
			\draw[thick] (1,3) -- (1,1.5);
			\draw[thick] (1,1.5) to[out=-90,in=90] (-1,-1);
			\draw[line width = 0.3cm, white] (-1,1.5) to[out=-90,in=90] (1,-1);
			\draw[thick] (-1,1.5) to[out=-90,in=90] (1,-1);
			\node at (-1, 3.6) {$ X $};
			\node at (1, 3.6) {$ Y $};
			\node[draw,outer sep=0,inner sep=2,minimum size=16,fill=white] at (-1,2) {$ \alpha $};
			\end{tikzpicture}
		\end{array}
	=
		\begin{array}{c}
			\begin{tikzpicture}[scale = 0.5,yscale=-1]
			\draw[thick] (-1,3) -- (-1,1.5);
			\draw[thick] (1,3) -- (1,1.5);
			\draw[thick] (-1,1.5) to[out=-90,in=90] (1,-1);
			\draw[line width = 0.3cm, white] (1,1.5) to[out=-90,in=90] (-1,-1);
			\draw[thick] (1,1.5) to[out=-90,in=90] (-1,-1);
			\node at (-1,-1.5) {$ X $};
			\node at (1,-1.5) {$ Y $};
			\node[draw,outer sep=0,inner sep=2,minimum size=16,fill=white] at (1,2) {$ \alpha $};
			\end{tikzpicture}
		\end{array}
	\quad \text{and} \quad
		\begin{array}{c}
			\begin{tikzpicture}[scale = 0.5]
			\draw[thick] (-1,3) -- (-1,1.5);
			\draw[thick] (1,3) -- (1,1.5);
			\draw[thick] (1,1.5) to[out=-90,in=90] (-1,-1);
			\draw[line width = 0.3cm, white] (-1,1.5) to[out=-90,in=90] (1,-1);
			\draw[thick] (-1,1.5) to[out=-90,in=90] (1,-1);
			\node at (-1, 3.6) {$ X $};
			\node at (1, 3.6) {$ Y $};
			\node[draw,outer sep=0,inner sep=2,minimum size=16,fill=white] at (1,2) {$ \beta $};
			\end{tikzpicture}
		\end{array}
	=
		\begin{array}{c}
			\begin{tikzpicture}[scale = 0.5,yscale=-1]
			\draw[thick] (-1,3) -- (-1,1.5);
			\draw[thick] (1,3) -- (1,1.5);
			\draw[thick] (-1,1.5) to[out=-90,in=90] (1,-1);
			\draw[line width = 0.3cm, white] (1,1.5) to[out=-90,in=90] (-1,-1);
			\draw[thick] (1,1.5) to[out=-90,in=90] (-1,-1);
			\node at (-1,-1.5) {$ X $};
			\node at (1,-1.5) {$ Y $};
			\node[draw,outer sep=0,inner sep=2,minimum size=16,fill=white] at (-1,2) {$ \beta $};
			\end{tikzpicture}
		\end{array}.
	\end{align*}
If $ \sigma $ is a braiding on a monoidal category $ \C $ then $ \overline{\sigma} $, given by $ \overline{\sigma}_{X,Y} \defeq (\sigma_{Y,X}^{-1}) $, also defines a braiding on $ \C $ called the \emph{opposite} braiding to $ \sigma $. In graphical notation the opposite braiding is depicted by the under-crossing,
	\begin{align*}
		\begin{array}{c}
			\begin{tikzpicture}[scale = 0.5]
			\draw[thick] (-1, 1) -- (1,-1);
			\draw[line width = 0.3cm, white] (1,1) -- (-1,-1);
			\draw[thick] (1,1) -- (-1,-1);
			\node at (-1, 1.6) {$ X $};
			\node at (1, 1.6) {$ Y $};
			\end{tikzpicture}
		\end{array}.
	\end{align*}
We therefore have the ``Reidemeister II" rule
	\begin{align*}
		\begin{array}{c}
			\begin{tikzpicture}[scale = 0.4]
			\draw[thick] (1,1) to[out=-90,in=90] (-1,-1);
			\draw[line width = 0.3cm, white] (-1, 1) to[out=-90,in=90] (1,-1);
			\draw[thick] (-1, -1) to[out=-90,in=90] (1,-3);
			\draw[line width = 0.3cm, white] (1,-1) to[out=-90,in=90] (-1,-3);
			\draw[thick] (1,-1) to[out=-90,in=90] (-1,-3);
			\draw[thick] (-1, 1) to[out=-90,in=90] (1,-1);
			\node at (-1, 1.7) {$ X $};
			\node at (1, 1.7) {$ Y $};
			\end{tikzpicture}
		\end{array}
	=
		\begin{array}{c}
			\begin{tikzpicture}[scale = 0.4]
			\draw[thick] (1,1) -- (1,-3);
			\draw[thick] (-1, 1) -- (-1,-3);
			\node at (-1, 1.7) {$ X $};
			\node at (1, 1.7) {$ Y $};
			\end{tikzpicture}
		\end{array}.
    \end{align*}
This, together with \eqref{eq:both_S} (which resembles a ``Reidemeister 0" rule), makes one wonder whether, for a spherical braided category, our goal has been achieved and graphical notation is well defined up to tangle isotopy. In general, however, this may still fail as we need to assume a certain compatibility between the braiding and the pivotal structure which may be described graphically as
	\begin{align*}
		\begin{array}{c}
			\begin{tikzpicture}[every node/.style={inner sep=0,outer sep=-1}]]
			\node (v5) at (-0.25,-0.375) {};
			\node (v6) at (0.125,-0.375) {};
			\node at (-0.25,0) {$ X $};
			\node at (0.25,0) {$ Y $};
			\node (v1) at (-0.25,-0.75) {};
			\node (v3) at (0.125,-0.625) {};
			\node (v2) at (0.375,-1.5) {};
			\node (v4) at (0.75,-1.375) {};
			\node (v7) at (0.75,-1.375) {};
			\node (v9) at (1,-1.375) {};
			\node (v8) at (1.375,-1.5) {};
			\node (v12) at (1,-1) {};
			\node (v10) at (1.375,-0.875) {};
			\node (v11) at (0.375,-0.875) {};
			\node (v14) at (0.75,-1) {};
			\node (v15) at (0,-1.875) {};
			\node (v16) at (-0.375,-1.75) {};
			\node (v13) at (0.75,-1) {};
			\node (v17) at (0,-2) {};
			\node (v18) at (-0.375,-2) {};
			\draw [thick] (v7) to[out=-90,in=-90] (v9);
			\draw [thick] (v10) to[out=90,in=90] (v11);
			\draw [thick] (v12) to[out=90,in=90] (v13);
			\draw [thick] (v14) to[out=-90,in=90] (v15);
			\draw [thick] (v11) to[out=-90,in=90] (v16);
			\draw[line width=4,white]  (v1) to[out=-90,in=90] (v2);
			\draw[line width=4,white]  (v3) to[out=-90,in=90] (v4);
			\draw[thick]  (v1) to[out=-90,in=90] (v2);
			\draw [line width=4,white] (v2) to[out=-90,in=-90] (v8);
			\draw [thick] (v2) to[out=-90,in=-90] (v8);
			\draw[thick]  (v3) to[out=-90,in=90] (v4);
			\draw [thick] (v5) edge (v1);
			\draw [thick] (v6) edge (v3);
			\draw [thick] (v4) edge (v7);
			\draw [thick] (v10) edge (v8);
			\draw [thick] (v12) edge (v9);
			\draw [thick] (v13) edge (v14);
			\draw [thick] (v15) edge (v17);
			\draw [thick] (v16) edge (v18);
			\end{tikzpicture}
		\end{array}
	=
		\begin{array}{c}
			\begin{tikzpicture}[scale = 0.5,every node/.style={inner sep=0,outer sep=-1}]]
			\draw[thick] (0.5,0.5) to[out=-90,in=90] (-1,-0.5);
			\draw[line width = 0.2cm, white] (-1,0.5) to[out=-90,in=90] (0.5,-0.5);
			\draw[thick] (0.5,-0.5) to[out=-90,in=90] (-1,-1.5) node (v1) {};
			\draw[thick] (-1,0.5) to[out=-90,in=90] (0.5,-0.5);
			\draw[line width = 0.2cm, white] (-1,-0.5) to[out=-90,in=90] (0.5,-1.5) node (v16) {};
			\draw[thick] (-1,-0.5) to[out=-90,in=90] (0.5,-1.5) node (v16) {};
			\node at (-1,1) {$ X $};
			\node at (0.5,1) {$ Y $};
			\node (v2) at (-1,-2) {};
			\node (v3) at (-0.5,-2.5) {};
			\node (v6) at (0,-2.5) {};
			\node (v8) at (0,-2) {};
			\node (v11) at (-0.5,-2) {};
			\node (v12) at (-1,-2.5) {};
			\node (v14) at (-1,-3) {};
			\node (v4) at (0.5,-2) {};
			\node (v13) at (0.5,-2.5) {};
			\node (v10) at (1,-2) {};
			\node (v5) at (1,-2.5) {};
			\node (v7) at (1.5,-2.5) {};
			\node (v9) at (1.5,-2) {};
			\node (v15) at (0.5,-3) {};
			\draw [thick] (v1) edge (v2);
			\draw [thick] (v3) to[out=-90,in=-90] (v6);
			\draw [thick] (v5) to[out=-90,in=-90] (v7);
			\draw [thick] (v6) edge (v8);
			\draw [thick] (v7) edge (v9);
			\draw [thick] (v9) to[out=90,in=90] (v10);
			\draw [thick] (v8) to[out=90,in=90] (v11);
			\draw [thick] (v11) to[out=-90,in=90] (v12);
			\draw [thick] (v10) to[out=-90,in=90] (v13);
			\draw [thick] (v12) edge (v14);
			\draw [thick] (v13) edge (v15);
			\draw [thick] (v16) edge (v4);
			\draw [line width = 0.2cm, white] (v2) to[out=-90,in=90] (v3);
			\draw [line width = 0.2cm, white] (v4) to[out=-90,in=90] (v5);
			\draw [thick] (v2) to[out=-90,in=90] (v3);
			\draw [thick] (v4) to[out=-90,in=90] (v5);
			\end{tikzpicture}
		\end{array}.
	\end{align*}
A braiding which satisfies this condition is called balanced; for the remainder of this article we suppose that all braidings on spherical categories are balanced. Within this context graphical notation is now well defined up to ribbon tangle isotopy~\cite{Res90}, for this reason we call a (balanced) braided spherical category a \emph{ribbon} category. Note that in general, our diagrams are not well-defined up to tangle isotopy, as
	\begin{align*}
		\begin{array}{c}
			\begin{tikzpicture}[scale = 0.3]
			\draw[thick] (-1,1) -- (-1,3);
			\draw[thick] (-1,-1) -- (-1,-3);
			\draw[thick] (1,1) to[out = 90, in = 90] (3,1);
			\draw[thick] (1,-1) to[out = -90, in = -90] (3,-1);
			\draw[thick] (3,1) -- (3,-1);
			\draw[thick] (1,1) to[out=-90,in=90] (-1,-1);
			\draw[line width = 0.3cm, white] (-1, 1) to[out=-90,in=90] (1,-1);
			\draw[thick] (-1, 1) to[out=-90,in=90] (1,-1);
			\node at (-1, 4) {$ X $};
			\end{tikzpicture}
		\end{array}
	\end{align*}
will often fail to denote $ \id_X \in \End_{\C}(X) $. Furthermore we shall see that this failure (which is precisely the failure of our graphical calculus to be well defined up to tangle isotopy) plays a crucial role in the appearance of certain interesting aspects of the theory.

\section{Preliminaries on Fusion Categories} \label{sec:prelim_fusion_cat}

A \emph{fusion category} is a Schurian, semisimple, rigid, monoidal category $ \C $ such that there is a complete set of simples $ \Irr(\C) $ satisfying $ |\Irr(\C)| < \infty $, $ \tid \in \Irr(\C) $ and $$ d(\C):= \sum_{I \in \Irr(\C)} d(I) \neq 0. $$ The condition on $ d(\C) $ is often ommitted as it holds automatically when the underlying field is algebraically closed~\cite[Proposition 8.20.16]{Etingof15}. A \emph{spherical fusion category} is simply a fusion category equipped with a spherical pivotal structure. 
 The aim of this section is to gather certain results on such categories. For the sake of expediency many of the proofs are omitted, however references are always provided. These results are all well known to experts.

So far we have only discussed graphical notation within the context of monoidal categories; in the more specific case of semisimple Schurian categories it is possible to decompose a given strand, as described by the following lemma.

    \begin{LEMMA} 	\label{lem:decompose_object}

    Let $ \C $ be a semisimple Schurian category and let $ X $ be in $ \C $. We have
        \begin{align} \label{eq:lemma_decompose}
            \begin{array}{c}
                \begin{tikzpicture}[scale=0.15,every node/.style={inner sep=0,outer sep=-1}]
                \node (v8) at (0,7.5) {};
                \node (v10) at (0,-7.5) {};
                \node at (-2,2) {$X$};
                \draw[thick]  (v8) edge (v10);
                \end{tikzpicture}
                \end{array}
                = \sum\limits_{R,b}
            \begin{array}{c}
                \begin{tikzpicture}[scale=0.15,every node/.style={inner sep=0,outer sep=-1}]
                \node (v2) at (0,7.5) {};
                \node (v4) at (0,-7.5) {};
                \node [draw,outer sep=0,inner sep=2,minimum width=14,minimum height = 12] (v6) at (0,3.5) {$ b^* $};
                \node [draw,outer sep=0,inner sep=2,minimum width=14,minimum height = 12] (v9) at (0,-3.5) {$ b $};
                \node at (-1.75,6.5) {$X$};
                \node at (-1.75,-6.75) {$X$};
                \node at (1.75,0) {$R$};
                \draw[thick]  (v6) edge (v9);
                \draw[thick]  (v2) edge (v6);
                \draw[thick]  (v9) edge (v4);
                \end{tikzpicture}
            \end{array}
        \end{align}
    where $ R $ ranges over a basis of $ \Irr(\C) $, $ \{b\} $ is a basis of $ \Hom_{\C}(R,X) $ and $ \{b^*\} $ is the corresponding dual basis with respect to the perfect pairing given by
        \begin{align} \label{eq:perfect_pairing}
        \begin{split}
        \Hom_{\C}(R,X) \otimes \Hom_{\C}(X,R) &\to \End_{\C}(R) \\
        f \otimes g &\mapsto g \circ f.
        \end{split}
        \end{align}

    \proof A proof is provided by Lemma 3.3 in~\cite{hardiman_king}. \endproof

    \end{LEMMA}

As decompositions of the form~\eqref{eq:lemma_decompose} come up fairly frequently, we now make the following conventions: unless otherwise specified, a sum over a variable object in $ \C $ ranges over $ \Irr(\C) $ and  a sum over a variable morphism in $ \C $ ranges over a basis of the appropriate $ \Hom $-space.

A standard technique when using graphical calculus in a spherical fusion category is to decompose a stand (using Lemma~\ref{lem:decompose_object}), rearrange the diagram (up to ribbon isotopy), and then reapply Lemma~\ref{lem:decompose_object} to pack the simple stands back into an object. During the rearranging step of this procedure the stands attached to certain morphisms may get pulled up or down. This can affect the final re-packaging step, as described by the following lemma.

	\begin{LEMMA} \label{lem:dual_decompose}

	Let $ \C $ be a spherical fusion category. Let $ X $ be in $ \C $ and $ S $ be in $ \Irr(\C) $. We have
		\begin{align*}
        \sum\limits_{T,b} d(T)\
            \begin{array}{c}
                \begin{tikzpicture}[scale=0.15,every node/.style={inner sep=0,outer sep=-1}]
                \node (v13) at (-0.625,1.375) {};
                \node (v130) at (-0.75,-1.875) {};
                \node (v6) at (0,7) {};
                \node (v2) at (0,-7.5) {};
                \node [draw,outer sep=0,inner sep=2,minimum width=14] (v9) at (0,3) {$ b $};
                \node [draw,outer sep=0,inner sep=2,minimum width=14] (v90) at (0,-3.5) {$ b^* $};
                \node at (3,-6.5) {$S$};
                \node at (-5.25,-6.5) {$X$};
                \node at (3,6) {$S$};
                \node at (3,-0.25) {$T$};
                \node at (-5.25,6) {$X$};
                \draw[thick]  (v6) edge (v9);
                \draw [thick] (v90) edge (v2);
                \draw [thick] plot[smooth, tension=.7] coordinates {(0.5441,1.6) (0.5441,-2.1)};
                \node (v1) at (-3,1.375) {};
                \node (v4) at (-3,-1.875) {};
                \node (v5) at (-3,-7.5) {};
                \node (v3) at (-3,7) {};
                \draw [thick] (v13) to[out=-90,in=-90] (v1);
                \draw [thick] (v3) edge (v1);
                \draw [thick] (v130) to[out=90,in=90] (v4);
                \draw [thick] (v5) edge (v4);
                \end{tikzpicture}
            \end{array}
		=
		d(S)
			\begin{array}{c}
				\begin{tikzpicture}[scale=0.15,every node/.style={inner sep=0,outer sep=-1}]
				\node (v5) at (-2,6) {};
				\node (v7) at (-2,-6) {};
				\node (v8) at (2,6) {};
				\node (v10) at (2,-6) {};
				\node at (-4,2) {$X$};
				\node at (3.5,2) {$S$};
				\draw[thick]  (v5) edge (v7);
				\draw[thick]  (v8) edge (v10);
				\end{tikzpicture}
			\end{array}
		\end{align*}
	and
		\begin{align*}
        \sum\limits_{T,b} d(T) \
			\begin{array}{c}
				\begin{tikzpicture}[scale=0.15,every node/.style={inner sep=0,outer sep=-1},yscale=1,xscale=-1]
				\node (v13) at (-0.7,1.5) {};
				\node (v130) at (-0.7,-2) {};
				\node (v6) at (0,7) {};
				\node (v2) at (0,-7.5) {};
				\node [draw,outer sep=0,inner sep=2,minimum width=14] (v9) at (0,3) {$ b $};
				\node [draw,outer sep=0,inner sep=2,minimum width=14] (v90) at (0,-3.5) {$ b^* $};
				\node at (2,-6.5) {$S$};
				\node at (-4.5,-6.5) {$X$};
				\node at (2,6) {$S$};
				\node at (2,-0.25) {$T$};
				\node at (-4.5,6) {$X$};
				\draw[thick]  (v6) edge (v9);
				\draw [thick] (v90) edge (v2);
				\draw [thick] plot[smooth, tension=.7] coordinates {(0.5441,1.6) (0.5441,-2.1)};
				\node (v1) at (-3,1.5) {};
				\node (v4) at (-3,-2) {};
				\node (v5) at (-3,-7.5) {};
				\node (v3) at (-3,7) {};
				\draw [thick] (v13) to[out=-90,in=-90] (v1);
				\draw [thick] (v3) edge (v1);
				\draw [thick] (v130) to[out=90,in=90] (v4);
				\draw [thick] (v5) edge (v4);
				\end{tikzpicture}
			\end{array}
		=
		d(S)
			\begin{array}{c}
				\begin{tikzpicture}[scale=0.15,every node/.style={inner sep=0,outer sep=-1}]
				\node (v5) at (-2,6) {};
				\node (v7) at (-2,-6) {};
				\node (v8) at (2,6) {};
				\node (v10) at (2,-6) {};
				\node at (-4,2) {$S$};
				\node at (3.5,2) {$X$};
				\draw[thick]  (v5) edge (v7);
				\draw[thick]  (v8) edge (v10);
				\end{tikzpicture}
			\end{array}.
		\end{align*}
	\proof A proof is provided by Lemma 5.1 in~\cite{MR2430629}, or alternatively, Lemma 3.11 in~\cite{hardiman_king}. \endproof

    \end{LEMMA}

An interesting consequence of this lemma is the following relationship between the dimension of simples objects and the dimension of the corresponding $ \Hom $-spaces. We consider the expression
    \begin{align*}
    \sum\limits_{S,T,b} d(S) d(T)
		\begin{array}{c}
			\begin{tikzpicture}[scale=0.15,every node/.style={inner sep=0,outer sep=-1}]
			\node (v3) at (-1,2.5) {};
			\node (v5) at (1,2.5) {};
			\node (v4) at (-1,-2.5) {};
			\node (v6) at (1,-2.5) {};
			\draw [thick] (v3) edge (v4);
			\draw [thick] (v5) edge (v6);
			\node [draw,outer sep=0,inner sep=2,minimum width=14,fill=white] (v1) at (0,3.5) {$ b $};
			\node [draw,outer sep=0,inner sep=2,minimum width=14,fill=white] (v2) at (0,-3.5) {$ b^* $};
			\draw[thick] (0,7.75) -- (v1);
			\draw[thick] (v2) -- (0,-7.75);
			\node at (2,6.75) {$ R $};
			\node at (2,-6.75) {$ R $};
			\node at (-2.5,0) {$ S $};
			\node at (2.5,0) {$ T $};
			\end{tikzpicture}
		\end{array}.
    \end{align*}
As, by definition, $ b $ and $ b^* $ compose to the identity this expression is equal to the scalar $ \sum_{S,T} \hom_{\C}(R,ST) d(S) d(T) $ in $ \End_{\C}(R) = \fld $, where $ \hom_{\C}(R,ST) $ denotes the dimension of $ \Hom_{\C}(R,ST) $. However, by rearranging the diagram as follows:
    \begin{align*}
\sum\limits_{S,T,b} d(S) d(T)
    	\begin{array}{c}
			\begin{tikzpicture}[scale=0.15,every node/.style={inner sep=0,outer sep=-1}]
			\node (v3) at (1,2.5) {};
			\node (v5) at (-5.5,5.75) {};
			\node (v4) at (1,-2.5) {};
			\node (v6) at (-5.5,-5.75) {};
			\draw [thick] (v3) edge (v4);
			\draw [thick] (v5) edge (v6);
			\node at (2,6.75) {$ R $};
			\node at (2,-6.75) {$ R $};
			\node at (-7,0) {$ S $};
			\node at (2.5,0) {$ T $};
			\node (v7) at (-1,2.25) {};
			\node (v8) at (-3.25,2.25) {};
			\node (v11) at (-3.25,5.75) {};
			\node (v12) at (-3.25,-5.75) {};
			\draw [thick] (v7) to[out=-90,in=-90] (v8);
			\draw [thick] (v11) to[out=90,in=90] (v5);
			\node (v10) at (-3.25,-2.25) {};
			\node (v9) at (-1,-2.25) {};
			\draw [thick] (v8) edge (v11);
			\draw [thick] (v9) to[out=90,in=90] (v10);
			\draw [thick] (v12) to[out=-90,in=-90] (v6);
			\draw [thick] (v10) edge (v12);
			\node [draw,outer sep=0,inner sep=2,minimum width=14,fill=white] (v1) at (0,3.5) {$ b $};
			\node [draw,outer sep=0,inner sep=2,minimum width=14,fill=white] (v2) at (0,-3.5) {$ b^* $};
			\draw[thick] (0,7.75) -- (v1);
			\draw[thick] (v2) -- (0,-7.75);
			\end{tikzpicture}
    	\end{array}
    \end{align*}
we may apply Lemma~\ref{lem:dual_decompose} and evaluate the expression to be $ d(R) d(\C) $ giving the following lemma.

	\begin{LEMMA} \label{lem:double_decompose}

	Let $ R $ be in $ \Irr(\C) $. Then
		\begin{align*}
		\sum\limits_{S,T} \hom_{\C}(R,ST) d(S) d(T) = d(R) d(\C).
		\end{align*}
    \end{LEMMA}

The only categorical structure described in Section~\ref{sec:prelim_monoidal_cat} still to be reintroduced is that of a braiding. To remedy this situation we introduce the notion of a \emph{pre-modular tensor category} (PTC) which is a (balanced) braided spherical fusion category. We note that such category is, in particular, a ribbon category.

This is clearly loaded terminology. So what additional property must $ \C $ satisfy to acheive the full status of modular tensor category? The answers lies in something called the \emph{modular data} of $ \C $. As $ \C $ is fusion it admits a finite complete set of simples $ \Irr(\C) $. We consider the following $ \Irr(\C) \times \Irr(\C) $ matrices,
	\begin{align} \label{eq:S_and_T_matrices}
	\Tc_{IJ} \defeq \delta_{I,J} \
		\begin{array}{c}
			\begin{tikzpicture}[scale = 0.5]
			\draw[thick] (-1,1) -- (-1,3);
			\draw[thick] (-1,-1) -- (-1,-3);
			\draw[thick] (1,1) to[out = 90, in = 90] (3,1);
			\draw[thick] (1,-1) to[out = -90, in = -90] (3,-1);
			\draw[thick] (3,1) -- (3,-1);
			\draw[thick] (1,1) to[out=-90,in=90] (-1,-1);
			\draw[line width = 0.3cm, white] (-1, 1) to[out=-90,in=90] (1,-1);
			\draw[thick] (-1, 1) to[out=-90,in=90] (1,-1);
			\node at (-1, 3.6) {$ I $};
			\end{tikzpicture}
		\end{array}
	\quad
	\Sc_{IJ} \defeq
		\begin{array}{c}
			\begin{tikzpicture}
			\draw [thick](-1.5,0) arc (180:0:1);
			\draw[white,line width = 0.3cm] (0.5,0) circle (1);
			\draw[thick] (0.5,0) circle (1);
			\node at (-1.6,1) {$ I $};
			\node at (1.4,1) {$ J $};
			\draw [white,line width = 0.3cm](.5,0) arc (0:-180:1);
			\draw [thick](.5,0) arc (0:-180:1);
			\end{tikzpicture}
		\end{array}.
	\end{align}

	\begin{REM} \label{rem:replace_strands}

	If we replace the $ I $ strand in the definition of $ \Tc_{IJ} $ with an $ I^\vee $ strand we obtain an equivalent definition~\cite[Section 8.10.]{Etingof15}. Similarly if we replace the $ I $ and $ J $ strands in the definition of $ \Sc_{IJ} $ with $ I^\vee $ and $ J^\vee $ strands respectively we also obtain an equivalent definition~\cite[Remark 8.1.12.3]{Etingof15}.

	\end{REM}

These matrices are called the T-matrix and the S-matrix respectively. Collectively they are know are the modular data of $ \C $. A \emph{modular tensor category} (MTC) is a pre-modular tensor category such that the above defined $ \Sc $ and $ \Tc $ matrices are invertible, i.e.\ the modular data is non-singular.

The main difference between graphical calculus in modular tensor categories (as opposed to pre-modular tensor categories) is the so called `killing ring' lemma. It goes as follows.

	\begin{LEMMA} \label{lemma:killing_ring}

	Let $ \C $ be an MTC and let $ R $ be in $ \Irr(\C) $. Then
		\begin{align*}
		\sum\limits_S d(S)
			\begin{array}{c}
				\begin{tikzpicture}[scale=0.15,every node/.style={inner sep=0,outer sep=-1}]
				\draw[thick] (-6,4) node (v4) {} to[out=90,in=90] (-3,4) node (v1) {};
				\draw[thick] (-3,-4) node (v3) {} to[out=-90,in=-90] (0,-4) node (v2) {};
				\draw[thick] (-2,1) -- (-3,0) -- (-6,-3) -- (-6,-6);
				\draw[thick] (-1,2) -- (0,3) -- (0,6);
				\draw [thick](v3) -- (-3,-3) -- (-4,-2) ;
				\draw [thick](v1) -- (-3,3) -- (0,0);
				\draw [thick](0,0) -- (0,-1.5) -- (v2);
				\node at (-7.5,2) {$ S $};
				\node at (2.2,-2) {$ S^\vee $};
				\node at (1.5,4.5) {$ R $};
				\draw[thick] (v4) -- (-6,0) -- (-5,-1);
				\end{tikzpicture}
			\end{array}
		= \delta_{R,\tid} d(\C)
		\end{align*}
	where $ \tid $ is the tensor identity.

	\proof

	A proof is provided by Corollary 3.1.11 in \cite{Bak01}. \endproof

	\end{LEMMA}

Heuristically speaking, the ring labelled $ S $ (as weighted by $ d(S) $) may be thought of as `killing' any strand that passes through it. This allows the ring to `slice horizontally' through a diagram, as illustrated by the following corollary.

	\begin{COR} \label{cor:killing_ring}

	Let $ \C $ be an MTC and let $ X $ and $ Y $ be in $ \C $. Then
		\begin{align*}
		\sum\limits_{S} d(S)
			\begin{array}{c}
				\begin{tikzpicture}[scale=0.15,every node/.style={inner sep=0,outer sep=-1}]
				\draw[thick] (-6,4) node (v4) {} to[out=90,in=90] (-3,4) node (v1) {};
				\draw[thick] (0,-7) node (v3) {} to[out=-90,in=-90] (3,-7) node (v2) {};
				\draw[thick] (1,-2) -- (0,-3) -- (-3,-6) -- (-3,-9);
				\draw[thick] (-2,1) -- (-3,0) -- (-7,-4) -- (-7,-9);
				\draw[thick] (2,-1) -- (4,1) -- (4,7);
				\draw[thick] (-1,2) -- (0,3) -- (0,7);
				\draw [thick] (-2,-4) -- (-4,-2) ;
				\draw [thick](v1) -- (-3,3) -- (3,-3);
				\draw [thick](3,-3) -- (3,-4.5) -- (v2);
				\node at (-7.5,2.5) {$ S $};
				\node at (5,-5) {$ S^\vee $};
				\node at (1.5,6.25) {$ X $};
				\node at (5.5,6.25) {$ Y $};
				\draw[thick] (v4) -- (-6,0) -- (-5,-1);
				\draw[thick] (-1,-5) -- (0,-6) -- (v3);
				\end{tikzpicture}
			\end{array}
		=  d(\C) \sum\limits_{T,b,c} \frac{1}{d(T)}
			\begin{array}{c}
				\begin{tikzpicture}[scale=0.15,every node/.style={inner sep=0,outer sep=-1}]
				\draw[thick] (-6,-7) node (v4) {} to[out=90,in=90] (-2,-7) node (v1) {};
				\draw[thick] (-6,-3) node (v3) {} to[out=-90,in=-90] (-2,-3) node (v2) {};
				\node at (-0.5,-5) {$ T^\vee $};
				\node at (-7,-5) {$ T $};
				\node at (-3.75,-12.75) {$ Y $};
				\node at (-7.75,-12.75) {$ X $};
				\node at (-3.75,2.75) {$ Y $};
				\node at (-7.75,2.75) {$ X $};
				\draw[thick] (v4) -- (-6,-13.5);
				\draw [thick](v1) -- (-2,-13.5);
				\draw [thick](v3) -- (-6,3.5);
				\draw [thick](v2) -- (-2,3.5);
				\node [draw,outer sep=0,inner sep=1,minimum width=14, minimum height = 13,fill=white] (v9) at (-6,-1) {$ b^* $};
				\node [draw,outer sep=0,inner sep=1,minimum width=14, minimum height = 13,fill=white] (v9) at (-2,-1) {$ c^* $};
				\node [draw,outer sep=0,inner sep=1,minimum width=14, minimum height = 13,fill=white] (v9) at (-2,-9) {$ c $};
				\node [draw,outer sep=0,inner sep=1,minimum width=14, minimum height = 13,fill=white] (v9) at (-6,-9) {$ b $};
				\end{tikzpicture}
			\end{array}.
		\end{align*}

	\proof A straightforward generalisation of the proof of Corollary 3.13 in \cite{hardiman_king} proves this result. \endproof

    \end{COR}

    \begin{REM}

    The statement of Corollary~\ref{cor:killing_ring} implicitly assumes that $ d(S) \neq  0 $ for all $ S \in \Irr(\C) $. This holds automatically for any pivotal semisimple category~\cite[Proposition 4.8.4.]{Etingof15}.

    \end{REM}

Although it requires a bit more work to see, the killing ring may also be used to slice a diagram vertically.

	\begin{PROP} \label{prop:twisted_S}

	Let $ \C $ be an MTC. We consider $ X,Y,A,B $ in $ \C $, and $ \alpha \in \Hom_{\C}(IJ,XY) $. Then
		\begin{align*}
		\sum\limits_{S} d(S)
			\begin{array}{c}
				\begin{tikzpicture}[scale=0.15,every node/.style={inner sep=0,outer sep=-1}]
				\node at (-7.5,2.5) {$ S $};
				\node at (5,-5) {$ S^\vee $};
				\node at (1.5,5) {$ X $};
				\node at (-8,-17) {$ A $};
				\node at (5.5,5) {$ Y $};
				\node at (-3.5,-17) {$ B $};
				\draw[thick] (-6,4) node (v4) {} to[out=90,in=90] (-3,4) node (v1) {};
				\draw[thick] (0,-16) node (v3) {} to[out=-90,in=-90] (3,-16) node (v2) {};
				\draw [thick](v1) -- (-3,3) -- (1,-1);
				\draw [thick]  (2,-2) -- (3,-3) -- (3,-4.5) -- (v2);
				\draw[thick] (-1,-13.5) -- (0,-14.5) -- (v3);
				\draw [thick] (-2,-12.5) -- (-4,-10.5) ;
				\draw[thick] (v4) -- (-6,-8.5) -- (-4,-10.5);
				\draw[thick] (-2,1) -- (-3,0) -- (-3,-4.5) node (v5) {};
				\draw[thick] (4,6) -- (4,1) -- (2,-1) -- (0,-3) -- (0,-4.5) node (v7) {};
				\draw[thick] (-1,2) -- (0,3) -- (0,6);
				\draw [thick] (-3,-6.5) -- (-3,-8.5) node (v12) {} -- (-4,-9.5);
				\draw [thick](-5,-10.5) -- (-6,-11.5) -- (-6,-18);
				\draw [thick](-2,-18) -- (-2,-13.5) -- (0,-11.5) -- (0,-6.5) node (v14) {};
				\node [draw,outer sep=0,inner sep=1,minimum width=22, minimum height = 13,fill=white] (v9) at (-1.5,-5.5) {$ \alpha $};
				\end{tikzpicture}
			\end{array}
		=   d(\C)\sum\limits_{T,b,c} \frac{1}{d(T)}
			\begin{array}{c}
				\begin{tikzpicture}[scale=0.15,every node/.style={inner sep=0,outer sep=-1}]
				\node (v5) at (-2,10) {};
				\node (v7) at (-2,-14) {};
				\node (v8) at (7,10) {};
				\node (v10) at (7,-14) {};
				\node at (-0.25,9.25) {$ X $};
				\node at (9,9.25) {$ Y $};
				\node at (0,-13) {$A$};
				\node at (9,-13) {$B$};
				\draw[thick]  (v5) edge (v7);
				\draw[thick]  (v8) edge (v10);
				\node (v1) at (1.25,-1) {};
				\node (v2) at (1.25,5.25) {};
				\node (v11) at (1.25,-3) {};
				\node (v6) at (1.25,-9) {};
				\node (v4) at (3.75,-9) {};
				\node (v3) at (3.75,5.25) {};
				\draw [thick] (v1) edge (v2);
				\draw [thick] (v2) to[out=90,in=90] (v3);
				\draw [thick] (v3) edge (v4);
				\draw [thick] (v4) to[out=-90,in=-90] (v6);
				\draw [thick] (v6) edge (v11);
				\node [draw,outer sep=0,inner sep=1,minimum width=22, minimum height = 13,fill=white] (v9) at (-0.25,-2) {$ \alpha $};
				\node [draw,outer sep=0,inner sep=1,minimum width=14, minimum height = 13,fill=white] (v9) at (1.25,-6.5) {$ b^* $};
				\node [draw,outer sep=0,inner sep=1,minimum width=14, minimum height = 13,fill=white] (v9) at (1.25,2.5) {$ c $};
				\node [draw,outer sep=0,inner sep=1,minimum width=14, minimum height = 13,fill=white] (v9) at (7,2.5) {$ c^* $};
				\node [draw,outer sep=0,inner sep=1,minimum width=14, minimum height = 13,fill=white] (v9) at (7,-6.5) {$ b $};
				\node at (8.75,-2.25) {$ T $};
				\end{tikzpicture}
			\end{array}
		\end{align*}
	\proof A straightforward generalisation of the proof of Proposition 3.14 in \cite{hardiman_king} proves this result. \endproof

	\end{PROP}

    This concludes the preliminary portion of this article. We now move on to the construction and characterisation of the tube category.

\section{The Tube Category} \label{sec:tube_category}

Let $ \C $ be a spherical fusion category. The \emph{tube category}, denoted $ \TC $, shares the same objects as $ \C $ but has more morphisms i.e.\ $ \Hom_{\C}(X,Y) \leq \Hom_{\TC}(X,Y) $.  The intuition is that whereas morphisms in $ \C $ may be represented graphically as diagrams drawn on a bounded region of the plane, morphisms in $ \TC $ are given  by diagrams drawn on a \emph{cylinder}. For example, for any $ f \in \Hom_{\C}(X,Y) $ diagrammatically represented by
	
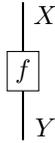
\begin{align*}
		\begin{tikzpicture}[scale=0.15,every node/.style={inner sep=0,outer sep=-1}]
		\node (v2) at (4,6) {};
		\node (v1) at (-4,6) {};
		\node (v4) at (4,-6) {};
		\node (v3) at (-4,-6) {};
		\node (v5) at (0,6) {};
		\node (v7) at (0,-6) {};
		\node at (2,5) {$X$};
		\node at (2,-5) {$Y$};
		\node [draw,outer sep=0,inner sep=3,minimum size=10] (v6) at (0,0) {$ f $};
		\draw[thick]  (v6) edge (v5);
		\draw[thick]  (v7) edge (v6);
		\end{tikzpicture}
	\end{align*}
there will be a morphism in $ \TC $ diagrammatically represented by
	\begin{align} \label{diag:morphism_on_cylinder}
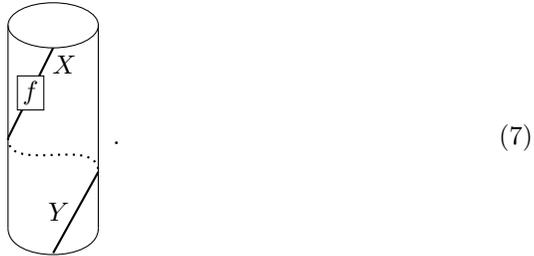
		\begin{array}{c}
			\begin{tikzpicture}[scale=0.15,every node/.style={inner sep=0,outer sep=-0.5}]
			\node (v2) at (4,12) {};
			\node (v1) at (-4,12) {};
			\node (v4) at (4,-6) {};
			\node (v3) at (-4,-6) {};
			\node (v5) at (0,10) {};
			\node (v7) at (-4,2) {};
			\node at (1,8.5) {$X$};
			\node at (0.5,-4.5) {$Y$};
			\node [draw,outer sep=0,inner sep=2,minimum size=10] (v6) at (-2,6) {$ f $};
			\draw[thick]  (v6) edge (v5);
			\draw[thick]  (v7) edge (v6);
			\draw  (0,12) ellipse (4 and 2);
			\draw  (v3) edge (v1);
			\draw  (v2) edge (v4);
			\draw (v3) to[out=-90,in=-90] (v4);
			\node [very thick] (v8) at (4,-1) {};
			\node [very thick] (v9) at (0,-8.2) {};
			\draw[thick,dotted]  (v7) to[in=90,out=-90] (v8);
			\draw[thick]  (v8) edge (v9);
			\end{tikzpicture}
		\end{array}.
	\end{align}
We capture such morphisms by drawing diagrams in a diamond and glueing the upper left and lower right edges. For example morphism \eqref{diag:morphism_on_cylinder} is represented by
	\begin{align*}
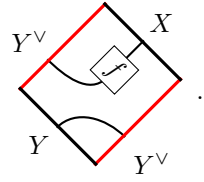
		\begin{array}{c}
			\begin{tikzpicture}[scale=0.5,every node/.style={inner sep=0,outer sep=-1}]
			\node (v1) at (0.75,2.75) {};
			\node (v4) at (0.5,-1.5) {};
			\node (v2) at (2.75,0.75) {};
			\node (v3) at (-1.5,0.5) {};
			\node (v9) at (1.75,1.75) {};
			\node (v6) at (-0.5,-0.5) {};
			\node (v8) at (-0.75,1.25) {};
			\node (v11) at (1.25,-0.75) {};
			\node at (-1,-1) {$ Y $};
			\node at (2.25,2.25) {$ X $};
			\node at (2,-1.5) {$ Y^\vee $};
			\node at (-1.25,1.75) {$ Y^\vee $};
			\node (v7) at (0.75,0.75) {};
			\draw [thick] (v7) to[out=-135,in=-45] (v8);
			\draw [thick] (v11) to[out=135,in=45] (v6);
			\node [draw,outer sep=0,rotate=-45,inner sep=2,minimum size=10,minimum width =13, fill=white] (v5) at (1,1) {\mbox{$ f $}};
			\draw[thick]  (v9) edge (v5);
			\draw[very thick, red]  (v1) edge (v3);
			\draw[very thick, red]  (v2) edge (v4);
			\draw[very thick]  (v1) edge (v2);
			\draw[very thick]  (v3) edge (v4);
			\end{tikzpicture}
		\end{array}.
	\end{align*}
We note that this diagram may also be read vertically and interpreted as an element in $ \Hom_{\C}(Y^\vee X, Y Y^\vee) $.
We also note that due to Lemma~\ref{lem:decompose_object} we may restrict ourselves to only gluing \emph{simple} strands. In this way morphism \eqref{diag:morphism_on_cylinder} would be represented as
	\begin{align*}
	\sum\limits_{R,b}
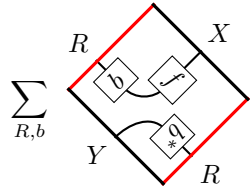
		\begin{array}{c}
				\begin{tikzpicture}[scale=0.5,every node/.style={inner sep=0,outer sep=-1}]
				\node (v1) at (0.5,3) {};
				\node (v4) at (0.75,-1.75) {};
				\node (v2) at (3,0.5) {};
				\node (v3) at (-1.75,0.75) {};
				\node (v9) at (1.75,1.75) {};
				\node (v6) at (-0.5,-0.5) {};
				\node (v8) at (-0.25,0.75) {};
				\node (v11) at (0.75,-0.25) {};
				\node at (-1,-1) {$ Y $};
				\node at (2.25,2.25) {$ X $};
				\node at (2,-1.5) {$ R $};
				\node at (-1.5,2) {$ R $};
				\node (v7) at (0.75,0.75) {};
				\draw [thick] (v7) to[out=-135,in=-45] (v8);
				\draw [thick] (v11) to[out=135,in=45] (v6);
				\node [draw,outer sep=0,rotate=-45,inner sep=2,minimum size=10,minimum width=13,fill=white] (v5) at (1,1) {\mbox{$ f $}};
				\node [draw,outer sep=0,rotate=45,inner sep=2,minimum height=10,minimum width=13,fill=white] (v10) at (-0.5,1) {\mbox{$ b $}};
				\node [draw,outer sep=0,rotate=-135,inner sep=2,minimum width=13,minimum height=10,fill=white] (v13) at (1,-0.5) {\mbox{$ b^* $}};
				\draw[thick]  (v9) edge (v5);
				\node (v12) at (-1,1.5) {};
				\node (v14) at (1.5,-1) {};
				\draw [thick] (v10) edge (v12);
				\draw [thick] (v13) edge (v14);
				\draw[very thick, red]  (v1) edge (v3);
				\draw[very thick, red]  (v2) edge (v4);
				\draw[very thick]  (v1) edge (v2);
				\draw[very thick]  (v3) edge (v4);
				\end{tikzpicture}
		\end{array}
	\end{align*}
where $ b $ ranges over a basis of $ \Hom_{\C}(R,Y^\vee) $. We note that each diagram may now be read vertically as an element in $ \Hom_{\C}(RX,YR) $. With this motivation in mind we may proceed with the definition of $ \TC $.

	\begin{DEF} \label{def:tube_category}

	Let $ \C $ be a spherical fusion category. The associated \emph{tube category}, denoted $ \TC $, is defined as the following category,

		\begin{enumerate}

		\item $ \Obj(\TC) \vcentcolon= \Obj(\C) $

		\item $ \Hom_{\TC}(X,Y) \vcentcolon= \bigoplus_{R} \Hom_{\C}(RX,YR) $

		\item Let $ f $ be in $ \Hom_{\TC}(X,Y)  $ and let $ g $ be in $ \Hom_{\TC}(Y,Z)  $. We define $ g \circ f $ as follows (using the diagrams explained above)
			\begin{align} \label{eq:composition_in_tc}
			g \circ f \defeq \bigoplus\limits_{T} \sum\limits_{S,R,b}
				\begin{array}{c}
					\begin{tikzpicture}[scale=0.5,every node/.style={inner sep=0,outer sep=-1}]
					\node (v1) at (-1.25,2.75) {};
					\node (v4) at (-1.25,-5.25) {};
					\node (v2) at (2.75,-1.25) {};
					\node (v3) at (-5.25,-1.25) {};
					\node (v9) at (0.75,0.75) {};
					\node (v6) at (-3.25,-3.25) {};
					\node (v70) at (0.75,-3.25) {};
					\node (v7) at (-3.25,0.75) {};
					\node at (0.75,-1.25) {$ R $};
					\node at (-3.25,-1) {$ S $};
					\node at (1.25,1.25) {$ X $};
					\node at (-3.75,-3.75) {$ Z $};
					\node at (-1.5607,-1.0119) {$ Y $};
					\node at (1.2,-3.7) {$ T $};
					\node at (-3.75,1.25) {$ T $};
					\node (v11) at (-2.625,-1.875) {};
					\node (v12) at (-1.875,-2.625) {};
					\draw [thick] (v11) edge (v12);
					\node (v14) at (-0.625,0.125) {};
					\node (v13) at (0.125,-0.625) {};
					\draw [thick] (v13) edge (v14);
					\node (v15) at (-2.45,-0.425) {};
					\node (v16) at (-2.075,-0.05) {};
					\draw [thick] (v15) to[out=-45,in=135] (v11);
					\draw [thick] (v16) to[out=-45,in=135] (v14);
					\node (v17) at (-0.05,-2.075) {};
					\node (v18) at (-0.425,-2.45) {};
					\draw [thick] (v13) to[out=-45,in=135] (v17);
					\draw [thick] (v12) to[out=-45,in=135] (v18);
					\node [diamond,draw,outer sep=0,inner sep=0.3,minimum size=25,fill=white] (v50) at (-2.25,-2.25) {\mbox{$ g_{S} $}};
					\node [diamond,draw,outer sep=0,inner sep=-0.2,minimum size=25,fill=white] (v5) at (-0.25,-0.25) {\mbox{$ f_R $}};
					\node [draw,rotate=45,outer sep=0,inner sep=2,minimum height=10,minimum width=13,fill=white] (v8) at (-2.5,0) {$ b $};
					\node [draw,rotate=45,outer sep=0,inner sep=2,minimum height=10,minimum width=13,fill=white] (v10) at (0,-2.5) {$ b^* $};
					\draw[thick]  (v50) edge (v5);
					\draw [thick] (v50) edge (v6);
					\draw[thick]  (v7) edge (v8);
					\draw[thick]  (v10) edge (v70);
					\draw[thick]  (v9) edge (v5);
					\draw[very thick, red]  (v1) edge (v3);
					\draw[very thick, red]  (v2) edge (v4);
					\draw[very thick]  (v1) edge (v2);
					\draw[very thick]  (v3) edge (v4);
					\end{tikzpicture}
				\end{array}
			\end{align}
		where $ f_R $ and $ g_S $ are the $ \Hom_{\C}(RX,YR) $ and $ \Hom_{\C}(SY,ZS) $ components of $ f $ and $ g $ respectively and $ b $ ranges over a basis of $ \Hom_{\C}(T,SR) $. We note that $ g \circ f \in \bigoplus_{T} \Hom_{\C}(TX,ZT) = \Hom_{\TC}(X,Z)  $ as desired.

		\end{enumerate}

	From Lemma~\ref{lem:decompose_object} we see that this definition agrees with the intuition that composition corresponds to vertically stacking the cylinders upon which the diagrams are drawn. This intuition, together with the associativity of the tensor product, makes it clear that composition in $ \TC $ is associative.

		\begin{REM}

		At this point, a careful reader reader might protest that the tensor product is merely weakly associative and yet composition in a category must be strongly associative. However, this is not an issue as the associator isomorphisms will simply modify the basis appearing in~\eqref{eq:composition_in_tc} leaving the composition unchanged.

		\end{REM}

	\end{DEF}

	\begin{REM}
		The summand indexed by $ \tid $ in $ \Hom_{\TC}(X,Y) $ is $ \Hom_{\C}(X,Y) $. This gives a map $ \Hom_{\C}(X,Y) \hookrightarrow \Hom_{\TC}(X,Y) $ such that

		\begin{equation*}
			\begin{tikzcd}[row sep=2cm,column sep=2cm,inner sep=1ex,labels={font=\normalfont}]
			\Hom_{\C}(X,Y) \otimes \Hom_{\C}(Y,Z) \arrow[swap]{d}[name=D]{\circ}   \arrow[hook]{r}{} &  \Hom_{\TC}(X,Y) \otimes \Hom_{\TC}(Y,Z) \arrow{d}[name=U]{\circ}
			\\
			\Hom_{\C}(X,Z) \arrow[hook]{r}{} & \Hom_{\TC}(X,Z)
			\arrow[to path={(U) node[scale=2,midway] {$\circlearrowleft$}  (D)}]{}
			\end{tikzcd}
	\end{equation*}
	commutes. In other words $ \C $ is a \emph{subcategory} of $ \TC $. In particular the identity in $ \End_{\TC}(X) $ is given by the image of $ \id_X \in \End_{\C}(X) $ under this embedding.

	\end{REM}

	\begin{REM}

	If we consider the algebra
		\begin{align*}
		\TA := \End_{\TC}\left(\bigoplus\limits_{S} S \right)
		\end{align*}
	we recover Ocneanu's \emph{tube algebra}~\cite{Ocneanu1993}. As $ \bigoplus\limits_{S} S $ is a projective generator in $ \TC $ the functor
		\begin{align*}
		\RTC &\to \LMod{\TA}\\
		F &\mapsto \Hom_{\RTC}\left(F,\bigoplus\limits_{S} S\right)
		\end{align*}
	gives an equivalence, i.e.\ $ \TC $ is Morita equivalent to $ \TA $.

	\end{REM}

	\begin{REM} \label{rem:end_alg_of_tid}

	The definition of Hom-spaces in $ \TC $ has the following interesting consequence. Let $ \K(\C) $ denote the Grothendieck ring of $ \C $ and let $ \K_{\fld}(\C) $ denote $ \K(\C) \otimes_{\mathbb{Z}} \fld $. Then $ \End_{\TC}(\tid) $ and $ \K_{\fld}(\C) $ are canonically isomorphic algebras. Indeed, $ \End_{\TC}(\tid) = \bigoplus_S \End(S) = \bigoplus_S \fld $ is precisely the underlying vector space of $ \K_{\fld}(\C) $. Furthermore, composition in $ \End_{\TC}(\tid) $ corresponds to the tensor product in $ \K_{\fld}(\C) $ by Lemma~\ref{lem:decompose_object}.

	\end{REM}

For $ X,Y,G $ in $ \C $ and $ \alpha \in \Hom_{\C}(GX,YG) $ we use
	\begin{align*}
	\alpha_G =
		\begin{array}{c}
			\begin{tikzpicture}[scale=0.25,every node/.style={inner sep=0,outer sep=-1}]
			\node (v1) at (0,4) {};
			\node (v4) at (0,-4) {};
			\node (v2) at (4,0) {};
			\node (v3) at (-4,0) {};
			\node (v5) at (-2,2) {};
			\node (v6) at (2,-2) {};
			\node (v7) at (2,2) {};
			\node (v11) at (-2,-2) {};
			\node [draw,diamond,outer sep=0,inner sep=.5,minimum size=22,fill=white] (v9) at (0,0) {$ \alpha $};
			\node at (3,3) {$ X $};
			\node at (-3,-3) {$ Y $};
			\node at (-3,3) {$ G $};
			\node at (3,-3) {$ G $};
			\draw [thick] (v9) edge (v6);
			\draw [thick] (v9) edge (v7);
			\draw [thick] (v9) edge (v11);
			\draw [thick] (v5) edge (v9);
			\draw[very thick, red]  (v1) edge (v3);
			\draw[very thick, red]  (v2) edge (v4);
			\draw[very thick]  (v1) edge (v2);
			\draw[very thick]  (v3) edge (v4);
			\end{tikzpicture}
		\end{array}
	\end{align*}
as shorthand for
	\begin{align*}
	\bigoplus\limits_S \sum\limits_b
		\begin{array}{c}
			\begin{tikzpicture}[scale=0.25,every node/.style={inner sep=0,outer sep=-1}]
			\node (v1) at (-1.5,5.5) {};
			\node (v4) at (1.5,-5.5) {};
			\node (v2) at (5.5,-1.5) {};
			\node (v3) at (-5.5,1.5) {};
			\node (v5) at (-3.5,3.5) {};
			\node (v6) at (3.5,-3.5) {};
			\node (v7) at (2,2) {};
			\node (v11) at (-2,-2) {};
			\node [draw,diamond,outer sep=0,inner sep=.5,minimum size=18,fill=white] (v9) at (0,0) {$ \alpha $};
			\node at (3,3) {$ X $};
			\node at (-3,-3) {$ Y $};
			\node at (-4.5,4.5) {$ S $};
			\node at (4.5,-4.5) {$ S $};
			\draw [thick] (v9) edge (v6);
			\draw [thick] (v9) edge (v7);
			\draw [thick] (v9) edge (v11);
			\draw [thick] (v5) edge (v9);
			\node [draw,rotate=45,outer sep=0,inner sep=2,minimum width=13,fill=white] (v8) at (-2,2) {$ b $};
			\node [draw,rotate=45,outer sep=0,inner sep=2,minimum width=13,fill=white] (v8) at (2,-2) {$ b^* $};
			\draw[very thick, red]  (v1) edge (v3);
			\draw[very thick, red]  (v2) edge (v4);
			\draw[very thick]  (v1) edge (v2);
			\draw[very thick]  (v3) edge (v4);
			\end{tikzpicture}
		\end{array} \in \bigoplus\limits_S \Hom_{\C}(SX,YS) = \Hom_{\TC}(X,Y).
	\end{align*}
	\begin{REM}

	This new notation may potentially create confusion with the pre-existing convention that, for $ f \in \Hom_{\TC}(X,Y) $ and $ R \in \Irr(\C) $, $ f_R $ denotes the $ \Hom_{C}(RX,YR) $ competent of $ f $. To avoid such confusion we will restrict our new notation to the letters $ G $ and $ H $.

	\end{REM}

We note that, using this new notation, we have
	\begin{align} \label{eq:pushing_map_across}
		\begin{array}{c}
			\begin{tikzpicture}[scale=0.25,every node/.style={inner sep=0,outer sep=-1}]
			\node (v1) at (-2,6) {};
			\node (v4) at (0,-4) {};
			\node (v2) at (4,0) {};
			\node (v3) at (-6,2) {};
			\node (v5) at (-4,4) {};
			\node (v6) at (2,-2) {};
			\node (v7) at (2,2) {};
			\node (v11) at (-2,-2) {};
			\node [draw,diamond,outer sep=0,inner sep=.5,minimum size=22,fill=white] (v9) at (0,0) {$ \alpha $};
			\node at (3,3) {$ X $};
			\node at (-3,-3) {$ Y $};
			\node at (-5,5) {$ G_1 $};
			\node at (3,-3) {$ G_1 $};
			\node at (-0.5,2.5) {$ G_2 $};
			\draw [thick] (v9) edge (v6);
			\draw [thick] (v9) edge (v7);
			\draw [thick] (v9) edge (v11);
			\draw [thick] (v5) edge (v9);
			\node [draw,rotate=45,outer sep=0,inner sep=2,minimum width=13,fill=white] (v8) at (-3,3) {$ g $};
			\draw[very thick, red]  (v1) edge (v3);
			\draw[very thick, red]  (v2) edge (v4);
			\draw[very thick]  (v1) edge (v2);
			\draw[very thick]  (v3) edge (v4);
			\end{tikzpicture}
		\end{array}
	=
		\begin{array}{c}
			\begin{tikzpicture}[scale=0.25,every node/.style={inner sep=0,outer sep=-1}]
			\node (v1) at (-2,6) {};
			\node (v4) at (0,-4) {};
			\node (v2) at (4,0) {};
			\node (v3) at (-6,2) {};
			\node (v5) at (-4,4) {};
			\node (v6) at (2,-2) {};
			\node (v7) at (0,4) {};
			\node (v11) at (-4,0) {};
			\node [draw,diamond,outer sep=0,inner sep=.5,minimum size=22,fill=white] (v9) at (-2,2) {$ \alpha $};
			\node at (1,5) {$ X $};
			\node at (-5,-1) {$ Y $};
			\node at (-5,5) {$ G_2 $};
			\node at (3,-3) {$ G_2 $};
			\node at (0.75,1.25) {$ G_1 $};
			\draw [thick] (v9) edge (v6);
			\draw [thick] (v9) edge (v7);
			\draw [thick] (v9) edge (v11);
			\draw [thick] (v5) edge (v9);
			\node [draw,rotate=45,outer sep=0,inner sep=2,minimum width=13,fill=white] (v8) at (1,-1) {$ g $};
			\draw[very thick, red]  (v1) edge (v3);
			\draw[very thick, red]  (v2) edge (v4);
			\draw[very thick]  (v1) edge (v2);
			\draw[very thick]  (v3) edge (v4);
			\end{tikzpicture}
		\end{array}.
	\end{align}
Indeed the $ S $-summand of the left hand side of~\eqref{eq:pushing_map_across} is
	\begin{align*}
	\sum\limits_b
		\begin{array}{c}
			\begin{tikzpicture}[scale=0.25,every node/.style={inner sep=0,outer sep=-1}]
			\node (v1) at (-1,3.5) {};
			\node (v3) at (1,3.5) {};
			\node (v2) at (-1,-7) {};
			\draw [thick] (v1) edge (v2);
			\node (v4) at (1,-7) {};
			\draw [thick] (v3) edge (v4);
			\node [draw,outer sep=0,inner sep=2,minimum width=13,minimum height=11,fill=white] (v8) at (-1,2) {$ b $};
			\node [draw,outer sep=0,inner sep=2,minimum width=13,minimum height=11,fill=white] (v8) at (-1,-0.5) {$ g $};
			\node [draw,outer sep=0,inner sep=2,minimum width=27,minimum height=11,fill=white] (v8) at (0,-3) {$ \alpha $};
			\node [draw,outer sep=0,inner sep=2,minimum width=13,minimum height=11,fill=white] (v8) at (1,-5.5) {$ b^* $};
			\node at (-1,4.5) {$ S $};
			\node at (1,4.5) {$ X $};
			\node at (1,-8) {$ S $};
			\node at (-1,-8) {$ Y $};
			\end{tikzpicture}
		\end{array}
	= \sum\limits_{b,\bar{b}}   \langle \bar{b}^*, g \circ b \rangle
		\begin{array}{c}
			\begin{tikzpicture}[scale=0.25,every node/.style={inner sep=0,outer sep=-1}]
			\node (v1) at (-1,1) {};
			\node (v3) at (1,1) {};
			\node (v2) at (-1,-7) {};
			\draw [thick] (v1) edge (v2);
			\node (v4) at (1,-7) {};
			\draw [thick] (v3) edge (v4);
			\node [draw,outer sep=0,inner sep=2,minimum width=13,minimum height=11,fill=white] (v8) at (-1,-0.5) {$ \bar{b} $};
			\node [draw,outer sep=0,inner sep=2,minimum width=27,minimum height=11,fill=white] (v8) at (0,-3) {$ \alpha $};
			\node [draw,outer sep=0,inner sep=2,minimum width=13,minimum height=11,fill=white] (v8) at (1,-5.5) {$ b^* $};
			\node at (-1,2) {$ S $};
			\node at (1,2) {$ X $};
			\node at (1,-8) {$ S $};
			\node at (-1,-8) {$ Y $};
			\end{tikzpicture}
		\end{array}
	\end{align*}
and similarly the right hand side of~\eqref{eq:pushing_map_across} is
	\begin{align*}
	\sum\limits_{\bar{b}}
		\begin{array}{c}
			\begin{tikzpicture}[scale=0.25,every node/.style={inner sep=0,outer sep=-1}]
			\node (v1) at (-1,1) {};
			\node (v3) at (1,1) {};
			\node (v2) at (-1,-9.5) {};
			\draw [thick] (v1) edge (v2);
			\node (v4) at (1,-9.5) {};
			\draw [thick] (v3) edge (v4);
			\node [draw,outer sep=0,inner sep=2,minimum width=13,minimum height=11,fill=white] (v8) at (-1,-0.5) {$ \bar{b} $};
			\node [draw,outer sep=0,inner sep=2,minimum width=13,minimum height=11,fill=white] (v8) at (1,-5.5) {$ g $};
			\node [draw,outer sep=0,inner sep=2,minimum width=27,minimum height=11,fill=white] (v8) at (0,-3) {$ \alpha $};
			\node [draw,outer sep=0,inner sep=2,minimum width=13,minimum height=11,fill=white] (v8) at (1,-8) {$ \bar{b}^* $};
			\node at (-1,2) {$ S $};
			\node at (1,2) {$ X $};
			\node at (1,-10.5) {$ S $};
			\node at (-1,-10.5) {$ Y $};
			\end{tikzpicture}
		\end{array}
	= \sum\limits_{b,\bar{b}}   \langle \bar{b}^*, g \circ b \rangle
		\begin{array}{c}
			\begin{tikzpicture}[scale=0.25,every node/.style={inner sep=0,outer sep=-1}]
			\node (v1) at (-1,1) {};
			\node (v3) at (1,1) {};
			\node (v2) at (-1,-7) {};
			\draw [thick] (v1) edge (v2);
			\node (v4) at (1,-7) {};
			\draw [thick] (v3) edge (v4);
			\node [draw,outer sep=0,inner sep=2,minimum width=13,minimum height=11,fill=white] (v8) at (-1,-0.5) {$ \bar{b} $};
			\node [draw,outer sep=0,inner sep=2,minimum width=27,minimum height=11,fill=white] (v8) at (0,-3) {$ \alpha $};
			\node [draw,outer sep=0,inner sep=2,minimum width=13,minimum height=11,fill=white] (v8) at (1,-5.5) {$ b^* $};
			\node at (-1,2) {$ S $};
			\node at (1,2) {$ X $};
			\node at (1,-8) {$ S $};
			\node at (-1,-8) {$ Y $};
			\end{tikzpicture}
		\end{array}.
	\end{align*}
This serves as an effective reality-check that Definition~\ref{def:tube_category} captures our original motivation of constructing an annular analogue of $ \C $.

\section{Representations of the Tube Category}

Let $ \C $ be a spherical fusion category let $ \RTC $ be the category of (contravariant) functors from $ \TC $ to $ \Vect $. As $ \C $ is a subcategory of $ \TC $ we have a canonical (covariant) functor
	\begin{align*}
	\RTC &\to \RC \\
	F &\mapsto \bar{F}
	\end{align*}
that simply restricts $ F $ to morphisms in $ \C $. A natural question now arises: for a given object $ \bar{F} $ in $ \RC $ what additional data could be provided to specify a unique extension to an object $ F $ in $ \RTC $? To answer this question we consider the following morphisms in $ \TC $:
	\begin{align*}
	c_{G,X} =
		\begin{array}{c}
		 	\begin{tikzpicture}[scale=0.2,every node/.style={inner sep=0,outer sep=-1}]
		 	\node (v1) at (0,5) {};
		 	\node (v2) at (5,0) {};
		 	\node (v3) at (-5,0) {};
		 	\node (v4) at (0,-5) {};
			\node (v5) at (3.5,1.5) {};
			\node (v6) at (2.5,-2.5) {};
			\node (v7) at (1.5,3.5) {};
			\node (v8) at (-1.5,-3.5) {};
			\node (v9) at (-2.5,2.5) {};
			\node (v10) at (-3.5,-1.5) {};
		 	\node at (2.75,4.75) {$ X $};
		 	\node at (-2.75,-4.75) {$ X $};
		 	\node at (4.75,2.75) {$ G $};
		 	\node at (-4.75,-2.75) {$ G $};
		 	\node at (3.75,-3.75) {$ G $};
		 	\node at (-3.75,3.75) {$ G $};
			\draw[thick]  (v5) to[out=-135,in=135] (v6);
			\draw[thick] (v7) to[out=-135,in=45] (v8);
			\draw[thick]  (v9) to[out=-45,in=45] (v10);
		 	\draw[very thick, red]  (v1) edge (v3);
		 	\draw[very thick, red]  (v2) edge (v4);
		 	\draw[very thick]  (v1) edge (v2);
		 	\draw[very thick]  (v3) edge (v4);
		 	\end{tikzpicture}
		\end{array}
	\quad \text{where $ G $ and $ X $ are in $ \C $.}
	\end{align*}
For $ f $ and $ g $ in $ \Hom_{\C}(X,Y) $ and $ \Hom_{\C}(G_1,G_2) $ respectively, we have
	\begin{align} \label{eq:naturality_of_cycle}
		\begin{split}
		(g \otimes f) \circ c_{G_1,X} &=
			\begin{array}{c}
				\begin{tikzpicture}[scale=0.2,every node/.style={inner sep=0,outer sep=-1}]
				\node (v1) at (0,5) {};
				\node (v2) at (5,0) {};
				\node (v3) at (-7,-2) {};
				\node (v4) at (-2,-7) {};
				\node (v5) at (3.5,1.5) {};
				\node (v6) at (3,-2) {};
				\node (v7) at (-2.5,2.5) {};
				\node (v8) at (-3.5,-1.5) {};
				\node (v9) at (1.5,3.5) {};
				\node (v10) at (-1.5,-3.5) {};
				\node at (2.75,4.75) {$ X $};
				\node at (4.75,2.75) {$ G_1 $};
				\node at (-6.75,-4.75) {$ G_2 $};
				\node at (-4.75,-6.75) {$ Y $};
				\node at (4.25,-3.25) {$ G_1 $};
				\node at (-3.75,3.75) {$ G_1 $};
				\draw[thick]  (v5) to[out=-135,in=135] (v6);
				\draw[thick] (v7) to[out=-45,in=45] (v8);
				\draw[thick]  (v9) to[out=-135,in=45] (v10);
				\node [draw,rotate=-45,outer sep=0,inner sep=2,minimum width=11,minimum height=13,fill=white] (v80) at (-4,-2) {$ g $};
				\node [draw,rotate=-45,outer sep=0,inner sep=2,minimum width=11,minimum height=13,fill=white] (v800) at (-2,-4) {$ f $};
				\draw[very thick, red]  (v1) edge (v3);
				\draw[very thick, red]  (v2) edge (v4);
				\draw[very thick]  (v1) edge (v2);
				\draw [thick] (v3) edge (v4);
				\node (v11) at (-5.5,-3.5) {};
				\node (v12) at (-3.5,-5.5) {};
				\draw [thick] (v80) edge (v11);
				\draw [thick] (v800) edge (v12);
				\end{tikzpicture}
			\end{array} \\
		&=
			\begin{array}{c}
				\begin{tikzpicture}[scale=0.2,every node/.style={inner sep=0,outer sep=-1},xscale=-1,yscale=-1]
				\node (v1) at (0,5) {};
				\node (v2) at (5,0) {};
				\node (v3) at (-7,-2) {};
				\node (v4) at (-2,-7) {};
				\node (v5) at (-3.5,-1.5) {};
				\node (v6) at (-2.5,2.5) {};
				\node (v7) at (3.5,1.5) {};
				\node (v8) at (2.5,-2.5) {};
				\node (v9) at (1.5,3.5) {};
				\node (v10) at (-1.5,-3.5) {};
				\node at (2.75,4.75) {$ Y $};
				\node at (4.75,2.75) {$ G_2 $};
				\node at (-6.75,-4.75) {$ G_1 $};
				\node at (-4.75,-6.75) {$ X $};
				\node at (3.75,-3.75) {$ G_2 $};
				\node at (-3.75,3.75) {$ G_2 $};
				\draw[thick]  (v5) to[out=45,in=-45] (v6);
				\draw[thick] (v7) to[out=-135,in=135] (v8);
				\draw[thick]  (v9) to[out=-135,in=45] (v10);
				\node [draw,rotate=-45,outer sep=0,inner sep=2,minimum width=11,minimum height=13,fill=white] (v80) at (-4,-2) {$ g $};
				\node [draw,rotate=-45,outer sep=0,inner sep=2,minimum width=11,minimum height=13,fill=white] (v800) at (-2,-4) {$ f $};
				\draw[very thick, red]  (v1) edge (v3);
				\draw[very thick, red]  (v2) edge (v4);
				\draw[very thick]  (v1) edge (v2);
				\draw [thick] (v3) edge (v4);
				\node (v11) at (-5.5,-3.5) {};
				\node (v12) at (-3.5,-5.5) {};
				\draw [thick] (v80) edge (v11);
				\draw [thick] (v800) edge (v12);
				\end{tikzpicture}
			\end{array}
		= c_{G_2,Y} \circ (f \otimes g).
		\end{split}
	\end{align}
Furthermore, $ c_{G,X} $ is an isomorphism and satisfies $ c_{G,HX} \circ c_{H,XG} = c_{GH,X}  $. The principal aim of this section is to prove that specifying how to evaluate $ F $ on $ c_{G,X} $ precisely captures how to extend $ \bar{F} $ to $ F $ uniquely and therefore provides an answer to the question mentioned above. Applying $ F $ to $ c_{G,X} $ gives a collection of maps
	\begin{align*}
	\kappa_{G,X}\colon \bar{F}(G X) \to \bar{F}(X G).
	\end{align*}
By~\eqref{eq:naturality_of_cycle} $ \kappa $ is natural in both $ X $ and $ G $. The additional properties of $ c_{G,X} $ listed above then imply that $ \kappa_{G,X} $ is an isomorphism and satisfies $ \kappa_{H,XG} \circ \kappa_{G,HX} = \kappa_{GH,X}  $. Suppose we start with an arbitrary object $ \bar{F} $ in $ \RC $ and isomorphisms $ \kappa_{G,X}\colon \bar{F}(G X) \to \bar{F}(X G) $ that satisfy naturality in $ G $ and $ X $ and $ \kappa_{H,XG} \circ \kappa_{G,HX} = \kappa_{GH,X}  $. We shall prove that there is a unique functor $ F $ in $ \RTC $ such that

	\begin{enumerate}[(i)]
	\item $ F(X) = \bar{F}(X) $ for all $ X $ in $ \C $ \label{ennum:objects}
	\item $ F(\alpha) = \bar{F}(\alpha) $ for all $ \alpha \in \Hom_{\C}(X,Y) $ \label{ennum:morphisms}
	\item $ F(c_{G,X}) = \kappa_{G,X} $ for all $ G,X $ in $ \C $. \label{ennum:twists}
	\end{enumerate}

If such a functor exists it is certainly unique as these conditions determine the functor on any morphism in $ \TC $. Indeed for any $ \alpha_G \in \Hom_{\TC}(Y,X) $ we have
	\begin{align} \label{eq:decomposing_morphism_in_tc}
	\alpha_G =
		\begin{array}{c}
			\begin{tikzpicture}[scale=0.2,every node/.style={inner sep=0,outer sep=-1}]
			\node (v1) at (0,5) {};
			\node (v4) at (0,-5) {};
			\node (v2) at (5,0) {};
			\node (v3) at (-5,0) {};
			\node (v5) at (1,4) {};
			\node (v7) at (4,1) {};
			\node (v8) at (-3.5,-1.5) {};
			\node (v10) at (-1.5,-3.5) {};
			\node (v6) at (-2.5,2.5) {};
			\node (v9) at (2.5,-2.5) {};
			\node (v12) at (-2.5,-2.5) {};
			\node (v14) at (2.5,2.5) {};
			\node at (2,5) {$ G $};
			\node at (5.5,2) {$ G^\vee $};
			\node at (-3.75,-3.75) {$ X $};
			\node at (3.75,3.25) {$ Y $};
			\node (v11) at (-1.0592,0.0122) {};
			\node (v15) at (1.0532,0.006) {};
			\node (v16) at (0,-1) {};
			\node (v13) at (-0.0183,1.0531) {};
			\draw[thick]  (v11) to[out=180,in=45] (v12);
			\draw [thick] (v13) to[out=90,in=-135] (v5);
			\draw [thick] (v14) to[out=-135,in=0] (v15);
			\draw [thick] plot[smooth, tension=.7] coordinates {(v16) (0.5,-2) (2,-1.5) (v7)};
			\node [draw,diamond,rotate=-45,outer sep=0,inner sep=1.8,minimum size=10,fill=white] at (0,0) {\mbox{$ \alpha $}};
			\draw[very thick, red]  (v1) edge (v3);
			\draw[very thick, red]  (v2) edge (v4);
			\draw[very thick]  (v1) edge (v2);
			\draw[very thick]  (v3) edge (v4);
			\end{tikzpicture}
		\end{array}
	\circ
		\begin{array}{c}
			\begin{tikzpicture}[scale=0.2,every node/.style={inner sep=0,outer sep=-1}]
			\node (v1) at (0,5) {};
			\node (v4) at (0,-5) {};
			\node (v2) at (5,0) {};
			\node (v3) at (-5,0) {};
			\node (v5) at (3.5,1.5) {};
			\node (v7) at (-2.5,2.5) {};
			\node (v8) at (-3.5,-1.5) {};
			\node (v10) at (-1.5,-3.5) {};
			\node (v6) at (2.5,-2.5) {};
			\node (v9) at (2.5,2.5) {};
			\node (v11) at (-2.5,-2.5) {};
			\node (v12) at (1.5,3.5) {};
			\node at (-2.75,-4.75) {$ G^\vee $};
			\node at (4.75,2.5) {$ G $};
			\node at (-5.25,-2.75) {$ G $};
			\node at (2.75,4.5) {$ Y $};
			\draw[thick]  (v5) to[out=-135,in=135] (v6);
			\draw[thick] (v7) to[out=-45,in=45] (v8);
			\draw[thick]  (v9) to[out=-135,in=45] (v10);
			\draw [thick] (v11) to[out=45,in=-135] (v12);
			\draw[very thick, red]  (v1) edge (v3);
			\draw[very thick, red]  (v2) edge (v4);
			\draw[very thick]  (v1) edge (v2);
			\draw[very thick]  (v3) edge (v4);
			\end{tikzpicture}
		\end{array}
	\circ
		\begin{array}{c}
			\begin{tikzpicture}[scale=0.2,every node/.style={inner sep=0,outer sep=-1}]
			\node (v1) at (-3.5,1.5) {};
			\node (v4) at (0,-5) {};
			\node (v2) at (1.5,-3.5) {};
			\node (v3) at (-5,0) {};
			\node (v5) at (-2.5,0.5) {};
			\node (v7) at (-1,-1) {};
			\node (v8) at (-2.5,-2.5) {};
			\node (v10) at (-4,-1) {};
			\node (v11) at (-1,-4) {};
			\node (v12) at (-2.5,0.5) {};
			\node at (-2.25,-5.25) {$ G $};
			\node at (-3.75,-3.75) {$ G^\vee $};
			\node at (-1.25,1.75) {$ Y $};
			\draw [thick] (v8) to[out=45,in=45] (v11);
			\draw [thick] (v10) edge (v12);
			\draw[very thick, red]  (v1) edge (v3);
			\draw[very thick, red]  (v2) edge (v4);
			\draw[very thick]  (v1) edge (v2);
			\draw[very thick]  (v3) edge (v4);
			\end{tikzpicture}
		\end{array}.
	\end{align}
The first and last terms are morphisms in $ \C $ and therefore determined by Condition \eqref{ennum:morphisms}. The middle term is $ c_{G,YG^\vee} $ and is therefore determined by Condition \eqref{ennum:twists}. The following proposition establishes existence.

	\begin{PROP} \label{prop:extending_functor}

	There exists a unque object $ F $ in $ \RTC $ that satisfies Conditions \eqref{ennum:objects}, \eqref{ennum:morphisms} and \eqref{ennum:twists}.

	\proof

	We first check that Conditions \eqref{ennum:morphisms} and \eqref{ennum:twists} don't contradict one another. The only case when $ c_{G,X} $ is a morphism in $ \C $ is when $ G = \tid $ and $ c_{G,X} = \id_X $. As we have
	  	\begin{align*}
	  	\kappa_{\tid,X} \circ \kappa_{\tid,X} = \kappa_{\tid,X}.
	  	\end{align*}
	and $ \kappa_{\tid,X} $ is an isomorphism this implies $ \kappa_{\tid,X} = \id_X $. Therefore Conditions \eqref{ennum:morphisms} and \eqref{ennum:twists} are equivalent in this case.

	To aid legibility when the domain of a map of the form $ \kappa_{G,X} $ is clear from the context we will suppress the second argument and simply write $ \kappa_G $. As discussed before the proof, any $ F $ that satisfies Conditions \eqref{ennum:objects}, \eqref{ennum:morphisms} and \eqref{ennum:twists} also satisfies
		\begin{align} \label{eq:extention}
		F(\alpha_G) = \bar{F} \left(
			\begin{array}{c}
				\begin{tikzpicture}[scale=0.25,xscale=-1]
				\draw [thick] (-1,0) to[out=90,in=90] (1,0);
				\draw [thick] (3,0) -- (3,2);
				\node at (-1,-1) {$G$};
				\node at (1,-1) {$G^\vee$};
				\node at (3,-1) {$Y$};
				\end{tikzpicture}
			\end{array}
		\right) \circ \kappa_{G} \circ \bar{F} \left(
			\begin{array}{c}
				\begin{tikzpicture}[scale=0.15,every node/.style={inner sep=0,outer sep=-1}]
				\node (v5) at (1.75,3.75) {};
				\node (v50) at (-1.75,3.75) {};
				\node (v7) at (1.75,-1) {};
				\node (v70) at (-1.75,-2) {};
				\node at (1.75,5.25) {$Y$};
				\node at (-1.75,5.25) {$G$};
				\node at (-1.75,-3.5) {$X$};
				\node at (5.25,5.25) {$G^\vee$};
				\draw[thick]  (1.75,-0.5) edge (v5);
				\draw[thick]  (-1.75,0.25) edge (v50);
				\draw[thick]  (-1.75,0.25) edge (v70);
				\draw[thick]  (1.75,-0.5) edge (v7);
				\node (v1) at (5,3.75) {};
				\node (v2) at (5,-1) {};
				\draw [thick] (v1) edge (v2);
				\draw [thick] (v7) to[out=-90,in=-90] (v2);
				\node [draw,diamond,outer sep=0,inner sep=3,minimum size=10,fill=white] (v6) at (0,0.75) {$ \alpha $};
				\end{tikzpicture}
			\end{array}
		\right)
		\end{align}
	for any $ \alpha_G \in \Hom_{\TC}(Y,X) $. We therefore only have to check that \eqref{eq:extention} does indeed define a functor. Let $ \beta_H $ be in $ \Hom_{\TC}(Z,Y) $. We have
		\begin{align*}
		&F(\beta_H) \circ F(\alpha_G) = \\ &\bar{F} \left(
			\begin{array}{c}
				\begin{tikzpicture}[scale=0.25]
				\draw [thick] (5,0) to[out=90,in=90] (7,0);
				\draw [thick] (3,0) -- (3,2);
				\node at (5,-1) {$H^\vee$};
				\node at (7,-1) {$H$};
				\node at (3,-1) {$Z$};
				\end{tikzpicture}
			\end{array}
		\right) \circ \kappa_{H} \circ \bar{F} \left(
			\begin{array}{c}
				\begin{tikzpicture}[scale=0.15,every node/.style={inner sep=0,outer sep=-1}]
				\node (v5) at (1.75,4) {};
				\node (v50) at (-1.75,4) {};
				\node (v7) at (1.75,-4) {};
				\node (v70) at (-1.75,-5) {};
				\node (v1) at (11.5,-5) {};
				\node (v4) at (8,-5) {};
				\node (v2) at (11.5,3) {};
				\node (v3) at (8,3) {};
				\node at (2,5.5) {$Z$};
				\node at (-1.875,5.5) {$H$};
				\node at (-2,-6.5) {$Y$};
				\node at (8,-6.5) {$G^\vee$};
				\node at (11.5,-6.5) {$G$};
				\node at (5.5,5.5) {$H^\vee$};
				\draw[thick]  (1.75,-0.75) edge (v5);
				\draw[thick]  (-1.75,-1.75) edge (v50);
				\draw[thick]  (-1.75,-1.75) edge (v70);
				\draw[thick]  (1.75,-0.75) edge (v7);
				\draw [thick] (5,4) edge (5,-4);
				\draw [thick] (v7) to[out=-90,in=-90] (5,-4);
				\draw[thick]  (v1) edge (v2);
				\draw[thick]  (v3) edge (v4);
				\draw[thick]  (v3) to[out=90,in=90] (v2);
				\node [draw,diamond,outer sep=0,inner sep=3,minimum size=10,fill=white] (v6) at (0,0) {$ \beta $};
				\end{tikzpicture}
			\end{array}
		\right)
		\circ \kappa_{G} \circ \bar{F} \left(
			\begin{array}{c}
				\begin{tikzpicture}[scale=0.15,every node/.style={inner sep=0,outer sep=-1}]
				\node (v5) at (1.75,4.25) {};
				\node (v50) at (-1.75,4.25) {};
				\node (v7) at (1.75,-2.25) {};
				\node (v70) at (-1.75,-4.25) {};
				\node at (1.75,5.75) {$Y$};
				\node at (-1.75,5.75) {$G$};
				\node at (-1.75,-5.75) {$X$};
				\node at (5.25,5.75) {$G^\vee$};
				\draw[thick]  (1.75,-0.5) edge (v5);
				\draw[thick]  (-1.75,0.25) edge (v50);
				\draw[thick]  (-1.75,0.25) edge (v70);
				\draw[thick]  (1.75,-0.5) edge (v7);
				\node (v1) at (5,4.25) {};
				\node (v2) at (5,-2.25) {};
				\draw [thick] (v1) edge (v2);
				\draw [thick] (v7) to[out=-90,in=-90] (v2);
				\node [draw,diamond,outer sep=0,inner sep=3,minimum size=10,fill=white] (v6) at (0,0.75) {$ \alpha $};
				\end{tikzpicture}
			\end{array}
		\right).
		\end{align*}
	By the naturality of $ \kappa_{H^\vee} $ and $ \kappa_{G^\vee}  $ this equation may be rearranged. The creation morphism of $ G $ in the middle term may be moved to after $ \kappa_{H} $, giving
		\begin{align*}
		\bar{F} \left(
			\begin{array}{c}
				\begin{tikzpicture}[scale=0.25]
				\draw [thick] (9,0) to[out=90,in=90] (11,0);
				\draw [thick] (7,0) to[out=90,in=90] (13,0);
				\draw [thick] (5,0) -- (5,2);
				\node at (7,-1) {$H^\vee$};
				\node at (13,-1) {$H$};
				\node at (9,-1) {$G^\vee$};
				\node at (11,-1) {$G$};
				\node at (5,-1) {$Z$};
				\end{tikzpicture}
			\end{array}
		\right) \circ \kappa_{H} \circ \bar{F} \left(
			\begin{array}{c}
				\begin{tikzpicture}[scale=0.15,every node/.style={inner sep=0,outer sep=-1}]
				\node (v5) at (1.75,4) {};
				\node (v50) at (-1.75,4) {};
				\node (v7) at (1.75,-4) {};
				\node (v70) at (-1.75,-5) {};
				\node (v1) at (11.5,-5) {};
				\node (v4) at (8,-5) {};
				\node (v2) at (11.5,4) {};
				\node (v3) at (8,4) {};
				\node at (2,5.5) {$Z$};
				\node at (-1.875,5.5) {$H$};
				\node at (-2,-6.5) {$Y$};
				\node at (8,-6.5) {$G^\vee$};
				\node at (11.5,-6.5) {$G$};
				\node at (5.5,5.5) {$H^\vee$};
				\draw[thick]  (1.75,-0.75) edge (v5);
				\draw[thick]  (-1.75,-1.75) edge (v50);
				\draw[thick]  (-1.75,-1.75) edge (v70);
				\draw[thick]  (1.75,-0.75) edge (v7);
				\draw [thick] (5,4) edge (5,-4);
				\draw [thick] (v7) to[out=-90,in=-90] (5,-4);
				\draw[thick]  (v1) edge (v2);
				\draw[thick]  (v3) edge (v4);
				\node [draw,diamond,outer sep=0,inner sep=3,minimum size=10,fill=white] (v6) at (0,0) {$ \beta $};
				\end{tikzpicture}
			\end{array}
		\right)
		\circ \kappa_{G} \circ \bar{F} \left(
			\begin{array}{c}
				\begin{tikzpicture}[scale=0.15,every node/.style={inner sep=0,outer sep=-1}]
				\node (v5) at (1.75,4.25) {};
				\node (v50) at (-1.75,4.25) {};
				\node (v7) at (1.75,-2.25) {};
				\node (v70) at (-1.75,-4.25) {};
				\node at (1.75,5.75) {$Y$};
				\node at (-1.75,5.75) {$G$};
				\node at (-1.75,-5.75) {$X$};
				\node at (5.25,5.75) {$G^\vee$};
				\draw[thick]  (1.75,-0.5) edge (v5);
				\draw[thick]  (-1.75,0.25) edge (v50);
				\draw[thick]  (-1.75,0.25) edge (v70);
				\draw[thick]  (1.75,-0.5) edge (v7);
				\node (v1) at (5,4.25) {};
				\node (v2) at (5,-2.25) {};
				\draw [thick] (v1) edge (v2);
				\draw [thick] (v7) to[out=-90,in=-90] (v2);
				\node [draw,diamond,outer sep=0,inner sep=3,minimum size=10,fill=white] (v6) at (0,0.75) {$ \alpha $};
				\end{tikzpicture}
			\end{array}
		\right).
		\end{align*}
	Then $ \beta $ together with the annihilation morphism of $ H $ may be moved to before $ \kappa_{G^\vee} $. This yields
		\begin{align*}
		\bar{F} \left(
			\begin{array}{c}
				\begin{tikzpicture}[scale=0.25]
				\draw [thick] (9,0) to[out=90,in=90] (11,0);
				\draw [thick] (7,0) to[out=90,in=90] (13,0);
				\draw [thick] (5,0) -- (5,2);
				\node at (7,-1) {$H^\vee$};
				\node at (13,-1) {$H$};
				\node at (9,-1) {$G^\vee$};
				\node at (11,-1) {$G$};
				\node at (5,-1) {$Z$};
				\end{tikzpicture}
			\end{array}
		\right) \circ \kappa_{H} \circ \kappa_{G} \circ \bar{F} \left(
			\begin{array}{c}
				\begin{tikzpicture}[scale=0.15,every node/.style={inner sep=0,outer sep=-1}]
				\node (v1) at (8.5,-5) {};
				\node (v4) at (1.5,4) {};
				\node (v2) at (8.5,4) {};
				\node (v3) at (1.5,-1) {};
				\node (v50) at (-4.5,4) {};
				\node (v500) at (-1.5,4) {};
				\node (v60) at (-4.5,-4.25) {};
				\node (v600) at (-4.5,-4.25) {};
				\node (v7) at (1.5,-5) {};
				\node (v70) at (-4.5,-8) {};
				\node (v10) at (5.5,4) {};
				\node (v20) at (5.5,-5) {};
				\node at (9.5,6) {$G^\vee$};
				\node at (1.5,6) {$Z$};
				\node at (-1.5,6) {$H$};
				\node at (-4.5,-10) {$X$};
				\node at (-4.5,6) {$G$};
				\node at (5.5,6) {$H^\vee$};
				\draw[thick]  (v60) edge (v50);
				\draw[thick]  (v600) edge (v70);
				\draw [thick] (v10) edge (v20);
				\draw [thick] (v7) to[out=-90,in=-90] (v20);
				\draw[thick]  (-1.5,0) edge (v500);
				\draw[thick]  (1.5,-1) edge (v7);
				\draw [thick] (-1.5,-5) node (v5) {} edge (-1.5,-3);
				\draw[thick]  (v1) edge (v2);
				\draw[thick]  (v3) edge (v4);
				\draw [thick] (v5) to[out=-90,in=-90] (v1);
				\node [draw,diamond,outer sep=0,inner sep=3.2,minimum size=10,fill=white] at (-3,-3) {$ \alpha $};
				\node [draw,diamond,outer sep=0,inner sep=1.9,minimum size=10,fill=white] at (0,0) {$ \beta $};
				\end{tikzpicture}
			\end{array}
		\right)
	\end{align*}
	\begin{align*}
		= \quad &\bar{F} \left(
			\begin{array}{c}
				\begin{tikzpicture}[scale=0.25]
				\draw [thick] (9,0) to[out=90,in=90] (11,0);
				\draw [thick] (7,0) to[out=90,in=90] (13,0);
				\draw [thick] (5,0) -- (5,2);
				\node at (7,-1) {$H^\vee$};
				\node at (13,-1) {$H$};
				\node at (9,-1) {$G^\vee$};
				\node at (11,-1) {$G$};
				\node at (5,-1) {$Z$};
				\end{tikzpicture}
			\end{array}
		\right) \circ \kappa_{GH} \circ \bar{F} \left(
			\begin{array}{c}
				\begin{tikzpicture}[scale=0.15,every node/.style={inner sep=0,outer sep=-1}]
				\node (v1) at (8.5,-5) {};
				\node (v4) at (1.5,4) {};
				\node (v2) at (8.5,4) {};
				\node (v3) at (1.5,-1) {};
				\node (v50) at (-4.5,4) {};
				\node (v500) at (-1.5,4) {};
				\node (v60) at (-4.5,-4.25) {};
				\node (v600) at (-4.5,-4.25) {};
				\node (v7) at (1.5,-5) {};
				\node (v70) at (-4.5,-8) {};
				\node (v10) at (5.5,4) {};
				\node (v20) at (5.5,-5) {};
				\node at (9.5,6) {$G^\vee$};
				\node at (1.5,6) {$Z$};
				\node at (-1.5,6) {$H$};
				\node at (-4.5,-10) {$X$};
				\node at (-4.5,6) {$G$};
				\node at (5.5,6) {$H^\vee$};
				\draw[thick]  (v60) edge (v50);
				\draw[thick]  (v600) edge (v70);
				\draw [thick] (v10) edge (v20);
				\draw [thick] (v7) to[out=-90,in=-90] (v20);
				\draw[thick]  (-1.5,0) edge (v500);
				\draw[thick]  (1.5,-1) edge (v7);
				\draw [thick] (-1.5,-5) node (v5) {} edge (-1.5,-3);
				\draw[thick]  (v1) edge (v2);
				\draw[thick]  (v3) edge (v4);
				\draw [thick] (v5) to[out=-90,in=-90] (v1);
				\node [draw,diamond,outer sep=0,inner sep=3.2,minimum size=10,fill=white] at (-3,-3) {$ \alpha $};
				\node [draw,diamond,outer sep=0,inner sep=1.9,minimum size=10,fill=white] at (0,0) {$ \beta $};
				\end{tikzpicture}
			\end{array}
		\right) \\
		= \quad &F(\alpha_G \circ \beta_H).
		\end{align*}
	\endproof

	\end{PROP}

	\begin{DEF}
	We call the functor constructed in Proposition \ref{prop:extending_functor} the \textit{extension of $ \bar{F} $ by $ \kappa $} and denote it $ (\bar{F},\kappa) $.
	\end{DEF}

	This new description of objects in $ \RTC $ then also yields a new description of morphisms as follows.

	\begin{PROP} \label{prop:morphisms_in_rtc}

	We consider $ F = (\bar{F}, \kappa) $ and $ F' = (\bar{F'}, \kappa') $ in $ \RTC $ then we have
		\begin{align*}
		\Hom_{\RTC}(F,F') = \{ \alpha \in \Hom_{\RC}(\bar{F},\bar{F'}) \mid \alpha_{XG} \circ \kappa_{G,X} = \kappa'_{G,X} \circ \alpha_{GX} \}.
		\end{align*}
	\proof

	Let $ \alpha \in \Hom_{\RC}(\bar{F},\bar{F'}) $ be such that $ \alpha_{XG} \circ \kappa_{G,X} = \kappa'_{G,X} \circ \alpha_{GX} $. As $ \alpha $ is in $ \Hom_{\RC}(\bar{F},\bar{F'}) $, $ \alpha $ is natural with respect to all morphism in $ \C $. Furthermore, the additional condition on $ \alpha $ implies that it is also natural with respect to all morphisms of the form $ c_{G,X} $. From~\eqref{eq:decomposing_morphism_in_tc} we see that any morphism in $ \TC $ may be written as a composition of morphisms in $ C $ and morphisms of the form $ c_{G,X} $. Therefore $ \alpha $ is natural with respect to all morphisms in $ \TC $.	\endproof

	\end{PROP}

\section{Equivalence with the Drinfeld Centre} \label{sec:equiv_with_zc}

	\begin{DEF}

	Let $ \C $ be a monoidal category and let $ X $ be an object in $ \C $. A \emph{half-braiding} on $ X $ is a collection of natural isomorphisms
		\begin{align*}
		\tau_G \colon G \otimes X \to X \otimes G
		\end{align*}
	such that
		\begin{align} \label{eq:half_braiding}
		\tau_{GH} = (\tau_G \otimes \id_H) \circ (\id_G \otimes \tau_H)
		\end{align}
	for all $ G,H $ in $ \C $. From a graphical perspective the condition that $ \tau_G $ is natural allows us to `push' morphisms though $ \tau $:
    \begin{align*}
        \begin{array}{c}
            \begin{tikzpicture}[scale = 0.5,every node/.style={inner sep=0,outer sep=-1}]
            \draw[thick] (-0.5,3) -- (-0.5,1.75);
            \draw[thick] (0.625,3) -- (0.625,1.625);
            \draw[thick] (0.625,1.625) to[out=-90,in=90] (0.25,0.625);
            \draw[thick] (-0.5,1.625) to[out=-90,in=90] (-0.25,0.625);
            \node at (-0.5,3.625) {$ G $};
            \node at (0.625,3.6) {$ X $};
            \node at (0.5,-1.5) {$ H $};
            \node at (-0.5,-1.5) {$ X $};
            \node (v1) at (-0.25,-0.25) {};
            \node (v3) at (0.25,-0.25) {};
            \node (v2) at (-0.5,-1) {};
            \node (v4) at (0.5,-1) {};
            \draw [thick] (v1) to[out=-90,in=90] (v2);
            \draw [thick] (v3) to[out=-90,in=90] (v4);
            \node[draw,outer sep=0,inner sep=2,minimum size=15,fill=white] at (-0.5,2.125) {$ \alpha $};
            \node[draw,outer sep=0,inner sep=2,minimum size=16,fill=white] at (0,0.25) {$ \tau_H $};
            \end{tikzpicture}
        \end{array}
    =
        \begin{array}{c}
			\begin{tikzpicture}[scale = 0.5,every node/.style={inner sep=0,outer sep=-1}]
			\draw[thick] (0.625,3.125) to[out=-90,in=90] (0.25,2.125);
			\draw[thick] (-0.5,3.125) to[out=-90,in=90] (-0.25,2.125);
			\node at (-0.5,3.625) {$ G $};
			\node at (0.625,3.6) {$ X $};
			\node at (0.5,-1.5) {$ H $};
			\node at (-0.5,-1.5) {$ X $};
			\node (v1) at (-0.25,1.25) {};
			\node (v3) at (0.25,1.25) {};
			\node (v2) at (-0.5,0.5) {};
			\node (v4) at (0.5,0.5) {};
			\draw [thick] (v1) to[out=-90,in=90] (v2);
			\draw [thick] (v3) to[out=-90,in=90] (v4);
			\node[draw,outer sep=0,inner sep=2,minimum size=15,fill=white] (v6) at (0.5,0) {$ \alpha $};
			\node[draw,outer sep=0,inner sep=2,minimum size=16,fill=white] at (0,1.75) {$ \tau_G $};
			\node (v7) at (0.5,-1) {};
			\node (v5) at (-0.5,-1) {};
			\draw [thick] (v2) edge (v5);
			\draw [thick] (v6) edge (v7);
			\end{tikzpicture}
        \end{array},
    \end{align*}
this motivates the name `half-braiding'. Pushing the graphical perspective further, we may rewrite~\eqref{eq:half_braiding} as
    \begin{align*}
        \begin{array}{c}
			\begin{tikzpicture}[scale = 0.5,every node/.style={inner sep=0,outer sep=-1}]
			\draw[thick] (-0.375,3) to[out=-90,in=90] (0,1.5);
			\draw[thick] (1.375,3) to[out=-90,in=90] (1,1.5);
			\node at (-0.375,3.625) {$ G $};
			\node at (0.5,3.625) {$ H $};
			\node at (1.375,3.625) {$ X $};
			\node at (1.375,-1.5) {$ H $};
			\node at (0.5,-1.5) {$ G $};
			\node at (-0.375,-1.5) {$ X $};
			\node (v1) at (0,0.625) {};
			\node (v3) at (0.5,0.625) {};
			\node (v2) at (-0.375,-1) {};
			\node (v4) at (0.5,-1) {};
			\draw [thick] (v1) to[out=-90,in=90] (v2);
			\draw [thick] (v3) to[out=-90,in=90] (v4);
			\node (v6) at (0.5,1.5) {};
			\node (v5) at (0.5,3) {};
			\draw [thick] (v5) to[out=-90,in=90] (v6);
			\node (v7) at (1,0.625) {};
			\node (v8) at (1.375,-1) {};
			\draw [thick] (v7) to[out=-90,in=90] (v8);
			\node[draw,outer sep=0,inner sep=2,minimum size=16,fill=white] at (0.5,1.125) {$ \tau_{GH} $};
			\end{tikzpicture}
        \end{array}
    =
        \begin{array}{c}
    		\begin{tikzpicture}[scale = 0.5,every node/.style={inner sep=0,outer sep=-1}]
    		\draw[thick] (-0.375,3) -- (-0.375,1.75);
    		\draw[thick] (1.375,3) -- (1.375,2.25);
    		\draw[thick] (0.875,1.5) to[out=-90,in=90] (0.5,0.625);
    		\draw[thick] (-0.375,1.75) to[out=-90,in=90] (0,0.625);
    		\node at (-0.375,3.625) {$ G $};
    		\node at (0.5,3.625) {$ H $};
    		\node at (1.375,3.625) {$ X $};
    		\node at (1.375,-1.5) {$ H $};
    		\node at (0.5,-1.5) {$ G $};
    		\node at (-0.375,-1.5) {$ X $};
    		\node (v1) at (0,-0.25) {};
    		\node (v3) at (0.5,-0.25) {};
    		\node (v2) at (-0.375,-1) {};
    		\node (v4) at (0.5,-1) {};
    		\draw [thick] (v1) to[out=-90,in=90] (v2);
    		\draw [thick] (v3) to[out=-90,in=90] (v4);
    		\node (v6) at (0.875,2.25) {};
    		\node (v5) at (0.5,3) {};
    		\draw [thick] (v5) to[out=-90,in=90] (v6);
    		\node (v7) at (1.375,1.5) {};
    		\node (v8) at (1.375,-1) {};
    		\draw [thick] (v7) edge (v8);
    		\node[draw,outer sep=0,inner sep=2,minimum size=15,fill=white] at (1.125,1.875) {$ \tau_H $};
    		\node[draw,outer sep=0,inner sep=2,minimum size=16,fill=white] at (0.25,0.25) {$ \tau_G $};
    		\end{tikzpicture}
        \end{array}.
    \end{align*}

	\end{DEF}

The (Drinfeld) \emph{centre} of $ \C $, denoted $ Z(\C) $, is a category with objects of the form $ (X,\tau) $ where $ X $ is in $ \C $ and $ \tau $ is a half braiding on $ X $. $ \Hom_{Z(\C)}((X,\tau),(Y,\gamma)) $ is then given by the subspace of $ \Hom_{\C}(X,Y) $ defined by the condition that $ f \in \Hom_{\C}(X,Y) $ satisfies
	\begin{align*}
	(f \otimes \id_G) \circ \tau_G = \gamma_G \circ (\id_G \otimes f)
	\end{align*}
for all $ G $ in $ \C $. This category is monoidal~\cite[Section 7.13]{Etingof15} with tensor product $ (X,\tau) \otimes (Y,\gamma) = (X \otimes Y,\iota) $ where $ \iota_G = (\id_X \otimes \gamma_G) \circ (\tau_G \otimes \id_Y) $. In fact $ Z(\C) $ also admits a natural braiding~\cite[Theorem XIII.4.2.]{MR1643398} given by
	\begin{align} \label{eq:def_braiding_on_zc}
	\sigma_{(X,\tau),(Y,\gamma)} =  \gamma_X.
	\end{align}

Now, for the remainder of this section, let $ \C $ be a spherical fusion category. Then $ \C $ satisfies the conditions of Proposition~\ref{prop:yon_equivalence} and so the Yoneda embedding gives an equivalence between $ \C $ and $ \RC $. This induces an equivalence
       \begin{align*}
       Z(\C) &\to Z(\RC) \\
       (X,\tau) &\mapsto (X\Sharp,\tau\Sharp).
       \end{align*}
Here $ \tau\Sharp $ is a natural isomorphism from $ (\Blank \otimes X)\Sharp $ to $ (X \otimes \Blank)\Sharp $ such that
    \begin{align*}
	\tau_{GH}\Sharp = (\tau_G \otimes \id_H)\Sharp \circ (\id_G \otimes \tau_H)\Sharp.
    \end{align*}
By Lemma~\ref{lem:duals_and_yoneda} this is equivalent to an isomorphism
    \begin{align*}
    \kappa_{G^\vee,Z} \colon X\Sharp(G^\vee Z) \to X\Sharp(Z G^\vee)
    \end{align*}
that is natural in both $ Z $ and $ G^\vee $ and satisfies
    \begin{align*}
	\kappa_{(HG)^\vee,Z} = \kappa_{H^\vee,ZG^\vee} \circ \kappa_{G^\vee,H^\vee Z}
    \end{align*}
or, more simply,
    \begin{align*}
	\kappa_{GH,Z} = \kappa_{H,ZG} \circ \kappa_{G,HZ}.
    \end{align*}
Therefore an object in $ Z(\RC) $ may simply be thought of as an object in $ \RC $ together with isomorphisms $ \kappa_{G,Z} $ as above, i.e.\ an object in $ \RTC $. Furthermore, we recall that a morphism in $ Z(\C) $ between $ (X,\tau) $ and $ (Y,\tau') $ is a map $ f \in \Hom_{\C}(X,Y) $ such that
    \begin{align*}
	(f \otimes \id_G) \circ \tau_G = \tau'_G \circ (\id_G \otimes f).
    \end{align*}
Applying the Yoneda embedding we get
    \begin{align*}
        \begin{array}{ccc}
		\Big( (f \otimes \id_G)\Sharp \circ \tau_G\Sharp \Big)_Z = \Big((\tau')_G\Sharp \circ (\id_G \otimes f)\Sharp \Big)_Z
        \end{array}
    \end{align*}
which is equivalent to
    \begin{align*}
    f\Sharp_{Z G^\vee} \circ \kappa_{G^\vee,Z} =  \kappa'_{G^\vee,Z} \circ f_{G^\vee Z}\Sharp
    \end{align*}
which is precisely the condition identified by Proposition~\ref{prop:morphisms_in_rtc} as characterising a morphism from $ (X\Sharp, \kappa_{G^\vee,Z}) $ to $ (Y\Sharp,\kappa'_{G^\vee,Z}) $ in $ \RTC $. Therefore $ \RTC $ and $ Z(\C) $ are equivalent. To see this equivalence more explicitly let $ (X,\tau) $ be in $ Z(\C) $. The corresponding object $ F $ in $ \RTC $ is given on objects by
	\begin{align*}
	F(Y) =  \Hom_{\C}(Y,X)
	\end{align*}
and on morphisms by
	\begin{align} \begin{split} \label{eq:functor_f}
    F(\alpha_G) \colon \Hom_{\C}(Y,X) &\to \Hom_{\C}(Z,X) \\
	g &\mapsto
		\begin{array}{c}
			\begin{tikzpicture}[scale=0.15,every node/.style={inner sep=0,outer sep=-1}]
			\node (v19) at (0.75,-3) {};
			\node (v20) at (0.75,-12.25) {};
			\node (v16) at (-2,-3) {};
			\node (v18) at (0.75,-3) {};
			\node (v17) at (0.75,3.25) {};
			\node (v15) at (-2,1.25) {};
			\draw [thick] (v15) edge (v16);
			\draw [thick] (v17) edge (v18);
			\draw [thick] (v19) edge (v20);
			\node (v5) at (0.75,3.25) {};
			\node (v50) at (-2,1.75) {};
			\node (v7) at (-5.25,-13.5) {};
			\node (v70) at (-2,-6) {};
			\node (v8) at (-5.25,1.75) {};
			\node at (0.75,5.25) {$Z$};
			\node at (-7.5,-2.5) {$G^\vee$};
			\node at (-5.25,-15.5) {$X$};
			\node at (2.5,-2.5) {$G$};
			\node [draw,diamond,outer sep=0,inner sep=2,minimum size=20,fill=white] (v6) at (-0.5,0) {$ \alpha $};
			\node [draw,outer sep=0,inner sep=1,minimum size=15,fill=white] (v700) at (-2,-5.75) {$ g $};
			\draw[thick]  (v6)+(1.25,3.25) node (v11) {} edge (v5);
			\draw[thick]  (v6)+(-1.5,1.25) node (v3) {} edge (v50);
			\draw[thick]  (v6)+(-1.5,-3) node (v9) {} edge (v700);
			\draw[thick]  (v6)+(-4.75,-11) node (v14) {} edge (v7);
			\draw[thick]  (v50) to[out=90,in=90] (v8);
			\node (v1) at (-5.25,-9.25) {};
			\draw [thick] (v8) edge (v1);
			\node (v4) at (-2,-10.25) {};
			\node (v10) at (-2,-12.25) {};
			\draw [thick] (v4) edge (v10);
			\draw [thick] (v10) to[out=-90,in=-90] (v20);
			\draw [thick] (v700) edge (v4);
			\node[draw,outer sep=0,inner sep=2.3,minimum width=27,,fill=white] (v2) at (-3.75,-10.25) {$ \tau_{G^\vee} $};
			\end{tikzpicture}
		\end{array}.
    \end{split}
    \end{align}
where $ \alpha_G \in \Hom_{\TC}(Z,Y) $. This functor may also be understood in terms of an idempotent in $ \TC $. Indeed, for $ (X,\tau) $ in $ Z(\C) $ we may consider the following endomorphism in $ \TC $,
	\begin{align*}
	\e_{\tau} = \frac{1}{d(\C)}
        \bigoplus\limits_S d(S)
		\begin{array}{c}
			\begin{tikzpicture}[scale=0.25,every node/.style={inner sep=0,outer sep=-1}]
			\node (v1) at (0,4) {};
			\node (v4) at (0,-4) {};
			\node (v2) at (4,0) {};
			\node (v3) at (-4,0) {};
			\node (v5) at (-2,2) {};
			\node (v6) at (2,-2) {};
			\node (v7) at (2,2) {};
			\node (v11) at (-2,-2) {};
			\node [draw,diamond,outer sep=0,inner sep=.5,minimum size=28,fill=white] (v9) at (0,0) {$ \tau_{S} $};
			\node at (3,3) {$ X $};
			\node at (-3,-3) {$ X $};
			\node at (-3,3) {$ S $};
			\node at (3,-3) {$ S $};
			\draw [thick] (v9) edge (v6);
			\draw [thick] (v9) edge (v7);
			\draw [thick] (v9) edge (v11);
			\draw [thick] (v5) edge (v9);
			\draw[very thick, red]  (v1) edge (v3);
			\draw[very thick, red]  (v2) edge (v4);
			\draw[very thick]  (v1) edge (v2);
			\draw[very thick]  (v3) edge (v4);
			\end{tikzpicture}
		\end{array}.
	\end{align*}

    \begin{PROP} \label{prop:handle_slide}

    For $ Y $ in $ \C $ and $ \alpha_G \in \Hom_{\TC}(Y,X) $ we have
        \begin{align*}
        \e_\tau \circ \alpha_G = \frac{1}{d(\C)} \bigoplus\limits_R d(R)
            \begin{array}{c}
                \begin{tikzpicture}[scale=0.45,every node/.style={inner sep=0,outer sep=-1}]
                \node (v1) at (-1.125,2.875) {};
                \node (v4) at (-1.25,-5.25) {};
                \node (v2) at (2.875,-1.125) {};
                \node (v3) at (-5.25,-1.25) {};
                \node (v9) at (0.875,0.875) {};
                \node (v6) at (-3.25,-3.25) {};
                \node (v7) at (-4.375,-0.375) {};
                \node at (1.375,-1.25) {$ G $};
                \node at (-1.25,1.375) {$ G $};
                \node at (1.5,1.5) {$ Y $};
                \node at (-3.875,-3.875) {$ X $};
                \node at (-1.5,-0.625) {$ X $};
                \node at (0.25,-5) {$ R $};
                \node at (-5,0.25) {$ R $};
                \node (v11) at (-2.75,-1.25) {};
                \node (v12) at (-1.25,-2.75) {};
                \draw [thick] (v11) edge (v12);
                \node (v15) at (-0.75,0.5) {};
                \draw [thick] (v15) to[out=135,in=135] (v11);
                \node (v18) at (0.5,-0.75) {};
                \draw [thick] (v12) to[out=-45,in=-45] (v18);
                \node (v8) at (-0.375,-4.375) {};
                \draw[thick]  (v7) edge (v8);
                \draw[very thick, red]  (v1) edge (v3);
                \draw[very thick, red]  (v2) edge (v4);
                \draw[very thick]  (v1) edge (v2);
                \draw[very thick]  (v3) edge (v4);
                \draw [thick] (v15) edge (v18);
                \node [diamond,draw,outer sep=0,inner sep=0.3,minimum size=28,fill=white] (v50) at (-2.25,-2.25) {\mbox{ $ \tau_{{RG^\vee}} $}};
                \node [diamond,draw,outer sep=0,inner sep=-0.2,minimum size=23,fill=white] (v5) at (-0.125,-0.125) {\mbox{$ \alpha $}};
                \draw[thick]  (v50) edge (v5);
                \draw [thick] (v50) edge (v6);
                \draw[thick]  (v9) edge (v5);
                \end{tikzpicture}
            \end{array}
        \end{align*}
    and, for $ \beta_G \in \Hom_{\TC}(X,Y) $, we have
        \begin{align*}
        \beta_G \circ \e_\tau  = \frac{1}{d(\C)} \bigoplus\limits_R d(R)
            \begin{array}{c}
				\begin{tikzpicture}[scale=0.45,every node/.style={inner sep=0,outer sep=-1}]
				\node (v1) at (-1.125,2.875) {};
				\node (v4) at (-1.25,-5.25) {};
				\node (v2) at (2.875,-1.125) {};
				\node (v3) at (-5.25,-1.25) {};
				\node (v9) at (0.875,0.875) {};
				\node (v6) at (-3.25,-3.25) {};
				\node (v7) at (-1.75,2.25) {};
				\node at (-0.875,-3.75) {$ G $};
				\node at (-3.75,-1.125) {$ G $};
				\node at (1.5,1.5) {$ X $};
				\node at (-3.875,-3.875) {$ Y $};
				\node at (-2,-0.875) {$ X $};
				\node at (2.875,-2.375) {$ R $};
				\node at (-2.375,2.875) {$ R $};
				\node (v11) at (-3.125,-1.625) {};
				\node (v12) at (-1.625,-3.125) {};
				\draw [thick] (v11) edge (v12);
				\node (v15) at (-1,0.25) {};
				\draw [thick] (v15) to[out=135,in=135] (v11);
				\node (v18) at (0.25,-1) {};
				\draw [thick] (v12) to[out=-45,in=-45] (v18);
				\node (v8) at (2.25,-1.75) {};
				\draw[thick]  (v7) edge (v8);
				\draw[very thick, red]  (v1) edge (v3);
				\draw[very thick, red]  (v2) edge (v4);
				\draw[very thick]  (v1) edge (v2);
				\draw[very thick]  (v3) edge (v4);
				\draw [thick] (v15) edge (v18);
				\draw [thick] (v9) edge (v6);
				\node [diamond,draw,outer sep=0,inner sep=0.3,minimum size=28,fill=white] (v50) at (-0.125,-0.125) {\mbox{ $ \tau_{G^\vee R} $}};
				\node [diamond,draw,outer sep=0,inner sep=-0.2,minimum size=23,fill=white] (v5) at (-2.375,-2.375) {\mbox{$ \beta $}};
				\end{tikzpicture}
            \end{array}.
        \end{align*}
    \proof

    By the definition of composition in $ \TC $, we have
        \begin{align*}
        \e_\tau \circ \alpha_G &= \frac{1}{d(\C)} \bigoplus\limits_{R}\sum\limits_{S,b} d(S)
            \begin{array}{c}
				\begin{tikzpicture}[scale=0.5,every node/.style={inner sep=0,outer sep=-1}]
				\node (v1) at (-1.25,2.75) {};
				\node (v4) at (-1.25,-5.25) {};
				\node (v2) at (2.75,-1.25) {};
				\node (v3) at (-5.25,-1.25) {};
				\node (v9) at (0.75,0.75) {};
				\node (v6) at (-3.25,-3.25) {};
				\node (v70) at (0.75,-3.25) {};
				\node (v7) at (-3.25,0.75) {};
				\node at (0.75,-1.25) {$ G $};
				\node at (-3.25,-1) {$ S $};
				\node at (1.25,1.25) {$ Y $};
				\node at (-3.75,-3.75) {$ X $};
				\node at (-1.625,-0.875) {$ X $};
				\node at (1.2,-3.7) {$ R $};
				\node at (-3.75,1.25) {$ R $};
				\node (v11) at (-2.625,-1.875) {};
				\node (v12) at (-1.875,-2.625) {};
				\draw [thick] (v11) edge (v12);
				\node (v14) at (-0.625,0.125) {};
				\node (v13) at (0.125,-0.625) {};
				\draw [thick] (v13) edge (v14);
				\node (v15) at (-2.45,-0.425) {};
				\node (v16) at (-2.075,-0.05) {};
				\draw [thick] (v15) to[out=-45,in=135] (v11);
				\draw [thick] (v16) to[out=-45,in=135] (v14);
				\node (v17) at (-0.05,-2.075) {};
				\node (v18) at (-0.425,-2.45) {};
				\draw [thick] (v13) to[out=-45,in=135] (v17);
				\draw [thick] (v12) to[out=-45,in=135] (v18);
				\node [diamond,draw,outer sep=0,inner sep=0.3,minimum size=25,fill=white] (v50) at (-2.25,-2.25) {\mbox{$ \tau_S $}};
				\node [diamond,draw,outer sep=0,inner sep=-0.2,minimum size=25,fill=white] (v5) at (-0.25,-0.25) {\mbox{$ \alpha $}};
				\node [draw,rotate=45,outer sep=0,inner sep=2,minimum height=10,minimum width=13,fill=white] (v8) at (-2.5,0) {$ b $};
				\node [draw,rotate=45,outer sep=0,inner sep=2,minimum height=10,minimum width=13,fill=white] (v10) at (0,-2.5) {$ b^* $};
				\draw[thick]  (v50) edge (v5);
				\draw [thick] (v50) edge (v6);
				\draw[thick]  (v7) edge (v8);
				\draw[thick]  (v10) edge (v70);
				\draw[thick]  (v9) edge (v5);
				\draw[very thick, red]  (v1) edge (v3);
				\draw[very thick, red]  (v2) edge (v4);
				\draw[very thick]  (v1) edge (v2);
				\draw[very thick]  (v3) edge (v4);
				\end{tikzpicture}
            \end{array}\\
        &= \frac{1}{d(\C)} \bigoplus\limits_{R}\sum\limits_{S,b} d(S)
            \begin{array}{c}
				\begin{tikzpicture}[scale=0.5,every node/.style={inner sep=0,outer sep=-1}]
				\node (v1) at (-1.25,2.75) {};
				\node (v4) at (-1.25,-5.25) {};
				\node (v2) at (2.75,-1.25) {};
				\node (v3) at (-5.25,-1.25) {};
				\node (v9) at (0.75,0.75) {};
				\node (v6) at (-3.25,-3.25) {};
				\node (v70) at (-0.375,-4.375) {};
				\node (v7) at (-4.375,-0.375) {};
				\node (v80) at (-3.75,-1) {};
				\node at (0.875,-1.375) {$ G $};
				\node at (-1.5,1) {$ G $};
				\node at (1.25,1.25) {$ Y $};
				\node at (-3.75,-3.75) {$ X $};
				\node at (-1,-1.375) {$ X $};
				\node at (0.075,-4.825) {$ R $};
				\node at (-4.875,0.125) {$ R $};
				\node (v11) at (-0.375,-3) {};
				\node (v12) at (-2.25,-2.875) {};
				\node (v14) at (-0.625,0.125) {};
				\node (v13) at (-0.75,-2.625) {};
				\node (v15) at (0.125,-0.625) {};
				\node (v16) at (-3.5,-0.625) {};
				\draw [thick] (v15) to[out=-45,in=-45] (v11);
				\draw [thick] (v16) to[out=135,in=135] (v14);
				\node (v17) at (-1.25,-3.125) {};
				\node (v18) at (-1.55,-3.575) {};
				\draw [thick] (v13) to[out=135,in=135] (v17);
				\draw [thick] (v12) to[out=-45,in=135] (v18);
				\draw[thick]  (v7) edge (v80);
				\draw[very thick, red]  (v1) edge (v3);
				\draw[very thick, red]  (v2) edge (v4);
				\draw[very thick]  (v1) edge (v2);
				\draw[very thick]  (v3) edge (v4);
				\node (v19) at (-3,-1.75) {};
				\node (v22) at (-2.625,-1.5) {};
				\node (v20) at (-2.375,-2.375) {};
				\draw [thick] (v19) to[out=-45,in=135] (v20);
				\node (v21) at (-1.75,-2.625) {};
				\node (v23) at (-1.25,-2.125) {};
				\draw [thick] (v21) to[out=-45,in=-45] (v23);
				\draw [thick] (v22) to[out=-45,in=135] (v23);
				\node (v24) at (-1.875,-2.5) {};
				\draw [thick] (v21) edge (v24);
				\node (v25) at (-1.125,-3.25) {};
				\draw [thick] (v17) edge (v25);
				\node [draw,rotate=45,outer sep=0,inner sep=2,minimum height=10,minimum width=13,fill=white] (v8) at (-2.25,-2.5) {$ b $};
				\node [diamond,draw,outer sep=0,inner sep=0.3,minimum size=25,fill=white] (v50) at (-3.25,-1.25) {\mbox{$ \tau_{RG^\vee} $}};
				\node [diamond,draw,outer sep=0,inner sep=-0.2,minimum size=25,fill=white] (v5) at (-0.25,-0.25) {\mbox{$ \alpha $}};
				\node [draw,rotate=45,outer sep=0,inner sep=2,minimum height=10,minimum width=13,fill=white] (v10) at (-1.125,-3.625) {$ b^* $};
				\draw[thick]  (v9) to[out=-135,in=45] (v5);
				\draw[thick]  (v50) to[out=45,in=-135] (v5);
				\draw [thick] (v50) to[out=-135,in=45] (v6);
				\draw[thick]  (v10) edge (v70);
				\draw [thick] (v13) edge (v11);
				\end{tikzpicture}
            \end{array}\\
        &= \frac{1}{d(\C)} \bigoplus\limits_R d(R)
            \begin{array}{c}
                \begin{tikzpicture}[scale=0.45,every node/.style={inner sep=0,outer sep=-1}]
                \node (v1) at (-1.125,2.875) {};
                \node (v4) at (-1.25,-5.25) {};
                \node (v2) at (2.875,-1.125) {};
                \node (v3) at (-5.25,-1.25) {};
                \node (v9) at (0.875,0.875) {};
                \node (v6) at (-3.25,-3.25) {};
                \node (v7) at (-4.375,-0.375) {};
                \node at (1.375,-1.25) {$ G $};
                \node at (-1.25,1.375) {$ G $};
                \node at (1.5,1.5) {$ Y $};
                \node at (-3.875,-3.875) {$ X $};
                \node at (-1.5,-0.625) {$ X $};
                \node at (0.25,-5) {$ R $};
                \node at (-5,0.25) {$ R $};
                \node (v11) at (-2.75,-1.25) {};
                \node (v12) at (-1.25,-2.75) {};
                \draw [thick] (v11) edge (v12);
                \node (v15) at (-0.75,0.5) {};
                \draw [thick] (v15) to[out=135,in=135] (v11);
                \node (v18) at (0.5,-0.75) {};
                \draw [thick] (v12) to[out=-45,in=-45] (v18);
                \node (v8) at (-0.375,-4.375) {};
                \draw[thick]  (v7) edge (v8);
                \draw[very thick, red]  (v1) edge (v3);
                \draw[very thick, red]  (v2) edge (v4);
                \draw[very thick]  (v1) edge (v2);
                \draw[very thick]  (v3) edge (v4);
                \draw [thick] (v15) edge (v18);
                \node [diamond,draw,outer sep=0,inner sep=0.3,minimum size=28,fill=white] (v50) at (-2.25,-2.25) {\mbox{ $ \tau_{{RG^\vee}} $}};
                \node [diamond,draw,outer sep=0,inner sep=-0.2,minimum size=23,fill=white] (v5) at (-0.125,-0.125) {\mbox{$ \alpha $}};
                \draw[thick]  (v50) edge (v5);
                \draw [thick] (v50) edge (v6);
                \draw[thick]  (v9) edge (v5);
                \end{tikzpicture}
            \end{array}
        \end{align*}
    where the first equality is achieved by pushing $ b $ though $ \tau $ and the second by applying Proposition~\ref{lem:dual_decompose}. This proves the first equality, the second is proved analogously. \endproof

    \end{PROP}

    \begin{COR}

    Let $ (X,\tau) $ be in $ Z(\C) $. Then $ \e_\tau $ is an idempotent.

    \proof By Proposition~\ref{prop:handle_slide}, we have
        \begin{align*}
        \e_\tau \circ \e_\tau &= \frac{1}{d(\C)^2} \bigoplus\limits_{R} \sum\limits_S d(R) d(S)
            \begin{array}{c}
                \begin{tikzpicture}[scale=0.45,every node/.style={inner sep=0,outer sep=-1}]
                \node (v1) at (-1.125,2.875) {};
                \node (v4) at (-1.25,-5.25) {};
                \node (v2) at (2.875,-1.125) {};
                \node (v3) at (-5.25,-1.25) {};
                \node (v9) at (0.875,0.875) {};
                \node (v6) at (-3.25,-3.25) {};
                \node (v7) at (-4.375,-0.375) {};
                \node at (1.375,-1.25) {$ S $};
                \node at (-1.25,1.375) {$ S $};
                \node at (1.5,1.5) {$ X $};
                \node at (-3.875,-3.875) {$ X $};
                \node at (-1.5,-0.625) {$ X $};
                \node at (0.25,-5) {$ R $};
                \node at (-5,0.25) {$ R $};
                \node (v11) at (-2.75,-1.25) {};
                \node (v12) at (-1.25,-2.75) {};
                \draw [thick] (v11) edge (v12);
                \node (v15) at (-0.75,0.5) {};
                \draw [thick] (v15) to[out=135,in=135] (v11);
                \node (v18) at (0.5,-0.75) {};
                \draw [thick] (v12) to[out=-45,in=-45] (v18);
                \node (v8) at (-0.375,-4.375) {};
                \draw[thick]  (v7) edge (v8);
                \draw[very thick, red]  (v1) edge (v3);
                \draw[very thick, red]  (v2) edge (v4);
                \draw[very thick]  (v1) edge (v2);
                \draw[very thick]  (v3) edge (v4);
                \draw [thick] (v15) edge (v18);
                \node [diamond,draw,outer sep=0,inner sep=0.3,minimum size=28,fill=white] (v50) at (-2.25,-2.25) {\mbox{ $ \tau_{{RS^\vee}} $}};
                \node [diamond,draw,outer sep=0,inner sep=-0.2,minimum size=23,fill=white] (v5) at (-0.125,-0.125) {\mbox{$ \tau_S $}};
                \draw[thick]  (v50) edge (v5);
                \draw [thick] (v50) edge (v6);
                \draw[thick]  (v9) edge (v5);
                \end{tikzpicture}
            \end{array}\\
        &= \frac{1}{d(\C)^2} \bigoplus\limits_{R} \sum\limits_S d(R) d(S)
            \begin{array}{c}
                \begin{tikzpicture}[scale=0.45,every node/.style={inner sep=0,outer sep=-1}]
                \node (v1) at (-2.75,1.25) {};
                \node (v4) at (-1.25,-5.25) {};
                \node (v2) at (1.25,-2.75) {};
                \node (v3) at (-5.25,-1.25) {};
                \node (v9) at (-0.75,-0.75) {};
                \node (v6) at (-3.25,-3.25) {};
                \node (v7) at (-4.375,-0.375) {};
                \node at (0,-2.75) {$ S $};
                \node at (-2.75,0) {$ S $};
                \node at (-0.125,-0.125) {$ X $};
                \node at (-3.875,-3.875) {$ X $};
                \node at (0.25,-4.875) {$ R $};
                \node at (-4.875,0.25) {$ R $};
                \node (v11) at (-3,-1) {};
                \node (v12) at (-1,-3) {};
                \draw [thick] (v11) edge (v12);
                \node (v15) at (-2.625,-0.625) {};
                \draw [thick] (v15) to[out=135,in=135] (v11);
                \node (v18) at (-0.625,-2.625) {};
                \draw [thick] (v12) to[out=-45,in=-45] (v18);
                \node (v8) at (-0.375,-4.375) {};
                \draw[thick]  (v7) edge (v8);
                \draw[very thick, red]  (v1) edge (v3);
                \draw[very thick, red]  (v2) edge (v4);
                \draw[very thick]  (v1) edge (v2);
                \draw[very thick]  (v3) edge (v4);
                \node (v5) at (-2.125,-1.125) {};
                \node (v10) at (-1.125,-2.125) {};
                \draw [thick] (v15) edge (v5);
                \draw [thick] (v18) edge (v10);
                \draw [thick] (v9) edge (v6);
                \node [diamond,draw,outer sep=0,inner sep=0.3,minimum size=32,fill=white] (v50) at (-2,-2) {\mbox{ $ \tau_{{RS^\vee S}} $}};
                \end{tikzpicture}
            \end{array} = \e_\tau.
        \end{align*}  \endproof
    \end{COR}

    \begin{PROP} \label{prop:equivalence_with_f}

    Let $ (X,\tau) $ be in $ Z(\C) $, let $ F $ in $ \RTC $ be given by~\eqref{eq:functor_f} and let $ (X,\e_\tau)\Sharp $ be as defined in Section~\ref{sec:prelim_lin_cat}. Then $ F \cong (X,\e_\tau)\Sharp $.

    \proof

    We consider the following two linear maps,
        \begin{align*}
        \Xi_Y \colon \Hom_{\C}(Y,X) &\to \Hom_{\TC}(Y,\e_\tau) \\
        \alpha &\mapsto \e_\tau \circ \alpha
        \end{align*}
    and
        \begin{align*}
        \Psi_Y \colon \Hom_{\TC}(Y,\e_\tau) &\to \Hom_{\C}(Y,X)  \\
        \beta_G &\mapsto
            \begin{array}{c}
    			\begin{tikzpicture}[scale=0.15,every node/.style={inner sep=0,outer sep=-1}]
    			\node (v19) at (1.25,-3) {};
    			\node (v20) at (1.25,-6.5) {};
    			\node (v16) at (-1.5,-4) {};
    			\node (v18) at (1.25,-3) {};
    			\node (v17) at (1.25,4.25) {};
    			\node (v15) at (-1.5,2.25) {};
    			\draw [thick] (v15) edge (v16);
    			\draw [thick] (v17) edge (v18);
    			\draw [thick] (v19) edge (v20);
    			\node (v5) at (1.25,4.25) {};
    			\node (v50) at (-1.5,2.75) {};
    			\node (v7) at (-4.75,-7.75) {};
    			\node (v8) at (-4.75,2.75) {};
    			\node at (1.25,6.25) {$Y$};
    			\node at (-7,-1.5) {$G^\vee$};
    			\node at (-4.75,-9.75) {$X$};
    			\node at (3,-1.5) {$G$};
    			\draw[thick]  (1.25,4.25) node (v11) {} edge (v5);
    			\draw[thick]  (-1.5,2.25) node (v3) {} edge (v50);
    			\draw[thick]  (-4.75,-5.25) node (v14) {} edge (v7);
    			\draw[thick]  (v50) to[out=90,in=90] (v8);
    			\node (v1) at (-4.75,-3.5) {};
    			\draw [thick] (v8) edge (v1);
    			\node (v4) at (-1.5,-4.5) {};
    			\node (v10) at (-1.5,-6.5) {};
    			\draw [thick] (v4) edge (v10);
    			\draw [thick] (v10) to[out=-90,in=-90] (v20);
    			\node[draw,outer sep=0,inner sep=2.3,minimum width=27,,fill=white] (v2) at (-3.25,-4.5) {$ \tau_{G^\vee} $};
    			\node [draw,diamond,outer sep=0,inner sep=2,minimum size=20,fill=white] (v6) at (0,1) {$ \beta $};
    			\end{tikzpicture}
            \end{array}.
        \end{align*}
    For $ \alpha \in \Hom_{\C}(Y,X) $, we have
        \begin{align*}
        \Psi \circ \Xi (\alpha) = \frac{1}{d(\C)} \sum\limits_S d(S)
            \begin{array}{c}
    			\begin{tikzpicture}[scale=0.15,every node/.style={inner sep=0,outer sep=-1}]
    			\node (v19) at (1.5,-2.75) {};
    			\node (v20) at (1.5,-6.25) {};
    			\node (v16) at (-1.25,-3.75) {};
    			\node (v18) at (1.5,-2.75) {};
    			\node (v17) at (1.5,5.5) {};
    			\node (v15) at (-1.25,2.5) {};
    			\draw [thick] (v15) edge (v16);
    			\draw [thick] (v17) edge (v18);
    			\draw [thick] (v19) edge (v20);
    			\node (v5) at (1.5,5.5) {};
    			\node (v50) at (-1.25,3) {};
    			\node (v7) at (-4.5,-7.5) {};
    			\node (v8) at (-4.5,3) {};
    			\node at (1.5,12) {$Y$};
    			\node at (-6.75,-1.25) {$S^\vee$};
    			\node at (-4.5,-9.5) {$X$};
    			\node at (3.25,-1.25) {$S$};
    			\draw[thick]  (1.5,10) node (v11) {} edge (v5);
    			\draw[thick]  (-1.25,2.5) node (v3) {} edge (v50);
    			\draw[thick]  (-4.5,-5) node (v14) {} edge (v7);
    			\draw[thick]  (v50) to[out=90,in=90] (v8);
    			\node (v1) at (-4.5,-3.25) {};
    			\draw [thick] (v8) edge (v1);
    			\node (v4) at (-1.25,-4.25) {};
    			\node (v10) at (-1.25,-6.25) {};
    			\draw [thick] (v4) edge (v10);
    			\draw [thick] (v10) to[out=-90,in=-90] (v20);
    			\node[draw,outer sep=0,inner sep=2.3,minimum width=27,,fill=white] (v2) at (-3,-4.25) {$ \tau_{S^\vee} $};
    			\node [draw,diamond,outer sep=0,inner sep=2,minimum size=20,fill=white] (v6) at (0.25,1.25) {$ \tau_S $};
    			\node [draw,outer sep=0,inner sep=2,minimum size=16,fill=white] (v6) at (1.5,7) {$ \alpha $};
    			\end{tikzpicture}
            \end{array}
        = \frac{1}{d(\C)} \sum\limits_S d(S)
            \begin{array}{c}
				\begin{tikzpicture}[scale=0.15,every node/.style={inner sep=0,outer sep=-1}]
				\node (v17) at (1.5,4.75) {};
				\node (v5) at (1.5,10) {};
				\node at (1.5,12) {$Y$};
				\node at (5,-3.25) {$S$};
				\node at (-2.5,-9.5) {$X$};
				\node at (3.5,2.25) {$X$};
				\node (v1) at (1.5,-0.25) {};
				\draw [thick] (v5) edge (v1);
				\node (v4) at (-1.5,1.75) {};
				\node (v7) at (-4,-0.75) {};
				\node (v12) at (-0.5,-4.25) {};
				\node (v10) at (2,-1.75) {};
				\node (v2) at (-2.5,-2.75) {};
				\node (v3) at (-2.5,-7.5) {};
				\draw [thick] (v2) edge (v3);
				\draw [thick] (v4) to[out=135,in=135] (v7);
				\node (v9) at (-0.5,0.75) {};
				\node (v8) at (-3,-1.75) {};
				\node (v11) at (3,-2.75) {};
				\node (v13) at (0.5,-5.25) {};
				\draw [thick] (v7) edge (v8);
				\draw [thick] (v9) edge (v4);
				\draw [thick] (v10) edge (v11);
				\draw [thick] (v12) edge (v13);
				\draw [thick] (v11) to[out=-45,in=-45] (v13);
				\node [draw,diamond,outer sep=0,inner sep=2,minimum size=20,fill=white] (v6) at (-0.5,-1.75) {$ \tau_{S^\vee S} $};
				\node [draw,outer sep=0,inner sep=2,minimum size=16,fill=white] (v6) at (1.5,6.25) {$ \alpha $};
				\end{tikzpicture}
            \end{array}
        = \alpha
        \end{align*}
    and, for $ \beta_G \in \Hom_{\TC}(Y,\e_\tau) $, we have
        \begin{align*}
        \beta_G = \e_\tau \circ \beta_G &= \frac{1}{d(\C)} \bigoplus\limits_S d(S)
            \begin{array}{c}
				\begin{tikzpicture}[scale=0.5,every node/.style={inner sep=0,outer sep=-1}]
				\node (v1) at (-1.125,2.875) {};
				\node (v4) at (-1.25,-5.25) {};
				\node (v2) at (2.875,-1.125) {};
				\node (v3) at (-5.25,-1.25) {};
				\node (v9) at (0.875,0.875) {};
				\node (v6) at (-3.25,-3.25) {};
				\node (v7) at (-4.375,-0.375) {};
				\node at (0.875,-1.125) {$ G $};
				\node at (-1.5,-0.75) {$ X $};
				\node at (-1.125,0.875) {$ G $};
				\node at (1.375,1.375) {$ Y $};
				\node at (-3.75,-3.75) {$ X $};
				\node at (0.125,-4.75) {$ S $};
				\node at (-4.75,0.125) {$ S $};
				\node (v11) at (-2.5,-1.75) {};
				\node (v12) at (-1.75,-2.375) {};
				\draw [thick] (v11) edge (v12);
				\node (v15) at (-0.5,0.25) {};
				\draw [thick] (v15) to[out=135,in=135] (v11);
				\node (v18) at (0.25,-0.5) {};
				\draw [thick] (v12) to[out=-45,in=-45] (v18);
				\node (v8) at (-0.375,-4.375) {};
				\draw[thick]  (v7) edge (v8);
				\draw[very thick, red]  (v1) edge (v3);
				\draw[very thick, red]  (v2) edge (v4);
				\draw[very thick]  (v1) edge (v2);
				\draw[very thick]  (v3) edge (v4);
				\node [diamond,draw,outer sep=0,inner sep=0.3,minimum size=28,fill=white] (v50) at (-2.25,-2.25) {\mbox{ $ \tau_{{SG^\vee}} $}};
				\node [diamond,draw,outer sep=0,inner sep=-0.2,minimum size=25,fill=white] (v5) at (-0.125,-0.125) {\mbox{$ \beta $}};
				\draw[thick]  (v50) edge (v5);
				\draw [thick] (v50) edge (v6);
				\draw[thick]  (v9) edge (v5);
				\end{tikzpicture}
            \end{array}
        = \Xi \circ \Psi(\beta_G).
        \end{align*}
    As $ \Xi_Y $ and $ \Psi_Y $ are inverse we only have to check naturality for one of them. For $ \alpha \in \Hom_{\C}(Z,X) $ and $ \beta_G \in \Hom_{\TC}(Y,Z) $, we have
        \begin{align*}
        (X,\e_\tau)\Sharp (\beta_G) \circ \Xi_Z (\alpha) = \frac{1}{d(\C)} \bigoplus\limits_S d(S) 
            \begin{array}{c}
				\begin{tikzpicture}[scale=0.5,every node/.style={inner sep=0,outer sep=-1}]
				\node (v1) at (-1.125,2.875) {};
				\node (v4) at (-1.375,-5.375) {};
				\node (v2) at (2.875,-1.125) {};
				\node (v3) at (-5.375,-1.375) {};
				\node (v9) at (0.875,0.875) {};
				\node (v6) at (-3.375,-3.375) {};
				\node (v7) at (-4.5,-0.5) {};
				\node at (0.875,-1.125) {$ G $};
				\node at (-1.125,0.875) {$ G $};
				\node at (1.375,1.375) {$ Y $};
				\node at (-3.875,-3.875) {$ X $};
				\node at (0,-4.875) {$ S $};
				\node at (-4.875,0) {$ S $};
				\node (v11) at (-2.625,-1.875) {};
				\node (v12) at (-1.875,-2.625) {};
				\draw [thick] (v11) edge (v12);
				\node (v15) at (-0.5,0.25) {};
				\draw [thick] (v15) to[out=135,in=135] (v11);
				\node (v18) at (0.25,-0.5) {};
				\draw [thick] (v12) to[out=-45,in=-45] (v18);
				\node (v8) at (-0.5,-4.5) {};
				\draw[thick]  (v7) edge (v8);
				\draw[very thick, red]  (v1) edge (v3);
				\draw[very thick, red]  (v2) edge (v4);
				\draw[very thick]  (v1) edge (v2);
				\draw[very thick]  (v3) edge (v4);
				\node [diamond,draw,outer sep=0,inner sep=0.3,minimum size=28,fill=white] (v50) at (-2.375,-2.375) {\mbox{ $ \tau_{{SG^\vee}} $}};
				\node [diamond,draw,outer sep=0,inner sep=-0.2,minimum size=25,fill=white] (v5) at (-0.125,-0.125) {\mbox{$ \beta $}};
				\draw[thick]  (v50) edge (v5);
				\draw [thick] (v50) edge (v6);
				\draw[thick]  (v9) edge (v5);
				\node [rotate=-45,draw,outer sep=0,inner sep=-0.2,minimum size=16,fill=white] (v5) at (-1.125,-1.125) {\mbox{$ \alpha $}};
				\end{tikzpicture}
            \end{array}
        = \Xi_Y \circ F(\beta_G)(\alpha).
        \end{align*}
            \endproof
    \end{PROP}

This proposition has an interesting consequence. As $ Z(\C) $ and $ \RTC $ are equivalent it also proves that every functor $ F $ in $ \RTC $ is represented by an idempotent (namely $ \e_\tau $ where $ (X,\tau) $ is the corresponding object in $ Z(\C) $). In summary, we have the following.

    \begin{COR}

    The Yoneda embedding $ \yen \colon \TC \to \RTC $ is an idempotent completion.

    \end{COR}

Informally, we may interpret this result as allowing us to study $ \RTC $ (and therefore also the centre of $ \C $) by simply working with idempotents in $ \TC $. This idea is precisely what is meant by the terms `graphical approach' which appear in the title of this article. To illustrate this approach we shall now describe an alternative proof of the equivalence between $ Z(\C) $ and $ \CC $ when $ \C $ is modular.

\section{Equivalence with $ \CC $}

We start by equipping $ \C $ with a (balanced) braiding $ \sigma $; $ \C $ is now a pre-modular tensor category. As we have now chosen a braiding we get a (covariant) braided monoidal functor
 	\begin{align*}
		\begin{split}
	 	\Phi \colon \CC  &\to Z(\C)\\
		X \boxtimes Y &\mapsto (XY,  (\id_X \otimes \bar{\sigma_Y}) \circ (\sigma_X \otimes \id_Y) )
		\end{split}
	\end{align*}
where $ \boxtimes $ denotes the Deligne tensor product and $ \overline{\C} $ is obtained by equipping $ \C $ with the opposite braiding. It is also known that this functor is an equivalence if and only if $ \C $ is modular (see \cite{MR1990929} or \cite[Proposition 8.20.12]{Etingof15}). Our aim in this section is to provide an alternative proof of this result by exploiting graphical calculus in the tube category. First we note that, by the results of the previous section, we have $ \Phi(X\boxtimes Y) = (XY,\e_X^Y)\Sharp $ where
	\begin{align*}
		\e_X^Y = \frac{1}{d(\C)} \bigoplus\limits_S d(S)
			\begin{array}{c}
				\begin{tikzpicture}[scale=0.2,every node/.style={inner sep=0,outer sep=-1}]
				\node (v1) at (0,5) {};
				\node (v4) at (0,-5) {};
				\node (v2) at (5,0) {};
				\node (v3) at (-5,0) {};
				\node (v5) at (-2.5,2.5) {};
				\node (v6) at (2.5,-2.5) {};
				\node (v7) at (1.5,3.5) {};
				\node (v11) at (3.5,1.5) {};
				\node (v12) at (-1.5,-3.5) {};
				\node (v10) at (-3.5,-1.5) {};
				\node (v13) at (2.5,-2.5) {};
				\draw [thick] (v7) edge (v10);
				\draw [line width =0.5em,white] (v5) edge (v13);
				\draw [thick] (v5) edge (v13);
				\draw [line width =0.5em,white] (v11) edge (v12);
				\draw [thick] (v11) edge (v12);
				\node at (2.75,4.75) {$ X $};
				\node at (4.75,2.75) {$ Y $};
				\node at (-3.75,3.75) {$ S $};
				\node at (3.75,-3.5) {$ S $};
				\draw[very thick, red]  (v1) edge (v3);
				\draw[very thick, red]  (v2) edge (v4);
				\draw[very thick]  (v1) edge (v2);
				\draw[very thick]  (v3) edge (v4);
				\end{tikzpicture}
			\end{array}.
	\end{align*}
When we suppose $ \C $ is modular the killing ring lemma (Lemma~\ref{lemma:killing_ring}) allows us to compute the $ \Hom $-spaces between these idempotents.

	\begin{PROP} \label{prop:homs_between_idems}

	Let $ \C $ be a modular tensor category and let $ X,Y,A,B $ be in $ \C $. We have
		\begin{align*}
		\Hom_{\TC}(\e_X^Y,\e_A^B) = \Hom_{C}(X,A) \otimes \Hom_{C}(Y,B).
		\end{align*}
	\proof

	We consider the maps
		\begin{align*}
		\phi \colon \Hom_{C}(X,A) \otimes \Hom_{C}(Y,B) &\to \Hom_{\TC}(\e_X^Y,\e_A^B)\\
		f \otimes g &\mapsto  \e_A^B \circ (f \otimes g) \circ \e_X^Y
		\end{align*}
	and
		\begin{align*}
		\varphi \colon \Hom_{\TC}(\e_X^Y,\e_A^B) &\to \Hom_{C}(X,A) \otimes \Hom_{C}(Y,B) \\
		\alpha_G &\mapsto \sum\limits_{T,b,c} \frac{1}{d(T)}
			\begin{array}{c}
				\begin{tikzpicture}[scale=0.15,every node/.style={inner sep=0,outer sep=-1}]
				\node (v2) at (0,9) {$ X $};
				\node (v1) at (-4.5,2) {};
				\node (v3) at (-4.5,-3) {};
				\node (v7) at (-1,-9) {};
				\node at (-6,-12.5) {$A$};
				\node (v10) at (2,0) {};
				\node (v11) at (-1,-1) {};
				\node (v12) at (1,-1) {};
				\node (v13) at (1,-3) {};
				\node (v9) at (2,6.5) {};
				\draw [thick] (v9) edge (v10);
				\draw [thick] (v12) edge (v13);
				\node (v16) at (-2,1) {};
				\node (v17) at (-2,2) {};
				\node (v15) at (-6,-3.5) {};
				\node (v14) at (-2.5,0) {};
				\draw [thick] (v14) to[out=-135,in=90] (v15);
				\draw [thick] (v16) edge (v17);
				\draw [thick] (v17) to[out=90,in=90] (v1);
				\draw [line width=0.2cm,white] (v1) edge (v3);
				\draw [thick] (v1) edge (v3);
				\draw [line width=0.2cm,white] (v3) to[out=-90,in=-90] (v13);
				\draw [thick] (v3) to[out=-90,in=-90] (v13);
				\draw [line width=0.2cm,white] (v11) edge (v7);
				\draw [thick] (v11) edge (v7);
				\node at (-6.5,2) {$ G^\vee $};
				\node (v5) at (5,6.5) {};
				\draw [thick] (v9) to[out=90,in=90] (v5);
				\node (v8) at (5,-9) {};
				\draw [thick] (v5) edge (v8);
				\draw [thick] (v7) to[out=-90,in=-90] (v8);
				\node (v19) at (0,2) {};
				\node (v18) at (0,7.5) {};
				\draw [thick] (v18) edge (v19);
				\node [draw,outer sep=0,inner sep=1,minimum width=12, minimum height = 13,fill=white] (v9) at (-1,-7) {$ b^* $};
				\node [draw,outer sep=0,inner sep=1,minimum width=12, minimum height = 13,fill=white] (v9) at (2,4.5) {$ c $};
				\node [diamond,draw,outer sep=0,inner sep=3,minimum height=10,minimum width=24,fill=white] (v6) at (-0.5,0) {$ \alpha $};
				\node (v22) at (-6,-11) {};
				\draw [thick] (v15) edge (v22);
				\end{tikzpicture}
			\end{array}
		\otimes \! \!
			\begin{array}{c}
				\begin{tikzpicture}[scale=0.15,every node/.style={inner sep=0,outer sep=-1}]
				\node at (11.5,-1.5) {$ T $};
				\node at (9,9) {$ Y $};
				\node at (9,-12.5) {$ B $};
				\node (v20) at (9,7.5) {};
				\node (v21) at (9,-11) {};
				\draw [thick] (v20) edge (v21);
				\node [draw,outer sep=0,inner sep=1,minimum width=12, minimum height = 13,fill=white] (v9) at (9,4.5) {$ c^* $};
				\node [draw,outer sep=0,inner sep=1,minimum width=12, minimum height = 13,fill=white] (v9) at (9,-7) {$ b $};
				\end{tikzpicture}
			\end{array}.
		\end{align*}
	We have
		\begin{align*}
		\varphi \circ \phi (f \otimes g) &= \varphi \left( \frac{1}{d(\C)^2}\sum\limits_{R,T} d(R) d(T)
			\begin{array}{c}
				\begin{tikzpicture}[scale=0.5,every node/.style={inner sep=0,outer sep=-1}]
				\node (v1) at (0.5,3.5) {};
				\node (v4) at (-0.5,-3.5) {};
				\node (v2) at (3.5,0.5) {};
				\node (v3) at (-3.5,-0.5) {};
				\node (v8) at (-1.625,-0.375) {};
				\node (v80) at (0.375,1.625) {};
				\node (v11) at (-0.25,-1.75) {};
				\node (v110) at (1.75,0.25) {};
				\node [every node] (v10) at (-0.75,-1.25) {};
				\node [every node] (v100) at (1.25,0.75) {};
				\draw [thick] (v8) edge (v10);
				\node (v27) at (-2.5,0.5) {};
				\node (v270) at (-0.5,2.5) {};
				\draw [thick]  (v8) to[out=135,in=-45] (v27);
				\node (v34) at (0.5,-2.5) {};
				\node (v340) at (2.5,-0.5) {};
				\draw [thick] (v11) to[out=-45,in=135] (v34);
				\node (v5) at (1,0) {};
				\node (v50) at (-0.125,0.875) {};
				\node (v7) at (-1.5,-2.5) {};
				\node (v70) at (1.465,2.465) {};
				\draw [thick] (v5) edge (v7);
				\node (v12) at (-0.125,0.875) {};
				\node (v120) at (1,0) {};
				\node (v13) at (-1.25,-0.25) {};
				\node (v130) at (0.25,1.25) {};
				\draw [thick] (v12) edge (v13);
				\node (v14) at (-1.75,-0.75) {};
				\node (v140) at (0.75,1.75) {};
				\node (v15) at (-2.5,-1.5) {};
				\node (v150) at (2.5,1.5) {};
				\draw [thick] (v14) edge (v15);
				\draw [thick] (v120) edge (v150);
				\draw [thick] (v50) edge (v130);
				\draw [thick] (v140) edge (v70);
				\draw [thick] (v270) to[out=-45,in=135] (v80);
				\draw [thick] (v80) edge (v100);
				\draw [thick] (v110) to[out=-45,in=135] (v340);
				\node at (-0.875,2.875) {$ R $};
				\node at (-2.875,0.875) {$ T $};
				\node at (2,2.875) {$ X $};
				\node at (-3,-1.875) {$ A $};
				\node at (2.85,1.99) {$ Y $};
				\node at (-2,-2.875) {$ B $};
				\node [draw,rotate=-45,outer sep=0,inner sep=1,minimum size=13,fill=white] (v30) at (-0.5,0.5) {$ f $};
				\node [draw,rotate=-45,outer sep=0,inner sep=1,minimum size=13,fill=white] (v30) at (0.5,-0.5) {$ g $};
				\draw[very thick, red]  (v1) edge (v3);
				\draw[very thick, red]  (v2) edge (v4);
				\draw[very thick]  (v1) edge (v2);
				\draw[very thick]  (v3) edge (v4);
				\end{tikzpicture}
			\end{array} \right) \\
		&= \varphi \left( \frac{1}{d(\C)}\sum\limits_{R,T} d(R)
			\begin{array}{c}
				\begin{tikzpicture}[scale=0.5,every node/.style={inner sep=0,outer sep=-1}]
				\node (v1) at (0.5,3.5) {};
				\node (v4) at (0.5,-2.5) {};
				\node (v2) at (3.5,0.5) {};
				\node (v3) at (-2.5,0.5) {};
				\node (v80) at (0.375,1.625) {};
				\node (v110) at (1.75,0.25) {};
				\node [every node] (v100) at (1.25,0.75) {};
				\node (v270) at (-0.5,2.5) {};
				\node (v340) at (2.5,-0.5) {};
				\node (v5) at (1,0) {};
				\node (v50) at (-0.125,0.875) {};
				\node (v7) at (-0.5,-1.5) {};
				\node (v70) at (1.465,2.465) {};
				\draw [thick] (v5) edge (v7);
				\node (v12) at (-0.125,0.875) {};
				\node (v120) at (1,0) {};
				\node (v13) at (0,1) {};
				\node (v130) at (0.25,1.25) {};
				\draw [thick] (v12) edge (v13);
				\node (v14) at (-0.75,0.25) {};
				\node (v140) at (0.75,1.75) {};
				\node (v15) at (-1.5,-0.5) {};
				\node (v150) at (2.5,1.5) {};
				\draw [thick] (v14) edge (v15);
				\draw [thick] (v120) edge (v150);
				\draw [thick] (v50) edge (v130);
				\draw [thick] (v140) edge (v70);
				\draw [thick] (v270) to[out=-45,in=135] (v80);
				\draw [thick] (v80) edge (v100);
				\draw [thick] (v110) to[out=-45,in=135] (v340);
				\node at (-0.875,2.875) {$ R $};
				\node at (2,2.875) {$ X $};
				\node at (-2,-0.875) {$ A $};
				\node at (2.85,1.99) {$ Y $};
				\node at (-1,-1.875) {$ B $};
				\node [draw,rotate=-45,outer sep=0,inner sep=1,minimum size=13,fill=white] (v30) at (-0.5,0.5) {$ f $};
				\node [draw,rotate=-45,outer sep=0,inner sep=1,minimum size=13,fill=white] (v30) at (0.5,-0.5) {$ g $};
				\draw[very thick, red]  (v1) edge (v3);
				\draw[very thick, red]  (v2) edge (v4);
				\draw[very thick]  (v1) edge (v2);
				\draw[very thick]  (v3) edge (v4);
				\end{tikzpicture}
			\end{array} \right)\\
		&= \sum\limits_{T,b,c} b^*(g \circ c) \ f \otimes (b \circ c^*) = f \otimes g
		\end{align*}
	and, for $ \alpha_G \in \Hom_{\TC}(\e_X^Y,\e_A^B) $,
		\begin{align*}
		\alpha_G =  \e_A^B \circ \alpha_G \circ \e_X^Y = \frac{1}{d(\C)^2} \bigoplus_{T} \sum\limits_{R,b} d(R) d(T)
			\begin{array}{c}
				\begin{tikzpicture}[scale=0.4,every node/.style={inner sep=0,outer sep=-1}]
				\node (v1) at (1.475,5.5) {};
				\node (v4) at (0,-4.05) {};
				\node (v2) at (5.5,1.5) {};
				\node (v3) at (-4.05,0) {};
				\node (v8) at (0.75,2.5) {};
				\node (v10) at (2.5,0.75) {};
				\node (v11) at (-0.225,-1.8) {};
				\node (v19) at (0.25,1.75) {};
				\node (v26) at (1.75,0.25) {};
				\node (v50) at (0,0.5) {};
				\node (v16) at (2.5,4.5) {};
				\draw [thick] (v50) to[out=45,in=-135] (v16);
				\draw [thick] (v11) to[out=-45,in=-45] (v10);
				\node (v5) at (0,-0.5) {};
				\node (v7) at (-1,-3) {};
				\node (v12) at (-0.5,0) {};
				\node (v120) at (0.5,0) {};
				\node (v13) at (-2.925,-1.125) {};
				\draw [thick] (v12) to[in=45,out=-135] (v13);
				\node (v23) at (-0.75,0.75) {};
				\node (v34) at (0.675,-0.675) {};
				\node (v22) at (-1.75,-0.25) {};
				\draw [thick] (v8) to[out=135,in=135] (v22);
				\draw [thick] (v19) to[out=135,in=135] (v23);
				\node (v270) at (0.5,4.5) {};
				\node (v340) at (4.5,0.5) {};
				\draw [line width=0.2cm,white] (v270) edge (v340);
				\draw [thick] (v270) edge (v340);
				\draw [line width=0.2cm,white] (v22) edge (v11);
				\draw [thick] (v22) edge (v11);
				\node (v29) at (4.5,2.5) {};
				\draw [thick] (v26) to[out=-45,in=-45] (v34);
				\draw [line width=0.2cm,white] (v19) edge (v26);
				\draw [thick] (v19) edge (v26);
				\draw [line width=0.2cm,white] (v10) edge (v8);
				\draw [thick] (v10) edge (v8);
				\draw [line width=0.2cm,white] (v120) to[out=45,in=-135] (v29);
				\draw [thick] (v120) to[out=45,in=-135] (v29);
				\draw [line width=0.2cm,white] (v5) to[in=45,out=-135] (v7);
				\draw [thick] (v5) to[in=45,out=-135] (v7);
				\node at (-2.5,0.75) {$ R $};
				\node at (3.25,5) {$ X $};
				\node at (-3.5,-1.75) {$ A $};
				\node at (5.25,3) {$ Y $};
				\node at (-1.5,-3.75) {$ B $};
				\node at (-0.25,5) {$ T $};
				\node at (5.25,-0.25) {$ T $};
				\draw[very thick, red]  (v1) edge (v3);
				\draw[very thick, red]  (v2) edge (v4);
				\draw[very thick]  (v1) edge (v2);
				\draw[very thick]  (v3) edge (v4);
				\draw [thick] (v23) edge (v34);
				\node [draw,diamond,outer sep=0,inner sep=.5,minimum size=25,fill=white] (v9) at (0,0) {$ \alpha $};
				\end{tikzpicture}
			\end{array}
		= \phi \circ \varphi (\alpha_G)
		\end{align*}
	where the penultimate equality uses Proposition~\ref{prop:handle_slide} and the final equality uses Proposition~\ref{prop:twisted_S}. \endproof

    \end{PROP}

	\begin{COR} \label{cor:primitive_idems}

	The set $ \{ \e_I^J \}_{I,J \in \Irr(\C)} $ is a set of orthogonal primitive idempotents (see Definition~\ref{def:primitive_idempotents}).

	\proof

	Let $ I,J,I',J' $ be in $ \Irr(\C) $. By Proposition~\ref{prop:homs_between_idems} we have
		\begin{align*}
		\Hom(\e_I^J,\e_{I'}^{J'}) =
			\begin{cases}
			\fld \quad &\text{if $ I = I' $ and $ J = J' $} \\
			0 \quad &\text{else}
			\end{cases}
		\end{align*}
	which proves the claim.

	\endproof

	\end{COR}

As the map from $ \Hom_{\CC}(X \boxtimes Y,A \boxtimes B) = \Hom_{C}(X,A) \otimes \Hom_{C}(Y,B) $ to $ \Hom_{\TC}(\e_X^Y,\e_A^B) $ induced by $ \Phi $ is precisely the map denoted $ \phi $ in the proof of Proposition~\ref{prop:homs_between_idems}, we have already shown that $ \Phi $ is fully faithful. It remains to be shown that $ \Phi $ is essentially surjective, i.e.\ that the set $ \{ \e_I^J \}_{I,J \in \Irr(\C)} $ forms a \emph{complete} set of orthogonal primitive idempotents in $ \TC $. A straightforward consideration of the dimension of $ \Hom_{\TC}(X,Y) $ achieves this.

	\begin{THM} \label{thm:primitive_idempotents}

	Let $ \C $ be an modular tensor category. We have
		\begin{align*}
		\Hom_{\TC}(X,Y) = \bigoplus_{IJ} \Hom_{\TC}(X,\e_I^J) \otimes \Hom_{\TC}(\e_I^J,Y),
		\end{align*}
	in other words, the set $ \{ \e_I^J \}_{I,J \in \Irr(\C)} $ forms a \emph{complete} set of orthogonal primitive idempotents in $ \TC $.

	\proof

	Our aim is to show that the map giving by composition
		\begin{align} \label{eq:composition}
		\bigoplus_{IJ} \Hom_{\TC}(X,\e_I^J) \otimes \Hom_{\TC}(\e_I^J,Y) \to \Hom_{\TC}(X,Y).
		\end{align}
    is an isomorphism. As before, let $ \hom_{\C}(X,Y) $ denote the dimension of $ \Hom_{\C}(X,Y) $. By Corollary~\ref{cor:primitive_idems}, \eqref{eq:composition} is injective and therefore
		\begin{align} \label{eq:dimension_check}
		\sum\limits_{IJ} \hom_{\TC}(X,\e_I^J) \hom_{\TC}(\e_I^J,Y) \leq \hom_{\TC}(X,Y)
		\end{align}
	with equality if and only if \eqref{eq:composition} is an isomorphism. Furthermore, by Proposition~\ref{prop:equivalence_with_f} we have
		\begin{align*}
		\hom_{\C}(X,IJ) = \hom_{\TC}(X,\e^I_J) \quad \text{and} \quad \hom_{\C}(IJ,Y) = \hom_{\TC}(\e^I_J,Y)
		\end{align*}
	which allows us to compute
		\begin{align*}
		\hom_{\TC}(X,Y) &= \sum\limits_{I^\vee} \hom_{C}(I^\vee X,YI^\vee) \\
		&= \sum\limits_{I^\vee,J}  \hom_{C}(I^\vee X, J) \hom_{C}(J,YI^\vee) \\
		&= \sum\limits_{I,J}  \hom_{C}(X,IJ) \hom_{C}(IJ,Y) \\
		&= \sum\limits_{IJ} \hom_{\TC}(X,\e_I^J) \hom_{\TC}(\e_I^J,Y),
		\end{align*}
	implying that \eqref{eq:composition} is an isomorphism.

	\endproof

	\end{THM}

    \begin{COR}

    Let $ \C $ be a modular tensor category. Then $ \Phi\colon \CC \to Z(\C) $ is an equivalence.

    \end{COR}

    \begin{REM}

    As mentioned at the beginning of this section, the converse statement is also true i.e.\ if $ \C $ fails to be modular then $ \Phi $ will also fail to be an equivalence. To prove this we require a converse of the killing ring lemma (Lemma~\ref{lemma:killing_ring}), which is provided by Theorem 8.20.7 in~\cite{Etingof15}. In particular, this theorem implies that if the S-matrix is degenerate, then there exists an object $ I \in \Irr(\C) $ such that
    	\begin{align*}
    		\begin{array}{c}
    			\begin{tikzpicture}[scale = 0.4]
    			\draw[thick] (1,1) to[out=-90,in=90] (-1,-1);
    			\draw[line width = 0.3cm, white] (-1, 1) to[out=-90,in=90] (1,-1);
    			\draw[thick] (1,-1) to[out=-90,in=90] (-1,-3);
    			\draw[thick] (-1, 1) to[out=-90,in=90] (1,-1);
     			\draw[line width = 0.3cm, white] (-1, -1) to[out=-90,in=90] (1,-3);
    			\draw[thick] (-1, -1) to[out=-90,in=90] (1,-3);
    			\node at (-1, 1.7) {$ I $};
    			\node at (1, 1.7) {$ X $};
    			\end{tikzpicture}
    		\end{array}
    	=
    		\begin{array}{c}
    			\begin{tikzpicture}[scale = 0.4]
    			\draw[thick] (1,1) -- (1,-3);
    			\draw[thick] (-1, 1) -- (-1,-3);
    			\node at (-1, 1.7) {$ I $};
    			\node at (1, 1.7) {$ X $};
    			\end{tikzpicture}
    		\end{array}
        \end{align*}
    for all $ X $ in $ \C $. From there one may check that
        \begin{align*}
        \frac{1}{d(\C)} \bigoplus_S d(S)
            \begin{array}{c}
    			\begin{tikzpicture}[scale=0.2,every node/.style={inner sep=0,outer sep=-1}]
    			\node (v1) at (0,5) {};
    			\node (v4) at (1.5,-3.5) {};
    			\node (v2) at (5,0) {};
    			\node (v3) at (-3.5,1.5) {};
    			\node (v5) at (-1.5,3.5) {};
    			\node (v6) at (3.5,-1.5) {};
    			\node (v7) at (1.5,3.5) {};
    			\node (v11) at (3.5,1.5) {};
    			\node (v12) at (1.25,-0.75) {};
    			\node (v10) at (-0.75,1.25) {};
    			\node (v13) at (3.5,-1.5) {};
    			\draw [thick] (v7) edge (v10);
    			\draw [line width =0.5em,white] (v5) edge (v13);
    			\draw [thick] (v5) edge (v13);
    			\draw [line width =0.5em,white] (v11) edge (v12);
    			\draw [thick] (v11) edge (v12);
    			\node at (2.75,4.75) {$ I $};
    			\node at (5,3) {$ I^\vee $};
    			\node at (-2.75,4.75) {$ S $};
    			\node at (4.75,-2.5) {$ S $};
    			\draw[very thick, red]  (v1) edge (v3);
    			\draw[very thick, red]  (v2) edge (v4);
    			\draw[very thick]  (v1) edge (v2);
    			\draw[very thick]  (v3) edge (v4);
    			\draw [thick] (v10) to[out=-135,in=-135] (v12);
    			\end{tikzpicture}
            \end{array}
        \end{align*}
    is a non-zero morphism between $ \e_{I^\vee}^I $ and $ \e_{\tid}^{\tid} $. As $ I \boxtimes I^\vee $ and $ \tid \boxtimes \tid $ are distinct simple objects in $ \CC $, $ \Phi $ is not an equivalence.

    \end{REM}



\bibliographystyle{alpha}
\bibliography{graphical_approach}

\begin{thebibliography}{EGNO15}

\bibitem[BK01]{Bak01}
Bojko Bakalov and Alexander Kirillov, Jr.
\newblock {\em Lectures on tensor categories and modular functors}, volume~21
  of {\em University Lecture Series}.
\newblock American Mathematical Society, Providence, RI, 2001.

\bibitem[EGNO15]{Etingof15}
Pavel Etingof, Shlomo Gelaki, Dmitri Nikshych, and Victor Ostrik.
\newblock {\em Tensor categories}, volume 205 of {\em Mathematical Surveys and
  Monographs}.
\newblock American Mathematical Society, Providence, RI, 2015.

\bibitem[HK19]{hardiman_king}
Leonard Hardiman and Alastair King.
\newblock Decomposing the tube category.
\newblock {\em Glasgow Mathematical Journal}, 2019.

\bibitem[Jr11]{alex2011stringnet}
Alexander~Kirillov Jr.
\newblock String-net model of turaev-viro invariants.
\newblock arXiv:1106.6033v1, 2011.

\bibitem[JS91]{MR1113284}
Andr\'{e} Joyal and Ross Street.
\newblock The geometry of tensor calculus. {I}.
\newblock {\em Adv. Math.}, 88(1):55--112, 1991.

\bibitem[Kas98]{MR1643398}
Ch. Kassel.
\newblock Quantum groups.
\newblock In {\em Algebra and operator theory ({T}ashkent, 1997)}, pages
  213--236. Kluwer Acad. Publ., Dordrecht, 1998.

\bibitem[Kon08]{MR2430629}
Liang Kong.
\newblock Cardy condition for open-closed field algebras.
\newblock {\em Comm. Math. Phys.}, 283(1):25--92, 2008.

\bibitem[Lev06]{MR2717368}
Michael~Aaron Levin.
\newblock {\em String-net condensation and topological phases in quantum spin
  systems}.
\newblock ProQuest LLC, Ann Arbor, MI, 2006.
\newblock Thesis (Ph.D.)--Massachusetts Institute of Technology.

\bibitem[Lur09]{Lurie09}
Jacob Lurie.
\newblock {\em Higher topos theory}, volume 170 of {\em Annals of Mathematics
  Studies}.
\newblock Princeton University Press, Princeton, NJ, 2009.

\bibitem[M{\"u}g03]{MR1990929}
Michael M{\"u}ger.
\newblock On the structure of modular categories.
\newblock {\em Proc. London Math. Soc. (3)}, 87(2):291--308, 2003.

\bibitem[Ocn94]{Ocneanu1993}
Adrian Ocneanu.
\newblock Chirality for operator algebras.
\newblock In {\em Subfactors ({K}yuzeso, 1993)}, pages 39--63. World Sci.
  Publ., River Edge, NJ, 1994.

\bibitem[Pen71]{MR0281657}
Roger Penrose.
\newblock Applications of negative dimensional tensors.
\newblock In {\em Combinatorial {M}athematics and its {A}pplications ({P}roc.
  {C}onf., {O}xford, 1969)}, pages 221--244. Academic Press, London, 1971.

\bibitem[PSV18]{MR3801484}
Sorin Popa, Dimitri Shlyakhtenko, and Stefaan Vaes.
\newblock Cohomology and {$L^2$}-{B}etti numbers for subfactors and
  quasi-regular inclusions.
\newblock {\em Int. Math. Res. Not. IMRN}, (8):2241--2331, 2018.

\bibitem[RT90]{Res90}
N.~Yu. Reshetikhin and V.~G. Turaev.
\newblock Ribbon graphs and their invariants derived from quantum groups.
\newblock {\em Comm. Math. Phys.}, 127(1):1--26, 1990.

\bibitem[ST09]{Snyder09}
Noah Snyder and Peter Tingley.
\newblock The half-twist for {$U_q(\mathfrak{g})$} representations.
\newblock {\em Algebra Number Theory}, 3(7):809--834, 2009.

\end{thebibliography}

\end{document}